**On generalised d'Alembert-type integral representations for damped wave equations on the quarter-plane**

Andreas Chatziafratis [1,*], Claudio Giorgi [2], Alain Miranville [3], and Federico Zullo [2]

[1] Department of Mathematics, National and Kapodistrian University of Athens, Greece
Department of Mathematics and Statistics, University of Cyprus

[2] Sezione di Matematica, Universita degli Studi di Brescia, Italy
Istituto Nazionale di Fisica Nucleare, Milano Bicocca, Italy

[3] Henan Normal University, School of Mathematics and Statistics, Xinxiang, China
Laboratoire de Mathématiques Appliquées du Havre (LMAH), Université Le Havre Normandie, France

**Abstract.** We rigorously construct and verify a posteriori new closed-form solutions for the forced *Maxwell-Cattaneo-Vernotte* equation (also broadly known as the damped wave equation, hyperbolic heat, and *telegrapher*'s equation on lossy transmission lines) posed on the spatiotemporal quarter-plane with general initial and boundary data in classical function spaces. For this purpose, the modern complex-analytic *unified transform method* of Fokas (originally developed for elliptic PDE and evolution equations with *polynomial* dispersion relations) is here, for the first time, extended for analysis of *hyperbolic-parabolic* problems on the semi-infinite interval. Importantly, we then establish theorems which pertain to regularity, boundary and asymptotic properties of the new analytical formulae as well as to well-posedness of the addressed boundary-value problems. Notably, the nature and generality of problems considered, combined with the semi-unboundedness of the domain, induce substantial analytic challenges which demand delicate treatment, both in appropriately interpreting oscillatory integral terms of the solution formulae and in proving the proposed results. In this process, crucially, certain compatibility conditions, between initial, boundary and forcing data at the origin, are revealed, which guarantee the existence of a smooth solution across the whole domain of interest. Our explicit integral representations are of direct utility for numercal benchmarking purposes, for exploring connections with modelling in continuum mechanics, mathematical physics, biology and the natural sciences, and for the investigation of well-posedness for nonlinear counterparts too.

## 1. Introduction

Damped wave equations constitute a fundamental class of evolution problems in partial differential equations, combining finite-speed propagation (causality) with dissipation mechanisms that induce decay of energy. They arise naturally in a wide variety of contexts and provide a canonical framework for the study of hyperbolic systems with relaxation. From a mathematical standpoint, such equations interpolate between conservative wave dynamics and diffusive behavior, and their analysis reveals phenomena absent in purely parabolic models.

In one spatial dimension, the prototypical example is the celebrated equation

$$\alpha u_{tt} + \beta u_t - u_{xx} = f \quad (\alpha, \beta > 0)$$

where $f$ represents an external forcing term. The hyperbolic structure underlying this equation ultimately traces back to d'Alembert's foundational contributions [1] upon which modern relaxation models may be viewed as dissipative perturbations or extensions.

The above equation is of paramount importance as reflected in the vast body of versatile literature it has motivated to date. Importantly, this model is mathematically equivalent, after rescaling, to the telegrapher equation

---

* *corresponding author*: chatziafrati@math.uoa.gr

[31] (developed to model electrical signal propagation along a transmission line with resistive damping) and the Maxwell-Cattaneo-Vernotte equation [4,44,56] (arising as a modification of Fourier-type heat conduction introduced to rectify the infinite propagation speed inherent in the classical heat equation). It has emerged in a plethora of natural and applied sciences, including modelling in electrodynamics, thermodynamics, biological systems and biomedicine, composite materials, mass transport phenomena, special relativity and so on.

More specifically, for instance, the telegraph equation is encountered in the modelling of biological wave phenomena in nerves or in the cardiac system (see e.g. [32,34,38,46,50]), in the description of pressure waves that occur in pipelines in transient flow regimes (for example when the fluid motion is forced to stop suddenly [17,48]), in the description of transport dynamics in magnetized plasma [30,52], in battery technology [42], while it also arises as a limiting case or governing equation in certain stochastic processes, particularly those involving random walks (see e.g. [57]). For additional connections to mechanics, thermodynamics and the applied sciences, the reader is referred to, e.g., [2,18,19,29,37,39,43,47,49,54,55].

Regarding its mathematical treatment when this equation is, in particular, formulated on the spatial half-line, the rigorous analysis of corresponding initial–boundary-value problems (IBVP) in full generality becomes rather unexpectedly delicate. More specifically, the presence of the boundary precludes the direct and effective use of traditional methods thereby necessitating a radically different approach in order to carefully incorporate boundary contributions and capture their interplay with damping and forcing mechanisms. As is well-known by now, conventional techniques can facilitate only partial results, namely they can be applied to certain limited types of problems (see e.g. the classical works [20,45], the comprehensive surveys [36] and the treatise [53] for fundamental solutions, Laplace-transform approaches and other special cases), therefore nonhomogeneous problems for the hyperbolic heat equation on semi-infinite intervals remain relatively incompletely understood.

The objective of this paper is to contribute new results concerning the analytical solution and rigorous qualitative theory of IBVP for the damped wave (or telegraph or hyperbolic heat) equation posed on the spacetime quadrant. To this end, we extend the applicability of the pioneering Fokas' unified transform method (UTM) towards a mixed hyperbolic-parabolic regime. The UTM has been established as powerful and efficient complex-analytic machinery for deriving explicit solution formulae for boundary-value problems in a unified and elegant fashion, with proven substantial advantages over other techniques. This methodology was originally introduced [21,22,28] a few decades ago for tackling elliptic planar problems and evolution PDE with polynomial dispersion relations, while various significant extensions have since been achieved, see e.g. [23-28,33,35,40, 41,51,58] and references therein. The subtle implementation of the UTM in the present work becomes lucid for the reader in view of the illustrative and self-contained demonstration provided in what follows. Subsequent interpretation and rigorous analysis of the resulting integral representations provide insights into various properties of the IBVP considered here; the obtained qualitative findings are herein discussed in detail. Analogous studies in a classically-minded spirit for other evolutionary configurations have been presented in [3,5-8,10-16].

Let us now consider the following:

***Problem*** *Rigorously solve*

$$(1.1)\qquad \begin{cases} \dfrac{\partial^2 u}{\partial t^2}+\dfrac{\partial u}{\partial t}-\dfrac{\partial^2 u}{\partial x^2}=f, \ (x,t)\in Q:=\mathbb{R}^+\times\mathbb{R}^+, \\ \lim\limits_{t\to 0^+} u(x,t)=u_0(x), \ x\in\mathbb{R}^+, \\ \lim\limits_{t\to 0^+} \dfrac{\partial u(x,t)}{\partial t}=u_1(x), \ x\in\mathbb{R}^+, \\ \lim\limits_{x\to 0^+} u(x,t)=g_0(t), \ t\in\mathbb{R}^+, \end{cases}$$

*for* $u=u(x,t)$, and analyse the solution's qualitative behaviour.

Throughtout this paper, we make the following assumptions on the data:

$$(1.2)\qquad u_0(x),u_1(x)\in\mathcal{S}([0,\infty)),\ g_0(t)\in C^\infty([0,\infty)) \text{ and } f=f(x,t)\in C^\infty(\overline{Q}) \text{ such that } f(\cdot,t)\in\mathcal{S}([0,\infty)).$$

More precisely the last assumption on the function $f(x,t)$ means that it is rapidly decreasing with respect to $x$, uniformy for $t$ in compact subsets of $[0,+\infty)$, i.e., for every $\ell,n\in\mathbb{N}\cup\{0\}$ and $t_0>0$,

$$\sup\left\{x^\ell\left|\frac{\partial^n f(x,t)}{\partial x^n}\right| : x\ge 0, 0\le t\le t_0\right\}<+\infty .$$

The function $u(x,t)=e^{i\lambda x-\omega(\lambda)t}$ satisfies the equation in (1.1) if $[\omega(\lambda)]^2-\omega(\lambda)+\lambda^2=0$. The roots of this quadratic equation are

$$\omega=\frac{1}{2}\pm i\sqrt{\lambda^2-\frac{1}{4}}=\frac{1}{2}\pm i\rho(\lambda)\text{, where }\rho(\lambda)=\sqrt{\lambda^2-\frac{1}{4}}.$$

We choose the following holomorphic branch of the square root:

$$\rho(\lambda):=\sqrt{\left|\lambda^2-\frac{1}{4}\right|}\exp\left[i\left(\frac{\theta_1(\lambda)+\theta_2(\lambda)}{2}\right)\right],\ \lambda\in\mathbb{C}-[-\tfrac{1}{2},\tfrac{1}{2}]\text{ (see Fig 1).}$$

We also set

$$\omega_1=\omega_1(\lambda)=\frac{1}{2}+i\rho(\lambda)\quad\text{and}\quad\omega_2=\omega_2(\lambda)=\frac{1}{2}-i\rho(\lambda).$$

Then the functions $\rho(\lambda)$, $\omega_1(\lambda)$, $\omega_2(\lambda)$ are analytic for $\lambda\in\mathbb{C}-[-\tfrac{1}{2},\tfrac{1}{2}]$.

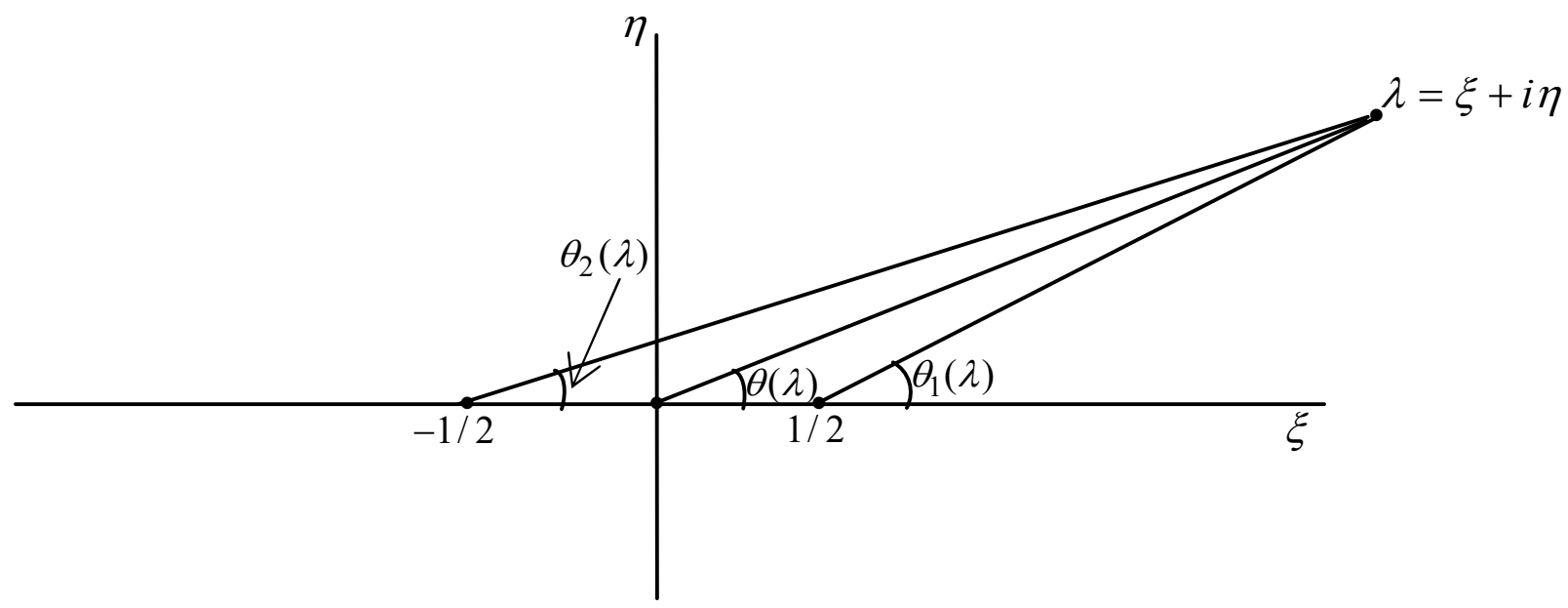

**Fig. 1** $\lambda=|\lambda|e^{i\theta(\lambda)}$, $\lambda-(1/2)=|\lambda-(1/2)|e^{i\theta_1(\lambda)}$, $\lambda-(-1/2)=|\lambda-(-1/2)|e^{i\theta_2(\lambda)}$, $-\pi<\theta(\lambda),\theta_1(\lambda),\theta_2(\lambda)\le\pi$.

We extend $\rho(\lambda)$ also for $\lambda\in[-\tfrac{1}{2},\tfrac{1}{2}]$ as follows:

$$\rho(\lambda):=\lim_{\substack{\mu\to\lambda\\ \operatorname{Im}\mu<0}}\rho(\mu)=-i\sqrt{\left|\lambda^2-\frac{1}{4}\right|},\ -\tfrac{1}{2}\le\lambda\le\tfrac{1}{2}.$$

Then $[\rho(\lambda)]^2=\lambda^2-(1/4)$, for every $\lambda\in\mathbb{C}$, and the restriction of $\rho(\lambda)$ to the (closed) lower half-plane $\{\lambda\in\mathbb{C}:\operatorname{Im}\lambda\le0\}$ is continuous – in fact it is $C^\infty$ in $\{\lambda\in\mathbb{C}:\operatorname{Im}\lambda\le0\}-\{-\tfrac{1}{2},\tfrac{1}{2}\}$.

Thus, the functions $\rho(\lambda)$, $\omega_1(\lambda)$, $\omega_2(\lambda)$ are defined for all $\lambda\in\mathbb{C}$.

We point out that $\omega_1(\lambda)-\omega_2(\lambda)=2i\rho(\lambda)$, for $\lambda\in\mathbb{C}$,

(1.3) $$\rho(-\lambda)=-\rho(\lambda),\ \omega_1(-\lambda)=\omega_2(\lambda),\ \omega_2(-\lambda)=\omega_1(\lambda)\text{ for }\lambda\in\mathbb{C}-[-\tfrac{1}{2},\tfrac{1}{2}]$$

and

(1.4) $$\rho(-\lambda)=\rho(\lambda),\ \omega_1(-\lambda)=\omega_1(\lambda),\ \omega_2(-\lambda)=\omega_2(\lambda)\text{ for }-\tfrac{1}{2}\le\lambda\le\tfrac{1}{2}.$$

We will use also the following notation:

$$\hat{u}_0(\lambda):=\int_0^\infty e^{-i\lambda x}u_0(x)dx\text{ and }\hat{u}_1(\lambda):=\int_0^\infty e^{-i\lambda x}u_1(x)dx,\ \lambda\in\mathbb{C}\text{ with }\operatorname{Im}\lambda\le0$$

(we recall that the restriction of $\lambda$ to $\operatorname{Im}\lambda\le0$ is *generically* necessary – see [9]),

$$\tilde{g}_0(\omega_j(\lambda),t):=\int_{\tau=0}^t e^{\omega_j(\lambda)t}g_0(\tau)d\tau,\ \lambda\in\mathbb{C},\ j=1,2,$$

$$\tilde{\hat{f}}(\lambda,\omega_j(\lambda),t):=\int_{\tau=0}^t e^{\omega_j(\lambda)\tau}\hat{f}(\lambda,\tau)d\tau\ (j=1,2),\text{ where }\hat{f}(\lambda,\tau):=\int_0^\infty e^{-i\lambda x}f(x,\tau)dx,\ \lambda\in\mathbb{C}\text{ with }\operatorname{Im}\lambda\le0.$$

With this notation the Fokas solution of Problem (1.1) is given by the following formula: For $x>0$, $t>0$, $x\ne t$,

(1.5) $$2\pi u(x,t)=\int_{-\infty}^{\infty}[\omega_1(\lambda)e^{i\lambda x-\omega_2(\lambda)t}-\omega_2(\lambda)e^{i\lambda x-\omega_1(\lambda)t}]\hat{u}_0(\lambda)\frac{d\lambda}{2i\rho(\lambda)}$$
$$+\int_{-\infty}^{\infty}[\omega_2(\lambda)e^{i\lambda x-\omega_1(\lambda)t}-\omega_1(\lambda)e^{i\lambda x-\omega_2(\lambda)t}]\hat{u}_0(-\lambda)\frac{d\lambda}{2i\rho(\lambda)}$$
$$+\int_{-\infty}^{\infty}[e^{i\lambda x-\omega_2(\lambda)t}-e^{i\lambda x-\omega_1(\lambda)t}]\hat{u}_1(\lambda)\frac{d\lambda}{2i\rho(\lambda)}+\int_{-\infty}^{\infty}[e^{i\lambda x-\omega_1(\lambda)t}-e^{i\lambda x-\omega_2(\lambda)t}]\hat{u}_1(-\lambda)\frac{d\lambda}{2i\rho(\lambda)}$$
$$+\int_{-\infty}^{\infty}[e^{i\lambda x-\omega_1(\lambda)t}\tilde{g}_0(\omega_1(\lambda),t)-e^{i\lambda x-\omega_2(\lambda)t}\tilde{g}_0(\omega_2(\lambda),t)]\frac{\lambda d\lambda}{\rho(\lambda)}$$
$$+\int_{-\infty}^{\infty}\left[e^{i\lambda x-\omega_2(\lambda)t}\tilde{\hat{f}}(\lambda,\omega_2(\lambda),t)-e^{i\lambda x-\omega_1(\lambda)t}\tilde{\hat{f}}(\lambda,\omega_1(\lambda),t)\right]\frac{d\lambda}{2i\rho(\lambda)}$$
$$+\int_{-\infty}^{\infty}\left[e^{i\lambda x-\omega_1(\lambda)t}\tilde{\hat{f}}(-\lambda,\omega_1(\lambda),t)-e^{i\lambda x-\omega_2(\lambda)t}\tilde{\hat{f}}(-\lambda,\omega_2(\lambda),t)\right]\frac{d\lambda}{2i\rho(\lambda)}.$$

Equivalently,

(1.6) $$2\pi u(x,t)=[\mathcal{U}_0^+(u_0)+\mathcal{U}_0^-(u_0)+\mathcal{U}_1^+(u_1)+\mathcal{U}_1^-(u_1)+\mathcal{G}(g_0)+\mathcal{F}^+(f)+\mathcal{F}^-(f)](x,t),$$

where

$$\mathcal{U}_0^+(u_0),\ \mathcal{U}_0^-(u_0),\ \mathcal{U}_1^+(u_1),\ \mathcal{U}_1^-(u_1),\ \mathcal{G}(g_0),\ \mathcal{F}^+(f),\ \mathcal{F}^-(f)$$

are the functions of $(x,t)$, defined by the integrals in (1.5) which involve the terms

$$\hat{u}_0(\lambda),\ \hat{u}_0(-\lambda),\ \hat{u}_1(\lambda),\ \hat{u}_1(-\lambda),\ \tilde{g}_0(\lambda,\cdot),\ \tilde{\hat{f}}(\lambda,\cdot,t),\ \tilde{\hat{f}}(-\lambda,\cdot,t),$$

respectively.

First we will prove the following main theorems for the function $u(x,t)$, defined for $(x,t)\in Q-\{x=t\}$ by (1.5). (More theorems on qualitative theory are formulated and established in subsequent sections.) We recall our ***main*** assumptions (1.2), which we make throughtout this paper.

**Theorem 1** *For $(x,t)\in Q$ with $x\neq t$, the integrals in (1.5) exist in the generalized sense, i.e, as limits of the form $\int_{-\infty}^{\infty}=\lim_{A\to\infty}\int_{-A}^{A}$, the function $u(x,t)$, defined this way, is $C^\infty$ for $(x,t)$ in the sets $\{x<t\}=Q-\{x\geq t\}$ and $\{x>t\}=Q-\{x\leq t\}$, and satisfies the differential equation $u_{tt}+u_t-u_{xx}=f$ in $Q-\{x=t\}$.*
*Moreover, the restriction $u(x,t)|_{\{x<t\}}$ extends to a $C^\infty$ function in $Q\cap\{x\leq t\}$ and the restriction $u(x,t)|_{\{x>t\}}$ extends to a $C^\infty$ function in $Q\cap\{x\geq t\}$, and, in particular, $u(x,t)$ is locally integrable in $Q$.*

**Theorem 2** *1st $\lim_{t\to0^+}u(x,t)=u_0(x)$ and $\lim_{t\to0^+}\frac{\partial u(x,t)}{\partial t}=u_1(x)$, uniformly for $x$ in compact sets of $\mathbb{R}^+$.*
*2nd $\lim_{x\to0^+}u(x,t)=g_0(t)$, uniformly for $t$ in compact sets of $\mathbb{R}^+$.*
*3rd If we assume, in addition to (1.2), that $u_0(0)=g_0(0)$ then $u(x,t)$ extends to a continuous function in all of $Q$, which is given by the following formula: For $(x,t)\in Q$,*

*(1.7)* $$2\pi u(x,t)=\int_{-1}^{1}[\omega_1(\lambda)e^{i\lambda x-\omega_2(\lambda)t}-\omega_2(\lambda)e^{i\lambda x-\omega_1(\lambda)t}]\hat{u}_0(\lambda)\frac{d\lambda}{2i\rho(\lambda)}$$
$$+\left(\int_{-\infty}^{-1}+\int_{1}^{\infty}\right)[\omega_1(\lambda)e^{i\lambda x-\omega_2(\lambda)t}-\omega_2(\lambda)e^{i\lambda x-\omega_1(\lambda)t}]\frac{(u_0^{(1)})\hat{}(\lambda)}{i\lambda}\frac{d\lambda}{2i\rho(\lambda)}$$
$$+\int_{-1}^{1}[\omega_2(\lambda)e^{i\lambda x-\omega_1(\lambda)t}-\omega_1(\lambda)e^{i\lambda x-\omega_2(\lambda)t}]\hat{u}_0(-\lambda)\frac{d\lambda}{2i\rho(\lambda)}$$

$$+\left(\int_{-\infty}^{-1}+\int_{1}^{\infty}\right)[\omega_2(\lambda)e^{i\lambda x-\omega_1(\lambda)t}-\omega_1(\lambda)e^{i\lambda x-\omega_2(\lambda)t}]\frac{(u_0^{(1)})\hat{}(-\lambda)}{-i\lambda}\frac{d\lambda}{2i\rho(\lambda)}$$

$$+\int_{-1}^{1}[e^{i\lambda x-\omega_1(\lambda)t}\tilde{g}_0(\omega_1(\lambda),t)-e^{i\lambda x-\omega_2(\lambda)t}\tilde{g}_0(\omega_2(\lambda),t)]\frac{\lambda d\lambda}{\rho(\lambda)}+g_0(t)\int_{\kappa^+}e^{i\lambda x}\left(\frac{1}{\omega_1(\lambda)}-\frac{1}{\omega_2(\lambda)}\right)d\lambda$$

$$-\left(\int_{-\infty}^{-1}+\int_{1}^{\infty}\right)\left[e^{i\lambda x-\omega_1(\lambda)t}\frac{1}{\omega_1(\lambda)}(g_0^{(1)})\tilde{}(\omega_1(\lambda),t)-e^{i\lambda x-\omega_2(\lambda)t}\frac{1}{\omega_2(\lambda)}(g_0^{(1)})\tilde{}(\omega_2(\lambda),t)\right]\frac{\lambda d\lambda}{\rho(\lambda)}$$

$$+[\mathcal{U}_1^+(u_1)+\mathcal{U}_1^-(u_1)+\mathcal{F}^+(f)+\mathcal{F}^-(f)](x,t),$$

*where all the above integrals converge absolutely. ($\kappa^+$ is the semicircle in the upper half-plane which starts at $\lambda=+1$ and ends at the point $\lambda=-1$. It is depicted in Fig. 2.)*

$4^{th}$ *If* $u_0(0)=g_0(0)$ *and* $u_1(0)=g_0'(0)$ *then* $u(x,t)$*, given by (1.7), is* $C^1$ *for* $(x,t)\in Q$.

$5^{th}$ *If we assume that*

*(1.8)* $$u_0(0)=g_0(0),\ u_1(0)=g_0'(0)\ \text{and}\ g_0''(0)+g_0'(0)-u_0''(0)=f(0,0),$$

*then* $u(x,t)$*, given by (1.7), is* $C^2$ *for* $(x,t)\in Q$ *and satisfies the equation* $u_{tt}+u_t-u_{xx}=f$ *in* $Q$.

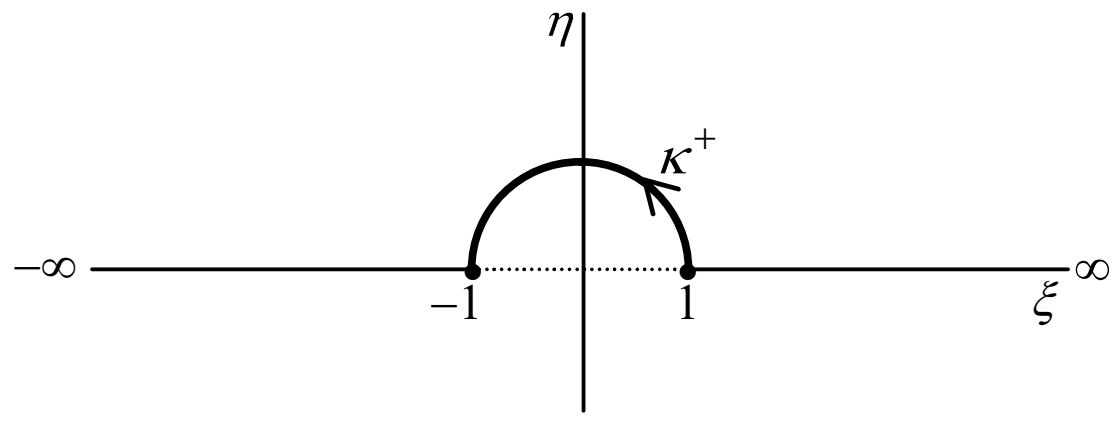

**Fig. 2** The semicircle $\kappa^+$

**Theorem 3** $1^{st}$ *For each fixed* $t>0$*, the limits*

$$g_{n,m}(t):=\lim_{x\to 0^+}\frac{\partial^{n+m}u(x,t)}{\partial x^n\partial t^m},\ n,m=0,1,2,...,$$

*exist. Moreover the functions* $g_{n,m}(t)$ *are* $C^\infty$ *for* $t\in(0,\infty)$. *Moreover, the above convergence is uniform for* $t$ *in compact subsets of* $(0,\infty)$.

$2^{nd}$ *For each fixed* $x>0$*, the limits*

$$u_{n,m}(x):=\lim_{t\to 0^+}\frac{\partial^{n+m}u(x,t)}{\partial x^n\partial t^m},\ n,m=0,1,2,...,$$

*exist. Moreover the functions* $u_{n,m}(x)$ *are* $C^\infty$ *for* $x\in(0,\infty)$. *Moreover, the above convergence is uniform for* $x$ *in compact subsets of* $(0,\infty)$.

$3^{rd}$ *The function* $u(x,t)$*, extended to* $\overline{Q}-\{(p,p):p\geq 0\}$ *by setting* $u(x,0):=u_0(x)$ *for* $x>0$ *and* $u(0,t):=g_0(t)$ *for* $t>0$*, is* $C^\infty$ *(in* $\overline{Q}-\{(p,p):p\geq 0\}$*).*

$4^{th}$ *For* $n=0,1,2,...$, and $x>0$,

$$\lim_{t\to 0^+}\frac{\partial^n u(x,t)}{\partial x^n}=\frac{d^n u_0(x)}{dx^n}\ \text{and}\ \lim_{t\to 0^+}\frac{\partial^{n+1}u(x,t)}{\partial x^n\partial t}=\frac{d^n u_1(x)}{dx^n},$$

*with the convergence being uniform for* $x$ *in compact subsets of* $(0,\infty)$.

$5^{th}$ *For* $m=0,1,2,...$, and $t>0$,

$$\lim_{x\to 0^+}\frac{\partial^m u(x,t)}{\partial t^m}=\frac{d^m g_0(t)}{dt^m}$$

*with the convergence being uniform for* $t$ *in compact subsets of* $(0,\infty)$.

$6^{th}$ *Under the further assumptions (1.8), the limits*

$$\lim_{\substack{(x,t)\to(0,0)\\(x,y)\in\overline{Q}-\{(0,0)\}}} \frac{\partial^{n+m}u(x,t)}{\partial x^n \partial t^m} \ \textit{exist, for } n+m\le 2 .$$

**Theorem 4** *The function* $u(x,t)$ *is rapidly decreasing with respect to* $x$. *More precisely,*

$$\lim_{x\to+\infty}\left[x^{\ell}\frac{\partial^n u(x,t)}{\partial x^n}\right]=0$$

*for nonnegative integers* $n$ *and* $\ell$, *uniformly for* $t$ *in compact subsets of* $[0,\infty)$.

Further results about asymptotic properties and well-posedness are presented in later sections of the manuscript.

## 2. Derivation of a solution formula

***Step 1*** The differential equation in (1.1) can be written in divergence form as follows:

$$(2.1)\qquad \frac{\partial}{\partial t}\left\{\left([1-\omega(\lambda)]u+u_t\right)e^{-i\lambda x+\omega(\lambda)t}\right\}-\frac{\partial}{\partial x}\left[(i\lambda u+u_x)e^{-i\lambda x+\omega(\lambda)t}\right]=f(x,t)e^{-i\lambda x+\omega(\lambda)t} .$$

Using Green's theorem, (2.1) gives

$$(2.2)\quad \int_{x=0}^{\infty}\left\{\left([1-\omega(\lambda)]u+u_t\right)e^{-i\lambda x+\omega(\lambda)t}\right\}\Big|_{t=0}dx-\int_{x=0}^{\infty}\left\{\left([1-\omega(\lambda)]u+u_t\right)e^{-i\lambda x+\omega(\lambda)\tau}\right\}\Big|_{\tau=t}dx$$

$$-\int_{\tau=0}^{t}\left[(i\lambda u+u_x)e^{-i\lambda x+\omega(\lambda)\tau}\right]\Big|_{x=0}d\tau=-\int_{\tau=0}^{t}e^{\omega(\lambda)\tau}\left[\int_{x=0}^{\infty}f(x,\tau)e^{-i\lambda x}dx\right]d\tau ,\quad \operatorname{Im}\lambda\le 0 .$$

Indeed, this follows by fixing $\mathrm{A}>0$ and $\mathrm{T}>0$, applying Green's formula in the set

$$\Pi_{\mathrm{A,T}}=\{(x,t)\in\mathbb{R}^2 : 0\le x\le \mathrm{A} \ \textit{and}\ 0\le t\le \mathrm{T}\}$$

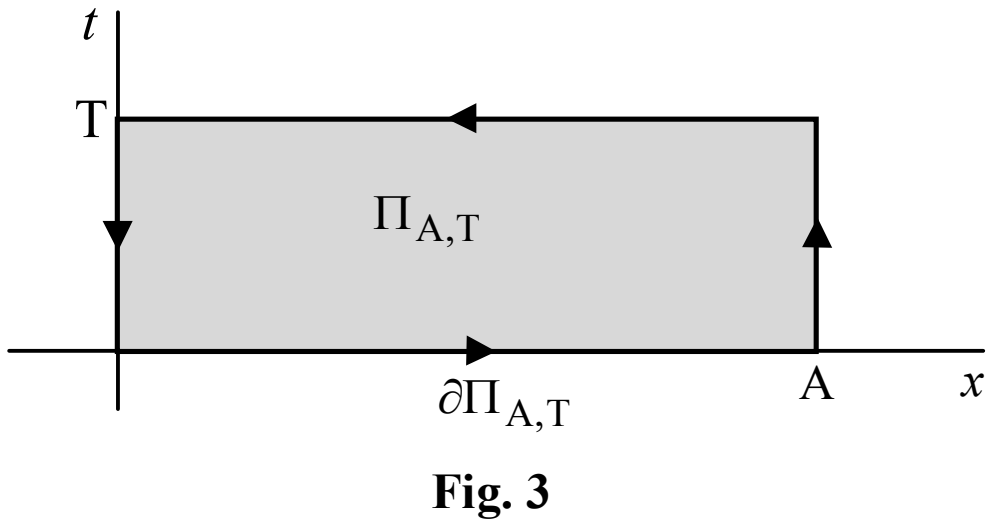

**Fig. 3**

and letting $\mathrm{A}\to\infty$. At this stage we have to assume that the behavior of the function $u(x,t)$ and its derivatives $u_x(x,t)$ and $u_{xx}(x,t)$, as $x\to\infty$, is such that

$$\lim_{\mathrm{A}\to\infty}\int_{t=0}^{\mathrm{T}}e^{-i\lambda \mathrm{A}+\omega(\lambda)t}[u_{xx}(\mathrm{A},t)+i\lambda u_x(\mathrm{A},t)+(i\lambda)^2u(\mathrm{A},t)]dt=0 .$$

Then (2.2) leads to the equation

$$\{[1-\omega(\lambda)]\hat{u}_0(\lambda)+\hat{u}_1(\lambda)\}-\{[1-\omega(\lambda)]\hat{u}(\lambda,t)e^{\omega(\lambda)t}+\hat{(u_t)}(\lambda,t)e^{\omega(\lambda)t}\}$$

$$=(i\lambda)\tilde{g}_0(\omega(\lambda),t)-\tilde{g}_1(\omega(\lambda),t)-\tilde{\hat{f}}(\lambda,\omega(\lambda),t) .$$

Multiplying by $e^{-\omega(\lambda)t}$, we have

$$(2.3)\quad \{[1-\omega(\lambda)]e^{-\omega(\lambda)t}\hat{u}_0(\lambda)+e^{-\omega(\lambda)t}\hat{u}_1(\lambda)\}-\{[1-\omega(\lambda)]\hat{u}(\lambda,t)+\hat{(u_t)}(\lambda,t)\}$$

$$=(i\lambda)e^{-\omega(\lambda)t}\tilde{g}_0(\omega(\lambda),t)-e^{-\omega(\lambda)t}\tilde{g}_1(\omega(\lambda),t)-e^{-\omega(\lambda)t}\tilde{\hat{f}}(\lambda,\omega(\lambda),t) .$$

Writing (2.3) for each $\omega$ with $\omega\in\{\omega_1,\omega_2\}$ and taking into consideration that $\omega_1(\lambda)+\omega_2(\lambda)=1$, we obtain

$$(2.4)\quad [\omega_2e^{-\omega_1t}\hat{u}_0(\lambda)+e^{-\omega_1t}\hat{u}_1(\lambda)]-[\omega_2\hat{u}(\lambda,t)+\hat{(u_t)}(\lambda,t)]=(i\lambda)e^{-\omega_1t}\tilde{g}_0(\omega_1,t)-e^{-\omega_1t}\tilde{g}_1(\omega_1,t)-e^{-\omega_1t}\tilde{\hat{f}}(\lambda,\omega_1,t)$$

and

(2.5) $\quad [\omega_1 e^{-\omega_2 t}\hat{u}_0(\lambda) + e^{-\omega_2 t}\hat{u}_1(\lambda)] - [\omega_1 \hat{u}(\lambda,t) + \hat{(u_t)}(\lambda,t)] = (i\lambda)e^{-\omega_2 t}\tilde{g}_0(\omega_2,t) - e^{-\omega_2 t}\tilde{g}_1(\omega_2,t) - e^{-\omega_2 t}\tilde{\hat{f}}(\lambda,\omega_2,t)$.

Subtracting (2.5) from (2.4), we have that, for $\operatorname{Im}\lambda \le 0$,

$$(\omega_2 e^{-\omega_1 t} - \omega_1 e^{-\omega_2 t})\hat{u}_0(\lambda) + (e^{-\omega_1 t} - e^{-\omega_2 t})\hat{u}_1(\lambda) + (\omega_1 - \omega_2)\hat{u}(\lambda,t)$$

$$= (i\lambda)e^{-\omega_1 t}\tilde{g}_0(\omega_1,t) - (i\lambda)e^{-\omega_2 t}\tilde{g}_0(\omega_2,t) - e^{-\omega_1 t}\tilde{g}_1(\omega_1,t) + e^{-\omega_2 t}\tilde{g}_1(\omega_2,t) - e^{-\omega_1 t}\tilde{\hat{f}}(\lambda,\omega_1,t) + e^{-\omega_2 t}\tilde{\hat{f}}(\lambda,\omega_2,t),$$

and, therefore,

(2.6)
$$\frac{1}{\omega_1-\omega_2}(\omega_2 e^{-\omega_1 t} - \omega_1 e^{-\omega_2 t})\hat{u}_0(\lambda) + \frac{1}{\omega_1-\omega_2}(e^{-\omega_1 t} - e^{-\omega_2 t})\hat{u}_1(\lambda) + \hat{u}(\lambda,t)$$
$$= \frac{i\lambda}{\omega_1-\omega_2}e^{-\omega_1 t}\tilde{g}_0(\omega_1,t) - \frac{i\lambda}{\omega_1-\omega_2}e^{-\omega_2 t}\tilde{g}_0(\omega_2,t) - \frac{1}{\omega_1-\omega_2}e^{-\omega_1 t}\tilde{g}_1(\omega_1,t) + \frac{1}{\omega_1-\omega_2}e^{-\omega_2 t}\tilde{g}_1(\omega_2,t)$$
$$- \frac{1}{\omega_1-\omega_2}e^{-\omega_1 t}\tilde{\hat{f}}(\lambda,\omega_1,t) + \frac{1}{\omega_1-\omega_2}e^{-\omega_2 t}\tilde{\hat{f}}(\lambda,\omega_2,t), \text{ for } \operatorname{Im}\lambda \le 0,\ \lambda \ne \pm\tfrac{1}{2}.$$

Let us keep in mind that $\omega_1(\lambda) - \omega_2(\lambda) = 2i\rho(\lambda)$ and $\omega_1(\lambda) - \omega_2(\lambda) \ne 0$ if and only if $\lambda \ne \pm\frac{1}{2}$.

***Step 2*** The following lemma is crucial for the construction of the kernels of the solution. According to this lemma, certain quantities, which appear in the integrals of (1.5), although they contain the function $\rho(\lambda)$, are, at least, $C^\infty$ for $\lambda \in (-\infty,\infty)$. This means, in particular, that the points $\lambda = \pm 1/2$ are, in a sense, "removable" singularities of certain combinations of functions which appear in the kernels, despite the fact that

$$[2i\rho(\lambda)]\big|_{\lambda=\pm 1/2} = [\omega_1(\lambda) - \omega_2(\lambda)]\big|_{\lambda=\pm 1/2} = 0.$$

***Lemma 1*** *For fixed $t$, the functions*

$$\mathcal{K}_0(\lambda,t) := \frac{\omega_1(\lambda)e^{-\omega_2(\lambda)t} - \omega_2(\lambda)e^{-\omega_1(\lambda)t}}{2i\rho(\lambda)},$$

$$\mathcal{K}_1(\lambda,t) := \frac{e^{-\omega_2(\lambda)t} - e^{-\omega_1(\lambda)t}}{2i\rho(\lambda)},$$

$$\mathcal{L}_{g_0}(\lambda,t) := \frac{e^{-\omega_1(\lambda)t}\tilde{g}_0(\omega_1(\lambda),t) - e^{-\omega_2(\lambda)t}\tilde{g}_0(\omega_2(\lambda),t)}{\rho(\lambda)}\lambda,$$

$$\mathcal{L}_{g_1}(\lambda,t) := \frac{e^{-\omega_1(\lambda)t}\tilde{g}_1(\omega_1(\lambda),t) - e^{-\omega_2(\lambda)t}\tilde{g}_1(\omega_2(\lambda),t)}{\rho(\lambda)},$$

$$\mathcal{M}_f^+(\lambda,t) := \frac{e^{-\omega_2(\lambda)t}\tilde{\hat{f}}(\lambda,\omega_2(\lambda),t) - e^{-\omega_1(\lambda)t}\tilde{\hat{f}}(\lambda,\omega_1(\lambda),t)}{\rho(\lambda)},$$

$$\mathcal{M}_f^-(\lambda,t) := \frac{e^{-\omega_1(\lambda)t}\tilde{\hat{f}}(-\lambda,\omega_1(\lambda),t) - e^{-\omega_2(\lambda)t}\tilde{\hat{f}}(-\lambda,\omega_2(\lambda),t)}{\rho(\lambda)}$$

*are $C^\infty$ for $\lambda \in (-\infty,\infty)$. In fact the functions $\mathcal{K}_0(\lambda,t)$, $\mathcal{K}_1(\lambda,t)$, $\mathcal{L}_{g_0}(\lambda,t)$, $\mathcal{L}_{g_1}(\lambda,t)$ are analytic for $\lambda \in \mathbb{C}$, the function $\mathcal{M}_f^+(\lambda,t)$ is analytic in $\{\lambda \in \mathbb{C} : \operatorname{Im}\lambda < 0\}$ and extends to be $C^\infty$ in $\{\lambda \in \mathbb{C} : \operatorname{Im}\lambda \le 0\}$, while the function $\mathcal{M}_f^-(\lambda,t)$ is analytic in $\{\lambda \in \mathbb{C} : \operatorname{Im}\lambda > 0\}$ and extends to be $C^\infty$ in $\{\lambda \in \mathbb{C} : \operatorname{Im}\lambda \ge 0\}$.*
*Furthermore, for $\lambda \in \mathbb{R}$,*

(2.7) $$\mathcal{K}_0(-\lambda,t) = \mathcal{K}_0(\lambda,t),\ \mathcal{K}_1(-\lambda,t) = \mathcal{K}_1(\lambda,t),\ \mathcal{L}_{g_0}(-\lambda,t) = \mathcal{L}_{g_0}(\lambda,t),\ \mathcal{L}_{g_1}(-\lambda,t) = \mathcal{L}_{g_1}(\lambda,t)$$

*and*

(2.8) $$\mathcal{M}_f^+(-\lambda,t) = \mathcal{M}_f^-(\lambda,t).$$

***Proof*** We have

$$(2.9)\quad \mathcal{K}_1(\lambda,t)=\frac{e^{-[\frac{1}{2}-i\rho(\lambda)]t}-e^{-[\frac{1}{2}+i\rho(\lambda)]t}}{2i\rho(\lambda)}=e^{-t/2}\frac{e^{i\rho(\lambda)t}-e^{-i\rho(\lambda)t}}{2i\rho(\lambda)}=e^{-t/2}\frac{\sin[\rho(\lambda)t]}{\rho(\lambda)}$$

$$=e^{-t/2}\frac{1}{\rho(\lambda)}\sum_{n=1}^{\infty}\frac{(-1)^{n-1}}{(2n-1)!}[t\rho(\lambda)]^{2n-1}=e^{-t/2}\sum_{n=1}^{\infty}\frac{1}{(2n-1)!}\left(\frac{1}{4}-\lambda^2\right)^{n-1}t^{2n-1},$$

and this proves the analyticity of $\mathcal{K}_1(\lambda,t)$ in $\lambda$.

Similarly,

$$(2.10)\qquad \mathcal{K}_0(\lambda,t)=\frac{[\frac{1}{2}+i\rho(\lambda)]e^{-[\frac{1}{2}-i\rho(\lambda)]t}-[\frac{1}{2}-i\rho(\lambda)]e^{-[\frac{1}{2}+i\rho(\lambda)]t}}{2i\rho(\lambda)}$$

$$=\frac{1}{2}e^{-t/2}\frac{e^{i\rho(\lambda)t}-e^{-i\rho(\lambda)t}}{2i\rho(\lambda)}+\frac{1}{2}e^{-t/2}\left[e^{i\rho(\lambda)t}+e^{-i\rho(\lambda)t}\right]$$

$$=\frac{1}{2}e^{-t/2}\frac{\sin[\rho(\lambda)t]}{\rho(\lambda)}+e^{-t/2}\cos[\rho(\lambda)t]$$

$$=\frac{1}{2}e^{-t/2}\sum_{n=1}^{\infty}\frac{1}{(2n-1)!}\left(\frac{1}{4}-\lambda^2\right)^{n-1}t^{2n-1}+e^{-t/2}\sum_{n=1}^{\infty}\frac{1}{(2n-2)!}\left(\frac{1}{4}-\lambda^2\right)^{n-1}t^{2n-2},$$

and this proves the analyticity of $\mathcal{K}_0(\lambda,t)$ in $\lambda$.

For the function $\mathcal{L}_{g_0}(\lambda,t)$, in view of (2.9), we have

$$(2.11)\qquad \mathcal{L}_{g_0}(\lambda,t)=-2i\int_{\tau=0}^{t}\left[e^{-(t-\tau)/2}\sum_{n=1}^{\infty}\frac{1}{(2n-1)!}\left(\frac{1}{4}-\lambda^2\right)^{n-1}\lambda(t-\tau)^{2n-1}\right]g_0(\tau)d\tau,$$

and its analyticity follows. Similar is the case of $\mathcal{L}_{g_1}(\lambda,t)$.

Finally, we find

$$(2.12)\qquad \mathcal{M}_f^+(\lambda,t)=\int_{\tau=0}^{t}\left[e^{-(t-\tau)/2}\sum_{n=1}^{\infty}\frac{1}{(2n-1)!}\left(\frac{1}{4}-\lambda^2\right)^{n-1}(t-\tau)^{2n-1}\right]\hat{f}(\lambda,\tau)d\tau$$

and

$$(2.13)\qquad \mathcal{M}_f^-(\lambda,t)=-\int_{\tau=0}^{t}\left[e^{-(t-\tau)/2}\sum_{n=1}^{\infty}\frac{1}{(2n-1)!}\left(\frac{1}{4}-\lambda^2\right)^{n-1}(t-\tau)^{2n-1}\right]\hat{f}(-\lambda,\tau)d\tau.$$

the corresponding assertions follow.

***Step 3*** Multiplying (2.6) by $e^{i\lambda x}$ and integrating, we obtain

$$(2.14)\quad \int_{-\infty}^{\infty}\frac{[\omega_2(\lambda)e^{i\lambda x-\omega_1(\lambda)t}-\omega_1(\lambda)e^{i\lambda x-\omega_2(\lambda)t}]}{\omega_1(\lambda)-\omega_2(\lambda)}\hat{u}_0(\lambda)d\lambda+\int_{-\infty}^{\infty}\frac{e^{i\lambda x-\omega_1(\lambda)t}-e^{i\lambda x-\omega_2(\lambda)t}}{\omega_1(\lambda)-\omega_2(\lambda)}\hat{u}_1(\lambda)d\lambda+\int_{-\infty}^{\infty}e^{i\lambda x}\hat{u}(\lambda,t)d\lambda$$

$$=\int_{-\infty}^{\infty}\left\{\frac{i\lambda e^{i\lambda x-\omega_1(\lambda)t}}{\omega_1(\lambda)-\omega_2(\lambda)}\widetilde{g}_0(\omega_1(\lambda),t)-\frac{i\lambda e^{i\lambda x-\omega_2(\lambda)t}}{\omega_1(\lambda)-\omega_2(\lambda)}\widetilde{g}_0(\omega_2(\lambda),t)\right\}d\lambda$$

$$-\int_{-\infty}^{\infty}\left\{\frac{e^{i\lambda x-\omega_1(\lambda)t}}{\omega_1(\lambda)-\omega_2(\lambda)}\widetilde{g}_1(\omega_1(\lambda),t)-\frac{e^{i\lambda x-\omega_2(\lambda)t}}{\omega_1(\lambda)-\omega_2(\lambda)}\widetilde{g}_1(\omega_2(\lambda),t)\right\}d\lambda$$

$$-\int_{-\infty}^{\infty}\left\{\frac{e^{i\lambda x-\omega_1(\lambda)t}}{\omega_1(\lambda)-\omega_2(\lambda)}\widetilde{\hat{f}}(\lambda,\omega_1(\lambda),t)-\frac{e^{i\lambda x-\omega_2(\lambda)t}}{\omega_1(\lambda)-\omega_2(\lambda)}\widetilde{\hat{f}}(\lambda,\omega_2(\lambda),t)\right\}d\lambda,\ \text{for } \operatorname{Im}\lambda\le 0.$$

***Step 4*** For $\lambda\in\mathbb{R}$, we have $\operatorname{Im}(-\lambda)\le 0$ and, therefore, (2.6) gives

$$(2.15)\quad \frac{1}{\omega_1(-\lambda)-\omega_2(-\lambda)}[\omega_2(-\lambda)e^{-\omega_1(-\lambda)t}-\omega_1(-\lambda)]e^{-\omega_2(-\lambda)t}]\hat{u}_0(-\lambda)$$

$$+\frac{1}{\omega_1(-\lambda)-\omega_2(-\lambda)}(e^{-\omega_1(-\lambda)t}-e^{-\omega_2(-\lambda)t})\hat{u}_1(-\lambda)+\hat{u}(-\lambda,t)$$

$$=\frac{-i\lambda}{\omega_1(-\lambda)-\omega_2(-\lambda)}e^{-\omega_1(-\lambda)t}\widetilde{g}_0(\omega_1(-\lambda),t)-\frac{-i\lambda}{\omega_1(-\lambda)-\omega_2(-\lambda)}e^{-\omega_2(-\lambda)t}\widetilde{g}_0(\omega_2(-\lambda),t)$$

$$-\frac{1}{\omega_1(-\lambda)-\omega_2(-\lambda)}e^{-\omega_1(-\lambda)t}\widetilde{g}_1(\omega_1(-\lambda),t)+\frac{1}{\omega_1(-\lambda)-\omega_2(-\lambda)}e^{-\omega_2(-\lambda)t}\widetilde{g}_1(\omega_2(-\lambda),t)$$

$$-\frac{1}{\omega_1(-\lambda)-\omega_2(-\lambda)}e^{-\omega_1(-\lambda)t}\tilde{\hat{f}}(-\lambda,\omega_1(-\lambda),t)+\frac{1}{\omega_1(-\lambda)-\omega_2(-\lambda)}e^{-\omega_2(-\lambda)t}\tilde{\hat{f}}(-\lambda,\omega_2(-\lambda),t)\text{, for }\lambda\in\mathbb{R}.$$

Taking into consideration (1.3) and (1.4), we see that, for $\lambda\in\mathbb{R}$, (2.15) can be written as follows:

$$(2.16)\quad \frac{1}{\omega_1(\lambda)-\omega_2(\lambda)}\omega_2(\lambda)e^{-\omega_1(\lambda)t}-\omega_1(\lambda)e^{-\omega_2(\lambda)t}\}\hat{u}_0(-\lambda)+\frac{1}{\omega_1(\lambda)-\omega_2(\lambda)}(e^{-\omega_1(\lambda)t}-e^{-\omega_2(\lambda)t})\hat{u}_1(-\lambda)+\hat{u}(-\lambda,t)$$

$$=\frac{-i\lambda}{\omega_1(\lambda)-\omega_2(\lambda)}e^{-\omega_1(\lambda)t}\widetilde{g}_0(\omega_1(\lambda),t)-\frac{-i\lambda}{\omega_1(\lambda)-\omega_2(\lambda)}e^{-\omega_2(\lambda)t}\widetilde{g}_0(\omega_2(\lambda),t)$$

$$-\frac{1}{\omega_1(\lambda)-\omega_2(\lambda)}e^{-\omega_1(\lambda)t}\widetilde{g}_1(\omega_1(\lambda),t)+\frac{1}{\omega_1(\lambda)-\omega_2(\lambda)}e^{-\omega_2(\lambda)t}\widetilde{g}_1(\omega_2(\lambda),t)$$

$$-\frac{1}{\omega_1(\lambda)-\omega_2(\lambda)}e^{-\omega_1(\lambda)t}\tilde{\hat{f}}(-\lambda,\omega_1(\lambda),t)+\frac{1}{\omega_1(\lambda)-\omega_2(\lambda)}e^{-\omega_2(\lambda)t}\tilde{\hat{f}}(-\lambda,\omega_2(\lambda),t)\text{, for }\lambda\in\mathbb{R}.$$

***Step 5*** Multiplying (2.16) by $e^{i\lambda x}$, integrating it from $\lambda=-\infty$ to $\lambda=\infty$, and subtracting the resulting equation from (2.14), we obtain:

$$(2.17)\quad \int_{-\infty}^{\infty}\frac{\omega_2(\lambda)e^{i\lambda x-\omega_1(\lambda)t}-\omega_1(\lambda)e^{i\lambda x-\omega_2(\lambda)t}}{\omega_1(\lambda)-\omega_2(\lambda)}\hat{u}_0(\lambda)d\lambda-\int_{-\infty}^{\infty}\frac{\omega_2(\lambda)e^{i\lambda x-\omega_1(\lambda)t}-\omega_1(\lambda)e^{i\lambda x-\omega_2(\lambda)t}}{\omega_1(\lambda)-\omega_2(\lambda)}\hat{u}_0(-\lambda)d\lambda$$

$$+\int_{-\infty}^{\infty}\frac{e^{i\lambda x-\omega_1(\lambda)t}-e^{i\lambda x-\omega_2(\lambda)t}}{\omega_1(\lambda)-\omega_2(\lambda)}\hat{u}_1(\lambda)d\lambda-\int_{-\infty}^{\infty}\frac{e^{i\lambda x-\omega_1(\lambda)t}-e^{i\lambda x-\omega_2(\lambda)t}}{\omega_1(\lambda)-\omega_2(\lambda)}\hat{u}_1(-\lambda)d\lambda+2\pi u(x,t)$$

$$=2i\int_{-\infty}^{\infty}\left\{\frac{\lambda e^{i\lambda x-\omega_1(\lambda)t}}{\omega_1(\lambda)-\omega_2(\lambda)}\widetilde{g}_0(\omega_1(\lambda),t)-\frac{\lambda e^{i\lambda x-\omega_2(\lambda)t}}{\omega_1(\lambda)-\omega_2(\lambda)}\widetilde{g}_0(\omega_2(\lambda),t)\right\}d\lambda$$

$$-\int_{-\infty}^{\infty}\left\{\frac{e^{i\lambda x-\omega_1(\lambda)t}}{\omega_1(\lambda)-\omega_2(\lambda)}\tilde{\hat{f}}(\lambda,\omega_1(\lambda),t)-\frac{e^{i\lambda x-\omega_2(\lambda)t}}{\omega_1(\lambda)-\omega_2(\lambda)}\tilde{\hat{f}}(\lambda,\omega_2(\lambda),t)\right\}d\lambda$$

$$+\int_{-\infty}^{\infty}\left\{\frac{e^{i\lambda x-\omega_1(\lambda)t}}{\omega_1(\lambda)-\omega_2(\lambda)}\tilde{\hat{f}}(-\lambda,\omega_1(\lambda),t)-\frac{e^{i\lambda x-\omega_2(\lambda)t}}{\omega_1(\lambda)-\omega_2(\lambda)}\tilde{\hat{f}}(-\lambda,\omega_2(\lambda),t)\right\}d\lambda,$$

where we used also the fact that

$$\int_{-\infty}^{\infty}e^{i\lambda x}\hat{u}(\lambda,t)=2\pi u(x,t)\text{ and }\int_{-\infty}^{\infty}e^{i\lambda x}\hat{u}(-\lambda,t)=0.$$

This gives (1.5) and completes the formal derivation.

***Remarks*** **(1)** As we will see in the next section, the integrals in (2.14) and (2.17) exist, for $(x,t)\in Q$, $x\neq t$, at least in the generalized sense, and this is crucial for the method to work.

**(2)** With the notation of Lemma 2.1, the solution formula (1.5) can be written as follows:

$$(2.18)\quad 2\pi u(x,t)=\int_{-\infty}^{\infty}\mathcal{K}_0(\lambda,t)[\hat{u}_0(\lambda)-\hat{u}_0(-\lambda)]e^{i\lambda x}d\lambda+\int_{-\infty}^{\infty}\mathcal{K}_1(\lambda,t)[\hat{u}_1(\lambda)-\hat{u}_1(-\lambda)]e^{i\lambda x}d\lambda$$

$$+\int_{-\infty}^{\infty}\mathcal{L}_{g_0}(\lambda,t)e^{i\lambda x}d\lambda+\int_{-\infty}^{\infty}\mathcal{M}_f^{+}(\lambda,t)e^{i\lambda x}d\lambda+\int_{-\infty}^{\infty}\mathcal{M}_f^{-}(\lambda,t)e^{i\lambda x}d\lambda,\ \text{for } (x,t)\in Q,\ x\neq t.$$

Let as note that, according to the proof of Lemma 1, the integrands in (2.18) do not contain the square root function $\rho(\lambda)=\sqrt{\lambda^2-(1/4)}$. As we will see, this fact and the formulas (2.14) – (2.17), obtained in the proof of Lemma 1, are crucial for the proof of Theorems 4, 7 and 8.

## 3. Interpretation of the integrals

***Preliminary remarks***

First, let us note that

$$(3.1)\quad \int_{-1}^{1}[e^{i\lambda x-\omega_1(\lambda)t}\tilde{g}_0(\omega_1(\lambda),t)-e^{i\lambda x-\omega_2(\lambda)t}\tilde{g}_0(\omega_2(\lambda),t)]\frac{\lambda d\lambda}{\rho(\lambda)}$$
$$=\int_{\kappa^-}[e^{i\lambda x-\omega_1(\lambda)t}\tilde{g}_0(\omega_1(\lambda),t)-e^{i\lambda x-\omega_2(\lambda)t}\tilde{g}_0(\omega_2(\lambda),t)]\frac{\lambda d\lambda}{\rho(\lambda)}$$
$$=\int_{\kappa^-}e^{i\lambda x-\omega_1(\lambda)t}\tilde{g}_0(\omega_1(\lambda),t)\frac{\lambda d\lambda}{\rho(\lambda)}-\int_{\kappa^-}e^{i\lambda x-\omega_2(\lambda)t}\tilde{g}_0(\omega_2(\lambda),t)]\frac{\lambda d\lambda}{\rho(\lambda)},$$

where $\kappa^-$ is the semicircle depicted in Fig 4. Furthermore, in the above integrals, the semicircle $\kappa^-$ can be replaced by $-\kappa^+$. Indeed, this follows from the analyticity of the function $e^{i\lambda x}\mathcal{L}_{g_0}(\lambda,t)\lambda$ with respect to $\lambda\in\mathbb{C}$ (by Lemma 1), which (by Cauchy's theorem) implies that $\int_{\kappa^-+\kappa^+}e^{i\lambda x}\mathcal{L}_{g_0}(\lambda,t)\lambda d\lambda=0$.

Also,

$$(3.2)\quad \int_{-1}^{1}\left[e^{i\lambda x-\omega_1(\lambda)t}\tilde{\hat{f}}(-\lambda,\omega_1(\lambda),t)-e^{i\lambda x-\omega_2(\lambda)t}\tilde{\hat{f}}(-\lambda,\omega_2(\lambda),t)\right]\frac{d\lambda}{2i\rho(\lambda)}$$
$$=\int_{-1}^{1}\left[e^{-i\lambda x-\omega_1(\lambda)t}\tilde{\hat{f}}(\lambda,\omega_1(\lambda),t)-e^{-i\lambda x-\omega_2(\lambda)t}\tilde{\hat{f}}(\lambda,\omega_2(\lambda),t)\right]\frac{d\lambda}{2i\rho(\lambda)}$$
$$=\int_{\kappa^-}e^{-i\lambda x-\omega_1(\lambda)t}\tilde{\hat{f}}(\lambda,\omega_1(\lambda),t)\frac{d\lambda}{2i\rho(\lambda)}-\int_{\kappa^-}e^{-i\lambda x-\omega_2(\lambda)t}\tilde{\hat{f}}(\lambda,\omega_2(\lambda),t)\frac{d\lambda}{2i\rho(\lambda)}$$

and analogous formulas hold for all the integrals in (1.5) "restricted" on the interval $[-1,1]$. Thus, the main difficulty in dealing with the integrals in (1.5) concerns their part taken on $\mathbb{L}:=(-\infty,-1]+[1,\infty)$. With this notation

$$\int_{\mathbb{L}}=\int_{-\infty}^{-1}+\int_{1}^{\infty},\quad \int_{-1}^{1}=\int_{\kappa^-}\ \text{and}\ \int_{-\infty}^{\infty}=\int_{-\infty}^{-1}+\int_{-1}^{1}+\int_{1}^{\infty}=\int_{-1}^{1}+\int_{\mathbb{L}}=\int_{\kappa^-}+\int_{\mathbb{L}},$$

where the last two equations hold for appropriate analytic integrands.

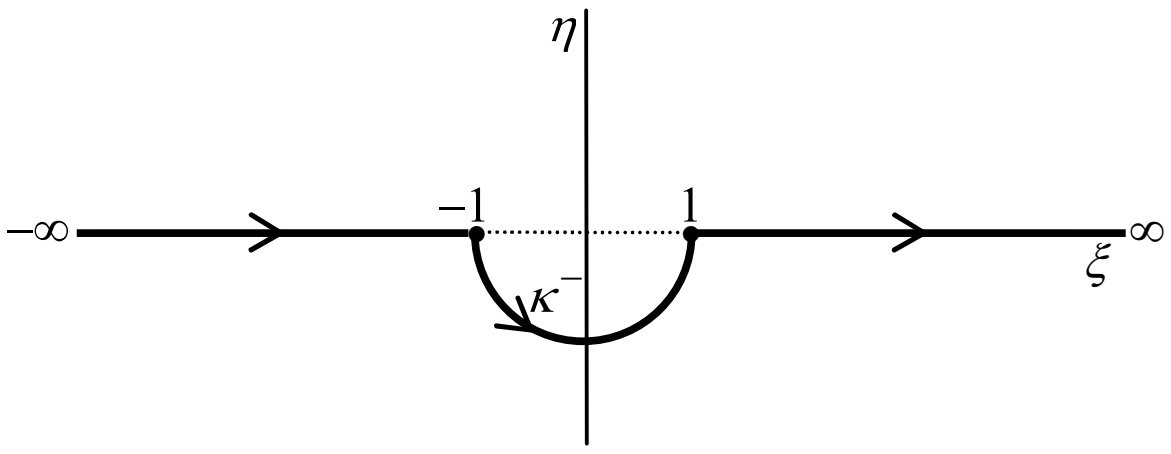

**Fig. 4** The two half-lines $\mathbb{L}=(-\infty,-1]+[1,\infty)$ and the semicircle $\kappa^-$

***Behavior of*** $\rho(\lambda)$, $\omega_1(\lambda)$ ***and*** $\omega_2(\lambda)$ ***as*** $\lambda\to\infty$ *with* $\lambda\in\mathbb{C}$

Let us note that, as $\lambda\to\infty$,

$$\sqrt{\left|\lambda^2-\frac{1}{4}\right|}\cong|\lambda| \text{ in the sense that } \lim_{\lambda\to\infty}\frac{1}{|\lambda|}\sqrt{\left|\lambda^2-\frac{1}{4}\right|}=1 \text{ and } \lim_{\lambda\to\infty}\left(\sqrt{\left|\lambda^2-\frac{1}{4}\right|}-|\lambda|\right)=0 .$$

Also, since

$$|\lambda||\sin[\theta_1(\lambda)-\theta(\lambda)]|<\frac{1}{2} \text{ and } |\lambda||\sin[\theta(\lambda)-\theta_2(\lambda)]|<\frac{1}{2} \text{ (see Fig 5),}$$

it follows that

$$\sup_\lambda\{|\lambda||\theta_1(\lambda)-\theta(\lambda)|\}<\infty \text{ and } \sup_\lambda\{|\lambda||\theta(\lambda)-\theta_2(\lambda)|\}<\infty .$$

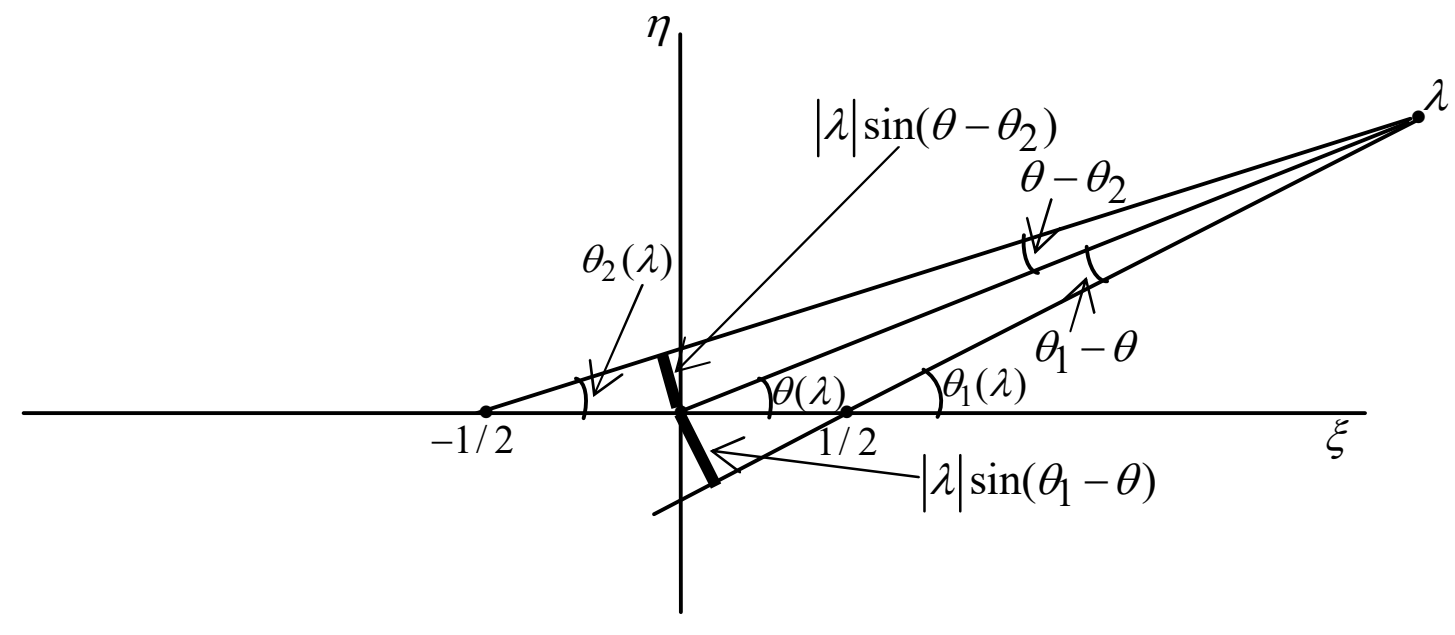

**Fig. 5**

Therefore,

$$\sup_\lambda\left||\lambda|\sin\left(\tfrac{\theta_1(\lambda)+\theta_2(\lambda)}{2}\right)-|\lambda|\sin\theta(\lambda)\right|<\infty .$$

Thus

$$\exp\left[t\sqrt{\left|\lambda^2-\tfrac{1}{4}\right|}\sin\left(\tfrac{\theta_1(\lambda)+\theta_2(\lambda)}{2}\right)\right]\approx e^{t|\lambda|\sin\theta(\lambda)} \tag{3.3}$$

in the sense that

$$\sup_\lambda\frac{\exp\left[t\sqrt{\left|\lambda^2-\tfrac{1}{4}\right|}\sin\left(\tfrac{\theta_1(\lambda)+\theta_2(\lambda)}{2}\right)\right]}{e^{t|\lambda|\sin\theta(\lambda)}}<\infty, \text{ uniformly for } t \text{ in compact subsets of } [0,\infty) .$$

Also, since

$$\tfrac{\theta_1(\lambda)+\theta_2(\lambda)}{2}-\theta(\lambda)\to 0, \text{ as } \lambda\to\infty ,$$

we have

$$\frac{\exp\left[i\left(\tfrac{\theta_1(\lambda)+\theta_2(\lambda)}{2}\right)\right]}{e^{i\theta(\lambda)}}\to 1 ,$$

and, therefore,

$$\rho(\lambda)=\sqrt{\left|\lambda^2-\tfrac{1}{4}\right|}\exp\left[i\left(\tfrac{\theta_1(\lambda)+\theta_2(\lambda)}{2}\right)\right]\approx\lambda, \text{ i.e., } \lim_{\lambda\to\infty}\frac{\rho(\lambda)}{\lambda}=1,\ \lim_{\lambda\to\infty}\frac{\omega_1(\lambda)}{\lambda}=i \text{ and } \lim_{\lambda\to\infty}\frac{\omega_2(\lambda)}{\lambda}=-i . \tag{3.4}$$

***Behavior of*** $e^{i\lambda x-\omega_j(\lambda)t}$ ***as*** $\lambda\to\infty$ In view of (3.3), we obtain the following lemma, which is crucial to the study of the behavior of the integrals of the solution.

***Lemma 2*** *Writing* $\lambda=\xi+i\eta=|\lambda|e^{i\theta}$, $\eta=|\lambda|\sin\theta$, *we have*

$$\left|e^{i\lambda x-\omega_1(\lambda)t}\right|=e^{-t/2}\exp\left[-\eta x+t\sqrt{\left|\lambda^2-\frac{1}{4}\right|}\sin\left(\frac{\theta_1+\theta_2}{2}\right)\right]\approx e^{-t/2}e^{-(x-t)|\lambda|\sin\theta} \tag{3.5}$$

*and*

$$\left|e^{i\lambda x-\omega_2(\lambda)t}\right|=e^{-t/2}\exp\left[-\eta x-t\sqrt{\left|\lambda^2-\frac{1}{4}\right|}\sin\left(\frac{\theta_1+\theta_2}{2}\right)\right]\approx e^{-t/2}e^{-(x+t)|\lambda|\sin\theta} , \tag{3.6}$$

*as* $\lambda\to\infty$ *(* $\lambda\in\mathbb{C}$ *), uniformly for* $t$ *in compact subsets of* $[0,\infty)$.

***Jordan-type lemma*** The following lemma will be used in the deformation of the contours of the various integrals. Its conclusions follow from (3.5) and (3.6).

***Lemma 3*** *Let $\varphi(\lambda)$ be a continuous function, defined for $\lambda \in \mathbb{C}$, such that $\lim_{|\lambda|\to\infty} \varphi(\lambda) = 0$. Then the following hold:*

***(1)*** *For $x > t \geq 0$,*

$$\lim_{A\to\infty} \int_{\{\lambda\in\mathbb{C}:\, \mathrm{Im}\lambda\geq 0,\, |\lambda|=A\}} e^{i\lambda x-\omega_1(\lambda)t}\varphi(\lambda)d\lambda = 0,$$

$$\lim_{A\to\infty} \int_{\{\lambda\in\mathbb{C}:\, \mathrm{Im}\lambda\geq 0,\, |\lambda|=A\}} e^{i\lambda x-\omega_2(\lambda)t}\varphi(\lambda)d\lambda = 0.$$

***(2)*** *For $0 \leq x < t$,*

$$\lim_{A\to\infty} \int_{\{\lambda\in\mathbb{C}:\, \mathrm{Im}\lambda\leq 0,\, |\lambda|=A\}} e^{i\lambda x-\omega_1(\lambda)t}\varphi(\lambda)d\lambda = 0,$$

$$\lim_{A\to\infty} \int_{\{\lambda\in\mathbb{C}:\, \mathrm{Im}\lambda\geq 0,\, |\lambda|=A\}} e^{i\lambda x-\omega_2(\lambda)t}\varphi(\lambda)d\lambda = 0.$$

*In particular, by Cauchy's theorem, we have:*

***(3)*** *If, in addition, $\varphi(\lambda)$ is analytic in an open neighborhood of $\{\lambda \in \mathbb{C} : \mathrm{Im}\lambda \geq 0, |\lambda| \geq 1\}$, then*

$$\int_{\mathbb{L}} e^{i\lambda x-\omega_1(\lambda)t}\varphi(\lambda)d\lambda = \int_{\kappa^+} e^{i\lambda x-\omega_1(\lambda)t}\varphi(\lambda)d\lambda, \text{ for } x > t \geq 0,$$

$$\int_{\mathbb{L}} e^{i\lambda x-\omega_2(\lambda)t}\varphi(\lambda)d\lambda = \int_{\kappa^+} e^{i\lambda x-\omega_2(\lambda)t}\varphi(\lambda)d\lambda, \text{ for } x \geq 0,\ t \geq 0,\ x+t > 0.$$

***(4)*** *If, in addition, $\varphi(\lambda)$ is analytic in an open neighborhood of $\{\lambda \in \mathbb{C} : \mathrm{Im}\lambda \leq 0, |\lambda| \geq 1\}$, then*

$$\int_{\mathbb{L}} e^{i\lambda x-\omega_1(\lambda)t}\varphi(\lambda)d\lambda = \int_{-\kappa^-} e^{i\lambda x-\omega_1(\lambda)t}\varphi(\lambda)d\lambda, \text{ for } t > x \geq 0.$$

***Some integration by parts formulas***

We have

(3.7) $$\hat{u}_0(\lambda) = \sum_{m=1}^{N} \frac{d^{m-1}u_0}{dx^{m-1}}(0)\frac{1}{(i\lambda)^m} + \frac{1}{(i\lambda)^N}\int_{y=0}^{\infty} e^{-i\lambda y}\frac{d^N u_0(y)}{dy^N}dy = O(1/\lambda^{N+1}), \text{ as } \lambda\to\infty \text{ with } \mathrm{Im}\,\lambda \leq 0.$$

In particular,

(3.8) $$\hat{u}_0(\lambda) = \frac{u_0(0)}{i\lambda} + \frac{(u_0')\hat{}(\lambda)}{i\lambda} = O(1/\lambda) \text{ as } \lambda\to\infty \text{ with } \lambda\in\mathbb{C},\ \mathrm{Im}\lambda \leq 0.$$

Since

(3.9) $$e^{-\omega(\lambda)t}\tilde{g}_0(\omega(\lambda),t) = \frac{1}{\omega(\lambda)}g_0(t) - \frac{1}{\omega(\lambda)}e^{-\omega(\lambda)t}g_0(0) - \frac{1}{\omega(\lambda)}e^{-\omega(\lambda)t}(g_0')\tilde{}(\omega(\lambda),t),$$

we have

(3.10) $e^{-\omega(\lambda)t}\tilde{g}_0(\omega(\lambda),t) = O(1/\omega(\lambda))$ as $\lambda\to\infty$, with $\mathrm{Re}\,\omega(\lambda)\geq 0$, uniformly for $t$ in compact subsets of $[0,\infty)$.

More generally, setting

$$\rho_N(\lambda,\omega(\lambda),t) := \frac{g_0(t)}{\omega(\lambda)} - \frac{g_0'(t)}{[\omega(\lambda)]^2} + \cdots + (-1)^{N-1}\frac{g_0^{(N-1)}(t)}{[\omega(\lambda)]^N} \quad (\lambda\in\mathbb{C},\ \omega(\lambda)\neq 0),$$

we have

(3.11) $$\tilde{g}_0(\omega(\lambda),t) = e^{\omega(\lambda)t}\rho_{\mathrm{N}}(\lambda,\omega(\lambda),t) - \rho_{\mathrm{N}}(\lambda,\omega(\lambda),0) + \frac{(-1)^{\mathrm{N}}}{[\omega(\lambda)]^{\mathrm{N}}}\int_{\tau=0}^{t} e^{\omega(\lambda)\tau} g_0^{(\mathrm{N})}(\tau)d\tau .$$

Also

$$\hat{f}(\lambda,t) = \int_{y=0}^{\infty} e^{-i\lambda y} f(y,t)dy = \frac{1}{i\lambda} f(0,t) + \frac{1}{i\lambda}\int_{y=0}^{\infty} e^{-i\lambda y}\frac{\partial f(y,t)}{\partial y}dy \quad (\lambda\in\mathbb{C}-\{0\}\,, \mathrm{Im}\,\lambda\le 0\,,\ t\ge 0)$$

and

(3.12) $$\begin{aligned}\tilde{\hat{f}}(\lambda,\omega(\lambda),t) &= \int_{\tau=0}^{t} e^{\omega(\lambda)\tau}\hat{f}(\lambda,\tau)d\tau \\ &= \frac{1}{i\lambda}\int_{\tau=0}^{t}\frac{\partial}{\partial\tau}\left(\frac{e^{\omega(\lambda)\tau}}{\omega(\lambda)}\right) f(0,\tau)d\tau + \frac{1}{i\lambda}\int_{\tau=0}^{t}\frac{\partial}{\partial\tau}\left(\frac{e^{\omega(\lambda)\tau}}{\omega(\lambda)}\right)\left(\int_{y=0}^{\infty} e^{-i\lambda y}\frac{\partial f(y,\tau)}{\partial y}dy\right)d\tau \\ &= \frac{1}{i\lambda}\frac{e^{\omega(\lambda)t}}{\omega(\lambda)} f(0,t) - \frac{1}{i\lambda}\frac{1}{\omega(\lambda)} f(0,0) - \frac{1}{i\lambda}\frac{1}{\omega(\lambda)}\int_{\tau=0}^{t} e^{\omega(\lambda)\tau}\frac{\partial}{\partial\tau}[f(0,\tau)]d\tau \\ &\quad + \frac{1}{i\lambda}\frac{e^{\omega(\lambda)t}}{\omega(\lambda)}\left(\int_{y=0}^{\infty} e^{-i\lambda y}\frac{\partial f(y,t)}{\partial y}dy\right) - \frac{1}{i\lambda}\frac{1}{\omega(\lambda)}\left(\int_{y=0}^{\infty} e^{-i\lambda y}\frac{\partial f(y,0)}{\partial y}dy\right) \\ &\quad - \frac{1}{i\lambda}\frac{1}{\omega(\lambda)}\int_{\tau=0}^{t} e^{\omega(\lambda)\tau}\left(\int_{y=0}^{\infty} e^{-i\lambda y}\frac{\partial^2 f(y,t)}{\partial\tau\partial y}dy\right)d\tau .\end{aligned}$$

It follows from (3.12) that

(3.13) $$e^{-\omega(\lambda)t}\tilde{\hat{f}}(\lambda,\omega(\lambda),t) = \mathrm{O}\big(1/[\lambda\omega(\lambda)]\big) \text{ for } \lambda\to\infty ,$$

with $\mathrm{Re}\,\omega(\lambda)\ge 0$ and $\mathrm{Im}\,\lambda\le 0$, uniformly for $t$ in compact subsets of $[0,\infty)$.

Finally, generalizing (3.7), we obtain, for $\mathrm{Im}\,\lambda\le 0$, $\lambda\ne 0$,

$$\hat{f}(\lambda,t) = h_{f,\mathrm{N}}(\lambda,t) + \frac{1}{(i\lambda)^N}\int_{y=0}^{\infty} e^{-i\lambda y}\frac{\partial^N f(y,t)}{\partial y^N}dy \text{ where } h_{f,N}(\lambda,t) := \sum_{m=1}^{N}\frac{1}{(i\lambda)^m}\left.\frac{\partial^{m-1} f(y,t)}{\partial y^{m-1}}\right|_{y=0} .$$

Therefore,

(3.14) $$\delta_{f,N}(\lambda,t) := \hat{f}(\lambda,t) - h_{f,N}(\lambda,t) = \mathrm{O}(1/\lambda^{N+1}) ,\text{ for } \lambda\to\infty ,\text{ with } \mathrm{Im}\,\lambda\le 0 .$$

Also

(3.15) $$\begin{aligned}\tilde{\hat{f}}(\lambda,\omega(\lambda),t) &= \int_{\tau=0}^{t} e^{\omega(\lambda)\tau}\hat{f}(\lambda,\tau)d\tau \\ &= \int_{\tau=0}^{t} e^{\omega(\lambda)\tau}\delta_{f,N}(\lambda,\tau)d\tau + \int_{\tau=0}^{t} e^{\omega(\lambda)\tau}h_{f,N}(\lambda,\tau)d\tau = \tilde{\delta}_{f,N}(\lambda,\omega(\lambda),t) + \tilde{h}_{f,N}(\lambda,\omega(\lambda),t) \ (\mathrm{Im}\,\lambda\le 0)\end{aligned}$$

and

(3.16) $$\tilde{h}_{f,N}(\lambda,\omega(\lambda),t) = e^{\omega(\lambda)t}\mu_{f,N,M}(\lambda,\omega(\lambda),t) - \mu_{f,N,M}(\lambda,\omega(\lambda),0) + \frac{(-1)^M}{[\omega(\lambda)]^M}\int_{\tau=0}^{t} e^{\omega(\lambda)\tau}h_{f,N}^{(M)}(\lambda,\tau)d\tau$$

where

$$\mu_{f,N,M}(\lambda,\omega(\lambda),t) := \frac{h_{f,N}(\lambda,t)}{\omega(\lambda)} - \frac{h_{f,N}^{(1)}(\lambda,t)}{[\omega(\lambda)]^2} + \cdots + (-1)^{\mathrm{M}-1}\frac{h_{f,N}^{(\mathrm{M}-1)}(\lambda,t)}{[\omega(\lambda)]^M} \quad (\lambda\omega(\lambda)\ne 0).$$

(The derivatives of $h_{f,N}(\lambda,t)$, in the above quantity, are taken with respect to $t$.)

Besides the above identities, crucial for the defomation and limiting processes are the estimates of the following lemma.

***Lemma 4*** *We have*

*(3.17)* $$e^{-\omega(\lambda)t}\tilde{\delta}_{f,N}(\lambda,\omega(\lambda),t) = \mathrm{O}\big(1/[\lambda^{N+1}\omega(\lambda)]\big),\ \textit{for } \lambda\to\infty \textit{ with } \mathrm{Re}\,\omega(\lambda)\ge 0 \textit{ and } \mathrm{Im}\,\lambda\le 0 ,$$

*and*

(3.18) $$\frac{(-1)^M}{[\omega(\lambda)]^M} e^{-\omega(\lambda)t} \int_{\tau=0}^{t} e^{\omega(\lambda)\tau} h_{f,N}^{(M)}(\lambda,\tau)d\tau = \mathrm{O}\left(1/[\lambda[\omega(\lambda)]^{M+1}]\right), \textit{ for } \lambda\to\infty \textit{ with } \operatorname{Re}\omega(\lambda)\ge 0,$$

*uniformly for* $t$ *in compact subsets of* $[0,\infty)$.

Now we proceed with the interpretation of the integrals in (1.5). The integrands of $\mathcal{U}_0^+(u_0)$, $\mathcal{U}_0^-(u_0)$ and $\mathcal{G}(g_0)$ are $\mathrm{O}(1/\lambda)$, as $\lambda\to\pm\infty$, and – in general – these integrals do not converge absolutely, hence they are considered in the generalized sense, i.e., as limits of the form $\int_{-\infty}^{\infty} = \lim_{A\to\infty}\int_{-A}^{A}$. (The existence of these limits will follow from the interpretation of the integrals that we will give below.) Moreover these integrals become worse if we differentiate their integrands with respect to $x$ or $t$. The following formulas, however, interpret the integrals in such a way that differentiation under the integral sign can be justified (once they have appropriately interpreted) and, thus, it is allowed. As a matter of fact the difficulty of differentiating under the integral sign occurs for all the integrals in (1.2). We will deal only with the more difficult cases – the remaining cases are either easier or similar to the treated ones. The formulas (3.3) – (3.18) will be used without any particular reference to them.

***In the case*** $x>t$, the following hold:

(3.19) $$\int_{\mathbb{L}} \omega_2(\lambda) e^{i\lambda x-\omega_1(\lambda)t}\hat{u}_0(\lambda)\frac{d\lambda}{2i\rho(\lambda)}$$
$$= \sum_{j=1}^{N} u_0^{(j-1)}(0)\int_{\kappa^+}\omega_2(\lambda)e^{i\lambda x-\omega_1(\lambda)t}\frac{1}{(i\lambda)^j}\frac{d\lambda}{2i\rho(\lambda)} + \int_{\mathbb{L}}\omega_2(\lambda)e^{i\lambda x-\omega_1(\lambda)t}\frac{1}{(i\lambda)^N}(u_0^{(N)})\hat{}(\lambda)\frac{d\lambda}{2i\rho(\lambda)}.$$

As a matter of fact, in formula (3.19), instead of $\kappa^+$, we can use any contour $\Gamma^+$ in the upper half-plane, starting at $\lambda=+1$ and ending at the point $\lambda=-1$. In particular $\Gamma^+$ can be the sum $\varepsilon_1+\varepsilon_2$ of the rays $\varepsilon_1$ and $\varepsilon_2$, depicted in Fig 6.

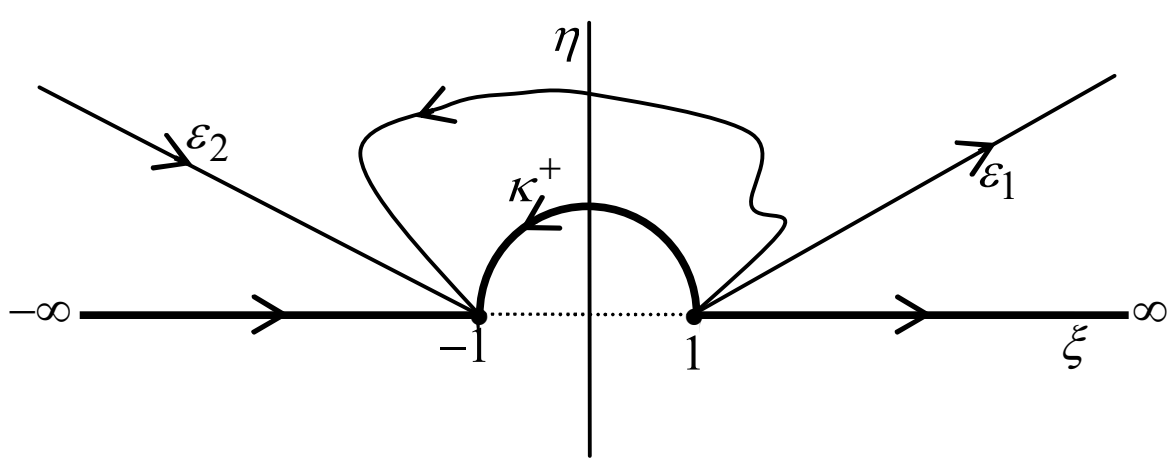

**Fig. 6** Various choices of $\Gamma^+$

(3.20) $$\int_{\mathbb{L}}\omega_2(\lambda)e^{i\lambda x-\omega_1(\lambda)t}\hat{u}_0(-\lambda)\frac{d\lambda}{2i\rho(\lambda)} = \int_{\kappa^+}\omega_2(\lambda)e^{i\lambda x-\omega_1(\lambda)t}\hat{u}_0(-\lambda)\frac{d\lambda}{2i\rho(\lambda)}.$$

(3.21) $$\int_{\mathbb{L}} e^{i\lambda x-\omega_1(\lambda)t}\tilde{g}_0(\omega_1(\lambda),t)\frac{\lambda d\lambda}{\rho(\lambda)} = \sum_{j=1}^{N}(-1)^{j-1}\left\{g_0^{(j-1)}(t)\int_{\kappa^+}\frac{e^{i\lambda x}}{[\omega_1(\lambda)]^j}\frac{\lambda d\lambda}{\rho(\lambda)} - g_0^{(j-1)}(0)\int_{\kappa^+}\frac{e^{i\lambda x-\omega_1(\lambda)t}}{[\omega_1(\lambda)]^j}\frac{\lambda d\lambda}{\rho(\lambda)}\right\}$$
$$+(-1)^N\int_{\mathbb{L}}\frac{1}{[\omega_1(\lambda)]^N}e^{i\lambda x-\omega_1(\lambda)t}(g_0^{(N)})\tilde{}(\omega_1(\lambda),t)\frac{\lambda d\lambda}{\rho(\lambda)}.$$

In the last integral in (3.21), $\mathbb{L}$ can be replaced by $\kappa^+$ or any $\Gamma^+$ ($N\ge 1$). Indeed, this follows from the fact that, in view of (3.5),

$$\left|e^{\omega_1(\lambda)\tau}\right| \approx e^{-t/2}e^{-\tau|\lambda|\sin\theta(\lambda)},$$

which implies that, for $\lambda\in\mathbb{C}$ with $\operatorname{Im}\lambda\ge 0$,

$$\left|(g_0^{(N)})\tilde{}(\omega_1(\lambda),t)\right| = \left|\int_{\tau=0}^{t} g_0^{(N)}(\tau)e^{\omega_1(\lambda)\tau}d\tau\right| \le \int_{\tau=0}^{t}\left|g_0^{(N)}(\tau)\right|\left|e^{\omega_1(\lambda)\tau}\right|d\tau$$

$$\preceq \int_{\tau=0}^{t} e^{\tau/2}\left|g_0^{(N)}(\tau)\right| e^{-\tau|\lambda|\sin\theta(\lambda)}d\tau \preceq \int_{\tau=0}^{t} e^{\tau/2}\left|g_0^{(N)}(\tau)\right| d\tau .$$

$$(3.22)\quad \int_{\mathbb{L}} e^{i\lambda x-\omega_1(\lambda)t}\widetilde{\hat{f}}(\pm\lambda,\omega_1(\lambda),t)\frac{d\lambda}{2i\rho(\lambda)} = \int_{\mathbb{L}} e^{i\lambda x-\omega_1(\lambda)t}\widetilde{\delta}_{f,N}(\pm\lambda,\omega_1(\lambda),t)\frac{d\lambda}{2i\rho(\lambda)}$$
$$+\int_{\kappa^+} e^{i\lambda x}\mu_{f,N,M}(\pm\lambda,\omega_1(\lambda),t)\frac{d\lambda}{2i\rho(\lambda)} - \int_{\kappa^+} e^{i\lambda x-\omega_1(\lambda)t}\mu_{f,N,M}(\pm\lambda,\omega_1(\lambda),0)\frac{d\lambda}{2i\rho(\lambda)}$$
$$+\int_{\mathbb{L}} e^{i\lambda x-\omega_1(\lambda)t}\frac{(-1)^M}{[\omega_1(\lambda)]^M}\int_{\tau=0}^{t} e^{\omega_1(\lambda)\tau}h_{f,N}^{(M)}(\pm\lambda,\tau)d\tau\frac{d\lambda}{2i\rho(\lambda)}.$$

In the last integral in (3.22), $\mathbb{L}$ can be replaced by $\kappa^+$.

Also,

$$\int_{\mathbb{L}} e^{i\lambda x-\omega_1(\lambda)t}\widetilde{\hat{f}}(-\lambda,\omega_1(\lambda),t)\frac{d\lambda}{2i\rho(\lambda)} = \int_{\kappa^+} e^{i\lambda x-\omega_1(\lambda)t}\widetilde{\hat{f}}(-\lambda,\omega_1(\lambda),t)\frac{d\lambda}{2i\rho(\lambda)}.$$

***In the case case*** $x<t$, the following hold:

$$(3.23)\quad \int_{\mathbb{L}} e^{i\lambda x-\omega_1(\lambda)t}\widetilde{g}_0(\omega_1(\lambda),t)\frac{\lambda d\lambda}{\rho(\lambda)} = \sum_{j=1}^{N}(-1)^{j-1}\left\{g_0^{(j-1)}(t)\int_{\kappa^+}\frac{e^{i\lambda x}}{[\omega_1(\lambda)]^j}\frac{\lambda d\lambda}{\rho(\lambda)} - g_0^{(j-1)}(0)\int_{-\kappa^-}\frac{e^{i\lambda x-\omega_1(\lambda)t}}{[\omega_1(\lambda)]^j}\frac{\lambda d\lambda}{\rho(\lambda)}\right\}$$
$$+(-1)^N\int_{\mathbb{L}}\frac{1}{[\omega_1(\lambda)]^N}e^{i\lambda x-\omega_1(\lambda)t}(g_0^{(N)}\widetilde{)}(\omega_1(\lambda),t)\frac{\lambda d\lambda}{\rho(\lambda)}$$

where $\kappa^-$ is the semicircle depicted Fig 4. (This time, in the last integral in (3.23), $\mathbb{L}$ ***cannot*** be replaced by $\kappa^+$ or $-\kappa^-$, in general.)

$$(3.24)\quad \int_{\mathbb{L}}\omega_2(\lambda)e^{i\lambda x-\omega_1(\lambda)t}\hat{u}_0(\lambda)\frac{d\lambda}{2i\rho(\lambda)} = \int_{-\kappa^-}\omega_2(\lambda)e^{i\lambda x-\omega_1(\lambda)t}\hat{u}_0(\lambda)\frac{d\lambda}{2i\rho(\lambda)}.$$

$$(3.25)\quad \int_{\mathbb{L}}\omega_2(\lambda)e^{i\lambda x-\omega_1(\lambda)t}\hat{u}_0(-\lambda)\frac{d\lambda}{2i\rho(\lambda)} = \sum_{j=1}^{N}u_0^{(j-1)}(0)\int_{-\kappa^-}\omega_2(\lambda)e^{i\lambda x-\omega_1(\lambda)t}\frac{1}{(-i\lambda)^j}\frac{d\lambda}{2i\rho(\lambda)}$$
$$+\int_{\mathbb{L}}\omega_2(\lambda)e^{i\lambda x-\omega_1(\lambda)t}\frac{1}{(-i\lambda)^N}(u_0^{(N)}\hat{)}(-\lambda)\frac{d\lambda}{2i\rho(\lambda)}.$$

$$(3.26)\quad \int_{\mathbb{L}} e^{i\lambda x-\omega_1(\lambda)t}\widetilde{\hat{f}}(\pm\lambda,\omega_1(\lambda),t)\frac{d\lambda}{2i\rho(\lambda)} = \int_{\mathbb{L}} e^{i\lambda x-\omega_1(\lambda)t}\widetilde{\delta}_{f,N}(\pm\lambda,\omega_1(\lambda),t)\frac{d\lambda}{2i\rho(\lambda)}$$
$$+\int_{\kappa^+} e^{i\lambda x}\mu_{f,N,M}(\pm\lambda,\omega_1(\lambda),t)\frac{d\lambda}{2i\rho(\lambda)} - \int_{-\kappa^-} e^{i\lambda x-\omega_1(\lambda)t}\mu_{f,N,M}(\pm\lambda,\omega_1(\lambda),0)\frac{d\lambda}{2i\rho(\lambda)}$$
$$+\int_{\mathbb{L}} e^{i\lambda x-\omega_1(\lambda)t}\frac{(-1)^M}{[\omega_1(\lambda)]^M}\int_{\tau=0}^{t} e^{\omega_1(\lambda)\tau}h_{f,N}^{(M)}(\pm\lambda,\tau)d\tau\frac{d\lambda}{2i\rho(\lambda)}.$$

***In any case:*** $x>t$ ***or*** $x=t$ ***or*** $x<t$, the following hold:

$$(3.27)\quad \int_{\mathbb{L}} e^{i\lambda x-\omega_2(\lambda)t}\widetilde{g}_0(\omega_2(\lambda),t)\frac{\lambda d\lambda}{\rho(\lambda)} = \int_{\kappa^+} e^{i\lambda x-\omega_2(\lambda)t}\widetilde{g}_0(\omega_2(\lambda),t)\frac{\lambda d\lambda}{\rho(\lambda)}.$$

$$(3.28)\quad \int_{\mathbb{L}}\omega_1(\lambda)e^{i\lambda x-\omega_2(\lambda)t}\hat{u}_0(\lambda)\frac{d\lambda}{2i\rho(\lambda)} = \sum_{j=1}^{N}u_0^{(j-1)}(0)\int_{\kappa^+}\omega_1(\lambda)e^{i\lambda x-\omega_2(\lambda)t}\frac{1}{(i\lambda)^j}\frac{d\lambda}{2i\rho(\lambda)}$$
$$+\int_{\mathbb{L}}\omega_1(\lambda)e^{i\lambda x-\omega_2(\lambda)t}\frac{1}{(i\lambda)^N}(u_0^{(N)}\hat{)}(\lambda)\frac{d\lambda}{2i\rho(\lambda)}.$$

(3.29) $$\int_{\mathbb{L}} \omega_1(\lambda) e^{i\lambda x-\omega_2(\lambda)t} \hat{u}_0(-\lambda) \frac{d\lambda}{2i\rho(\lambda)} = \int_{\kappa^+} \omega_1(\lambda) e^{i\lambda x-\omega_2(\lambda)t} \hat{u}_0(-\lambda) \frac{d\lambda}{2i\rho(\lambda)}.$$

(3.30) $$\int_{\mathbb{L}} e^{i\lambda x-\omega_2(\lambda)t} \tilde{\hat{f}}(\lambda,\omega_2(\lambda),t) \frac{d\lambda}{2i\rho(\lambda)} = \int_{\mathbb{L}} e^{i\lambda x-\omega_2(\lambda)t} \tilde{\delta}_{f,N}(\lambda,\omega_2(\lambda),t) \frac{d\lambda}{2i\rho(\lambda)}$$
$$+ \int_{\kappa^+} e^{i\lambda x} \mu_{f,N,M}(\lambda,\omega_2(\lambda),t) \frac{d\lambda}{2i\rho(\lambda)} - \int_{\kappa^+} e^{i\lambda x-\omega_2(\lambda)t} \mu_{f,N,M}(\lambda,\omega_2(\lambda),0) \frac{d\lambda}{2i\rho(\lambda)}$$
$$+ \int_{\mathbb{L}} e^{i\lambda x-\omega_2(\lambda)t} \frac{(-1)^M}{[\omega_2(\lambda)]^M} \int_{\tau=0}^{t} e^{\omega_2(\lambda)\tau} h_{f,N}^{(M)}(\lambda,\tau) d\tau \frac{d\lambda}{2i\rho(\lambda)}.$$

(3.31) $$\int_{\mathbb{L}} e^{i\lambda x-\omega_2(\lambda)t} \tilde{\hat{f}}(-\lambda,\omega_2(\lambda),t) \frac{d\lambda}{2i\rho(\lambda)} = \int_{\kappa^+} e^{i\lambda x-\omega_2(\lambda)t} \tilde{\hat{f}}(-\lambda,\omega_2(\lambda),t) \frac{d\lambda}{2i\rho(\lambda)}.$$

***Remarks*** **(1)** The above equations give the interpretation of the integrals – in the cases the integral does not converge absolutely.

**(2)** Formulas, similar to the case of $\mathcal{U}_0^+(u_0)$, $\mathcal{U}_0^-(u_0)$, can be written also in the case of the integrals $\mathcal{U}_1^+(u_1)$, $\mathcal{U}_1^-(u_1)$, respectively. (The integrals $\mathcal{U}_1^+(u_1)(x,t)$, $\mathcal{U}_1^-(u_1)(x,t)$ converge absolutely, but we need these formulas in order to deal also with their derivatives $\partial^{n+m}[\mathcal{U}_1^{\pm}(u_1)(x,t)]/\partial x^n \partial t^m$.)

**(3)** In all the above formulas (3.19) – (3.31), the semicircle $\kappa^+$ can be replaced by any contour $\Gamma^+$. Similarly, $-\kappa^-$ can be replaced by any contour $\Gamma^-$ in the ***lower*** half-plane, starting at $\lambda=+1$ and ending at the point $\lambda=-1$. (This is analogous to $\Gamma^+$ of Fig 6.)

**(4)** We claim that for $x>t$,

(3.32) $$\int_{-\infty}^{\infty} [e^{i\lambda x-\omega_1(\lambda)t} \tilde{g}_0(\omega_1(\lambda),t) - e^{i\lambda x-\omega_2(\lambda)t} \tilde{g}_0(\omega_2(\lambda),t)] \frac{\lambda d\lambda}{\rho(\lambda)} = 0.$$

Indeed, with the notation of Lemma 1, we write the above integral as

$$\int_{-\infty}^{\infty} e^{i\lambda x} \mathcal{L}_{g_0}(\lambda,t) d\lambda = \int_{-1}^{1} e^{i\lambda x} \mathcal{L}_{g_0}(\lambda,t) d\lambda + \int_{\mathbb{L}} e^{i\lambda x} \mathcal{L}_{g_0}(\lambda,t) d\lambda$$
$$= \int_{-1}^{1} e^{i\lambda x} \mathcal{L}_{g_0}(\lambda,t) d\lambda + \int_{\kappa^+} e^{i\lambda x} \mathcal{L}_{g_0}(\lambda,t) d\lambda = 0,$$

where the last equation follows from that fact that $e^{i\lambda x}\mathcal{L}_{g_0}(\lambda,t)$ is an entire function of $\lambda$ (by Lemma 1).

## 4. Proof of Theorem 1

***Part 1*** We start with the important observation that all the integrals in (1.5) exist at least in the generalized sense, provided that $(x,t)\in Q$ and $x\neq t$. For example, for $x>t>0$, in the spirit of (3.20),

(4.1) $$\int_{\mathbb{L}} \omega_2(\lambda) e^{i\lambda x-\omega_1(\lambda)t} \hat{u}_0(\lambda) \frac{d\lambda}{2i\rho(\lambda)} = \int_{\kappa^+} \omega_2(\lambda) e^{i\lambda x-\omega_1(\lambda)t} \frac{u_0(0)}{i\lambda} \frac{d\lambda}{2i\rho(\lambda)} + \int_{\mathbb{L}} \omega_2(\lambda) e^{i\lambda x-\omega_1(\lambda)t} \frac{1}{i\lambda} (u_0^{(1)})\hat{}(\lambda) \frac{d\lambda}{2i\rho(\lambda)}$$
$$= \lim_{A\to\infty} \left( \int_{-A}^{-1} + \int_{1}^{A} \right) \omega_2(\lambda) e^{i\lambda x-\omega_1(\lambda)t} \hat{u}_0(\lambda) \frac{d\lambda}{2i\rho(\lambda)}.$$

The point here is that, although the integral in the LHS of (4.1) does not converge absolutely (in general), it exists in the generalized sense, since

$$\int_{\kappa^+} \omega_2(\lambda) e^{i\lambda x-\omega_1(\lambda)t} \frac{u_0(0)}{i\lambda} \frac{d\lambda}{2i\rho(\lambda)} = \lim_{A\to\infty} \left( \int_{-A}^{-1} + \int_{1}^{A} \right) \omega_2(\lambda) e^{i\lambda x-\omega_1(\lambda)t} \frac{u_0(0)}{i\lambda} \frac{d\lambda}{2i\rho(\lambda)}$$

and the integral

$$\int_{\mathbb{L}} \omega_2(\lambda) e^{i\lambda x - \omega_1(\lambda)t} \frac{1}{i\lambda} (u_0^{(1)}\hat{)}(\lambda) \frac{d\lambda}{2i\rho(\lambda)} \text{ converges absolutely.}$$

Writing an equation analogous to (4.1) also for the integral $\int_{\mathbb{L}} \omega_1(\lambda) e^{i\lambda x-\omega_2(\lambda)t} \hat{u}_0(\lambda) d\lambda / 2i\rho(\lambda)$ and working as before, we obtain that the integral

$$\text{(4.2)} \quad \int_{-\infty}^{\infty} [\omega_1(\lambda) e^{i\lambda x-\omega_2(\lambda)t} - \omega_2(\lambda) e^{i\lambda x-\omega_1(\lambda)t}] \hat{u}_0(\lambda) \frac{d\lambda}{2i\rho(\lambda)}$$

$$= \lim_{\mathrm{A}\to\infty} \int_{-\mathrm{A}}^{\mathrm{A}} [\omega_1(\lambda) e^{i\lambda x-\omega_2(\lambda)t} - \omega_2(\lambda) e^{i\lambda x-\omega_1(\lambda)t}] \hat{u}_0(\lambda) \frac{d\lambda}{2i\rho(\lambda)}$$

exists, also in view of the fact that $\mathcal{K}_0(\lambda,t)$ is $C^\infty$ for $\lambda \in (-\infty,\infty)$. (See Lemma 1.) In fact, the above equation may be considered as the definition of the integral $\int_{-\infty}^{\infty}$ in (4.2), for $x > t > 0$. Also, we point out that this definition is in agreement with (3.19) and (3.28).

Conclusions like the above can be drawn for ***all*** the integrals in (1.5) in both cases: either $x > t > 0$ or $0 < x < t$, which may not converge absolutely.

***Part 2*** Next we prove that the function $u(x,t)$ is $C^\infty$ in $Q - \{x = t\}$. Indeed, in view of the previous observations and explanations, all the equations (3.19) – (3.31) as well as (3.1), (3.2) and their analogues, hold true, and since the RHSs of these equations define $C^\infty$ functions in $Q - \{x = t\}$, the required smoothness of $u(x,t)$ follows. The main observation here is that the integrals in the RHSs of the equations (3.19) – (3.31) are taken either on the finite contours $\kappa^+$ and $-\kappa^-$, or on the infinite contour $\mathbb{L}$, in which case differentiation under the integral sign is allowed, provided that the positive integers $N$ and $M$, involved in these integrals, are chosen sufficiently large.

***Part 3*** Now we prove that $u_{tt} + u_t - u_{xx} = f$ in $Q - \{x = t\}$. Keeping in mind that

$$\text{(4.3)} \quad \left(\frac{\partial^2}{\partial t^2} + \frac{\partial}{\partial t} - \frac{\partial^2}{\partial x^2}\right) e^{i\lambda x - \omega_j(\lambda)t} = 0$$

and in view of (3.20), (3.21), (3.24), (3.25), (3.28) and (3.29), it is easy to see that it suffices to show that

$$\text{(4.4)} \quad \left(\frac{\partial^2}{\partial t^2} + \frac{\partial}{\partial t} - \frac{\partial^2}{\partial x^2}\right) \mathcal{G}(g_0)(x,t) = 0,$$

$$\text{(4.5)} \quad \left(\frac{\partial^2}{\partial t^2} + \frac{\partial}{\partial t} - \frac{\partial^2}{\partial x^2}\right) \mathcal{F}^+(f)(x,t) = f(x,t)$$

and

$$\text{(4.6)} \quad \left(\frac{\partial^2}{\partial t^2} + \frac{\partial}{\partial t} - \frac{\partial^2}{\partial x^2}\right) \mathcal{F}^-(f)(x,t) = 0.$$

*Proof of (4.4) in the case* $x > t > 0$*:* This is trivial in view of (3.32).

*Proof of (4.4) in the case* $0 < x < t$*:* In view of (3.11) and the fact that

$$\left(\int_{\mathbb{L}} + \int_{\kappa^-}\right)\left[\frac{e^{i\lambda x-\omega_1(\lambda)t}}{[\omega_1(\lambda)]^j} \frac{\lambda d\lambda}{\rho(\lambda)}\right] = 0, \text{ for } j \geq 1,$$

we have

$$\text{(4.7)} \quad \mathcal{G}(g_0)(x,t) = \left\{\sum_{j=1}^{3} (-1)^{j-1} g_0^{(j-1)}(t) \int_{\kappa^+ + \kappa^-} \left[\frac{1}{[\omega_1(\lambda)]^j} - \frac{1}{[\omega_2(\lambda)]^j}\right] \frac{\lambda e^{i\lambda x} d\lambda}{\rho(\lambda)}\right\}$$

$$+ \left\{\sum_{j=1}^{3} (-1)^{j-1} g_0^{(j-1)}(0) \int_{\kappa^+ + \kappa^-} \frac{e^{i\lambda x - \omega_2(\lambda)t}}{[\omega_2(\lambda)]^j} \frac{\lambda d\lambda}{\rho(\lambda)}\right\}$$

$$- \left\{\left(\int_{\mathbb{L}} + \int_{\kappa^-}\right)\left(\frac{1}{[\omega_1(\lambda)]^3} e^{i\lambda x - \omega_1(\lambda)t} (g_0^{(3)}\tilde{)}(\omega_1(\lambda),t) - \frac{1}{[\omega_2(\lambda)]^3} e^{i\lambda x - \omega_2(\lambda)t} (g_0^{(3)}\tilde{)}(\omega_2(\lambda),t)\right) \frac{\lambda d\lambda}{\rho(\lambda)}\right\}.$$

In view of (4.3), the middle term in the RHS of (4.7) is annihilated by $\dfrac{\partial^2}{\partial t^2}+\dfrac{\partial}{\partial t}-\dfrac{\partial^2}{\partial x^2}$ and

$$(4.8)\quad \left(\frac{\partial^2}{\partial t^2}+\frac{\partial}{\partial t}-\frac{\partial^2}{\partial x^2}\right)\mathcal{G}(g_0)(x,t)=\left(\frac{\partial^2}{\partial t^2}+\frac{\partial}{\partial t}-\frac{\partial^2}{\partial x^2}\right)\left\{\sum_{j=1}^{3}(-1)^{j-1}g_0^{(j-1)}(t)\int_{\kappa^++\kappa^-}\left[\frac{1}{[\omega_1(\lambda)]^j}-\frac{1}{[\omega_2(\lambda)]^j}\right]\frac{\lambda e^{i\lambda x}d\lambda}{\rho(\lambda)}\right\}$$
$$-\left\{\left(\int_{\mathbb{L}}+\int_{\kappa^-}\right)\left(\frac{1}{[\omega_1(\lambda)]^3}e^{i\lambda x-\omega_1(\lambda)t}\frac{\partial^2}{\partial t^2}[(g_0^{(3)})\tilde{}(\omega_1(\lambda),t)]-\frac{1}{[\omega_2(\lambda)]^3}e^{i\lambda x-\omega_2(\lambda)t}\frac{\partial^2}{\partial t^2}[(g_0^{(3)})\tilde{}(\omega_2(\lambda),t)]\frac{\lambda d\lambda}{\rho(\lambda)}\right)\right\}$$
$$-\left\{\left(\int_{\mathbb{L}}+\int_{\kappa^-}\right)\left(\frac{1}{[\omega_1(\lambda)]^3}\frac{\partial}{\partial t}[e^{i\lambda x-\omega_1(\lambda)t}]\frac{\partial}{\partial t}[(g_0^{(3)})\tilde{}(\omega_1(\lambda),t)]-\frac{1}{[\omega_2(\lambda)]^3}\frac{\partial}{\partial t}[e^{i\lambda x-\omega_2(\lambda)t}]\frac{\partial}{\partial t}[(g_0^{(3)})\tilde{}(\omega_2(\lambda),t)]\frac{\lambda d\lambda}{\rho(\lambda)}\right)\right\}$$
$$-\left\{\left(\int_{\mathbb{L}}+\int_{\kappa^-}\right)\left(\frac{1}{[\omega_1(\lambda)]^3}e^{i\lambda x-\omega_1(\lambda)t}\frac{\partial}{\partial t}[(g_0^{(3)})\tilde{}(\omega_1(\lambda),t)]-\frac{1}{[\omega_2(\lambda)]^3}e^{i\lambda x-\omega_2(\lambda)t}\frac{\partial}{\partial t}[(g_0^{(3)})\tilde{}(\omega_2(\lambda),t)]\frac{\lambda d\lambda}{\rho(\lambda)}\right)\right\}.$$

Carrying out the differentiations in the RHS of (4.8), we find that

$$(4.9)\quad \text{LHS of (4.8)}=\left\{\sum_{j=1}^{3}(-1)^{j-1}g_0^{(j-1)}(t)\int_{\kappa^++\kappa^-}\left[\frac{1}{[\omega_1(\lambda)]^j}-\frac{1}{[\omega_2(\lambda)]^j}\right]\frac{\lambda^3 e^{i\lambda x}d\lambda}{\rho(\lambda)}\right\}$$
$$+\left\{\sum_{j=1}^{3}(-1)^{j-1}[g_0^{(j+1)}(t)+g_0^{(j)}(t)]\int_{\kappa^++\kappa^-}\left[\frac{1}{[\omega_1(\lambda)]^j}-\frac{1}{[\omega_2(\lambda)]^j}\right]\frac{\lambda e^{i\lambda x}d\lambda}{\rho(\lambda)}\right\}$$
$$-\left\{\left(\int_{\mathbb{L}}+\int_{\kappa^-}\right)\left[\frac{1}{[\omega_1(\lambda)]^3}e^{i\lambda x}g_0^{(4)}(t)-\frac{1}{[\omega_2(\lambda)]^3}e^{i\lambda x}g_0^{(4)}(t)\right]\frac{\lambda d\lambda}{\rho(\lambda)}\right\}$$
$$+\left\{\left(\int_{\mathbb{L}}+\int_{\kappa^-}\right)\left[\frac{1}{[\omega_1(\lambda)]^2}e^{i\lambda x}g_0^{(3)}(t)-\frac{1}{[\omega_2(\lambda)]^2}e^{i\lambda x}g_0^{(3)}(t)\right]\frac{\lambda d\lambda}{\rho(\lambda)}\right\}$$
$$-\left\{\left(\int_{\mathbb{L}}+\int_{\kappa^-}\right)\left[\frac{1}{[\omega_1(\lambda)]^3}e^{i\lambda x}g_0^{(3)}(t)-\frac{1}{[\omega_2(\lambda)]^3}e^{i\lambda x}g_0^{(3)}(t)\right]\frac{\lambda d\lambda}{\rho(\lambda)}\right\}.$$

Now, the contour $\mathbb{L}$, in the RHS of (4.9), can be replaced by $\kappa^+$. Thus, replacing $\mathbb{L}$ by $\kappa^+$ and taking into consideration the identities

$$\left(\frac{1}{\omega_1(\lambda)}-\frac{1}{\omega_2(\lambda)}\right)\frac{1}{\rho(\lambda)}=-\frac{1}{\lambda^2},\ \left(\frac{1}{[\omega_1(\lambda)]^2}-\frac{1}{[\omega_2(\lambda)]^2}\right)\frac{1}{\rho(\lambda)}=-\frac{1}{\lambda^4},\ \left(\frac{1}{[\omega_1(\lambda)]^3}-\frac{1}{[\omega_2(\lambda)]^3}\right)\frac{1}{\rho(\lambda)}=-\frac{1-\lambda^2}{\lambda^4},$$

it is easy to see that the RHS of (4.9) vanishes.

Therefore, (4.4), in the case $0<x<t$, follows from (4.8) and (4.9).

*Proof of (4.5):* Using the identity

$$e^{-\omega_j(\lambda)t}\tilde{\hat{f}}(\lambda,\omega_j(\lambda),t)=e^{-\omega_j(\lambda)t}\int_0^t e^{\omega_j(\lambda)\tau}\hat{f}(\lambda,\tau)d\tau$$
$$=\frac{\hat{f}(\lambda,t)}{\omega_j(\lambda)}-e^{-\omega_j(\lambda)t}\frac{\hat{f}(\lambda,0)}{\omega_j(\lambda)}-\frac{1}{\omega_j(\lambda)}e^{-\omega_j(\lambda)t}\int_0^t e^{\omega_j(\lambda)\tau}(f^{(1)})\hat{}(\lambda,\tau)d\tau\qquad(\omega_j(\lambda)\neq 0,\ \operatorname{Im}\lambda\leq 0)$$

and, setting $\gamma^-=(-\infty,-1]+\kappa^-+[1,\infty)$ so that $\int_{-\infty}^{-1}+\int_1^{\infty}+\int_{\kappa^-}=\int_{\gamma^-}$, we have

$$(4.7)\quad \mathcal{F}^+(f)(x,t)=-\int_{-\infty}^{\infty}\left[e^{i\lambda x-\omega_1(\lambda)t}\tilde{\hat{f}}(\lambda,\omega_1(\lambda),t)-e^{i\lambda x-\omega_2(\lambda)t}\tilde{\hat{f}}(\lambda,\omega_2(\lambda),t)\right]\frac{d\lambda}{2i\rho(\lambda)}$$
$$=-\int_{\gamma^-}\left[\frac{1}{\omega_1(\lambda)}-\frac{1}{\omega_2(\lambda)}\right]e^{i\lambda x}\hat{f}(\lambda,t)\frac{d\lambda}{2i\rho(\lambda)}+\int_{\gamma^-}\left[\frac{e^{i\lambda x-\omega_1(\lambda)t}}{\omega_1(\lambda)}-\frac{e^{i\lambda x-\omega_2(\lambda)t}}{\omega_2(\lambda)}\right]\hat{f}(\lambda,0)\frac{d\lambda}{2i\rho(\lambda)}$$

$$+\int_{\gamma^-}\left[e^{i\lambda x-\omega_1(\lambda)t}\frac{1}{\omega_1(\lambda)}\int_0^t e^{\omega_1(\lambda)\tau}(f^{(1)}\hat{)}(\lambda,\tau)d\tau-e^{i\lambda x-\omega_2(\lambda)t}\frac{1}{\omega_2(\lambda)}\int_0^t e^{\omega_2(\lambda)\tau}(f^{(1)}\hat{)}(\lambda,\tau)d\tau\right]\frac{d\lambda}{2i\rho(\lambda)}.$$

Kepping in mind the interpretation of the various integrals and working with $x\neq t$, $(x,t)\in Q$, we have

$$(4.8)\quad \frac{\partial^2}{\partial x^2}\int_{\gamma^-}\left[\frac{1}{\omega_1(\lambda)}-\frac{1}{\omega_2(\lambda)}\right]e^{i\lambda x}\hat{f}(\lambda,t)\frac{d\lambda}{2i\rho(\lambda)}$$

$$=-\frac{\partial^2}{\partial x^2}\left\{\int_{\kappa^++\kappa^-}\frac{1}{\lambda^2}e^{i\lambda x}\frac{f(0,t)}{i\lambda}d\lambda+\int_{\gamma^-}\left[\frac{1}{\lambda^2}e^{i\lambda x}\frac{(f^{(1)}\hat{)}(\lambda,t)}{i\lambda}d\lambda\right]\right\}$$

$$=\int_{\kappa^++\kappa^-}e^{i\lambda x}\frac{f(0,t)}{i\lambda}d\lambda+\int_{\gamma^-}\left[e^{i\lambda x}\frac{(f^{(1)}\hat{)}(\lambda,t)}{i\lambda}d\lambda\right]$$

$$=\int_{\gamma^-}e^{i\lambda x}\frac{f(0,t)}{i\lambda}d\lambda+\int_{\gamma^-}\left[e^{i\lambda x}\frac{(f^{(1)}\hat{)}(\lambda,t)}{i\lambda}d\lambda\right]=\int_{-\infty}^{\infty}e^{i\lambda x}\hat{f}(\lambda,t)d\lambda=2\pi f(x,t)$$

and

$$(4.9)\quad \left(\frac{\partial^2}{\partial t^2}+\frac{\partial}{\partial t}\right)\int_{\gamma^-}\left[\frac{1}{\omega_1(\lambda)}-\frac{1}{\omega_2(\lambda)}\right]e^{i\lambda x}\hat{f}(\lambda,t)\frac{d\lambda}{2i\rho(\lambda)}=\int_{\gamma^-}\left[\frac{1}{\omega_1(\lambda)}-\frac{1}{\omega_2(\lambda)}\right]e^{i\lambda x}[(f^{(1)}\hat{)}(\lambda,t)+(f^{(2)}\hat{)}(\lambda,t)]\frac{d\lambda}{2i\rho(\lambda)}.$$

Similarly,

$$(4.10)\quad \left(\frac{\partial^2}{\partial t^2}+\frac{\partial}{\partial t}-\frac{\partial^2}{\partial x^2}\right)\int_{\gamma^-}\left[\frac{e^{i\lambda x-\omega_1(\lambda)t}}{\omega_1(\lambda)}-\frac{e^{i\lambda x-\omega_2(\lambda)t}}{\omega_2(\lambda)}\right]\hat{f}(\lambda,0)\frac{d\lambda}{2i\rho(\lambda)}=0,$$

$$(4.11)\quad \int_{\gamma^-}\left\{e^{i\lambda x-\omega_1(\lambda)t}\frac{1}{\omega_1(\lambda)}\frac{\partial}{\partial t}\left[\int_0^t e^{\omega_1(\lambda)\tau}(f^{(1)}\hat{)}(\lambda,\tau)d\tau\right]-e^{i\lambda x-\omega_2(\lambda)t}\frac{1}{\omega_2(\lambda)}\frac{\partial}{\partial t}\left[\int_0^t e^{\omega_2(\lambda)\tau}(f^{(1)}\hat{)}(\lambda,\tau)d\tau\right]\right\}\frac{d\lambda}{2i\rho(\lambda)}$$

$$=\int_{\gamma^-}\left[\frac{1}{\omega_1(\lambda)}-\frac{1}{\omega_2(\lambda)}\right]e^{i\lambda x}(f^{(1)}\hat{)}(\lambda,t)\frac{d\lambda}{2i\rho(\lambda)},$$

$$(4.12)\quad \int_{\gamma^-}\left\{\frac{\partial}{\partial t}\left[e^{i\lambda x-\omega_1(\lambda)t}\frac{1}{\omega_1(\lambda)}\right]\frac{\partial}{\partial t}\left[\int_0^t e^{\omega_1(\lambda)\tau}(f^{(1)}\hat{)}(\lambda,\tau)d\tau\right]\right.$$

$$\left.-\frac{\partial}{\partial t}\left[e^{i\lambda x-\omega_2(\lambda)t}\frac{1}{\omega_2(\lambda)}\frac{\partial}{\partial t}\right]\left[\int_0^t e^{\omega_2(\lambda)\tau}(f^{(1)}\hat{)}(\lambda,\tau)d\tau\right]\right\}\frac{d\lambda}{2i\rho(\lambda)}=0$$

and

$$(4.13)\quad \int_{\gamma^-}\left\{e^{i\lambda x-\omega_1(\lambda)t}\frac{1}{\omega_1(\lambda)}\frac{\partial^2}{\partial t^2}\left[\int_0^t e^{\omega_1(\lambda)\tau}(f^{(1)}\hat{)}(\lambda,\tau)d\tau\right]-e^{i\lambda x-\omega_2(\lambda)t}\frac{1}{\omega_2(\lambda)}\frac{\partial^2}{\partial t^2}\left[\int_0^t e^{\omega_2(\lambda)\tau}(f^{(1)}\hat{)}(\lambda,\tau)d\tau\right]\right\}\frac{d\lambda}{2i\rho(\lambda)}$$

$$=\int_{\gamma^-}\left[\frac{1}{\omega_1(\lambda)}-\frac{1}{\omega_2(\lambda)}\right]e^{i\lambda x}(f^{(2)}\hat{)}(\lambda,t)\frac{d\lambda}{2i\rho(\lambda)}.$$

Now (4.5) follows from (4.7) – (4.13).

The proof of (4.6) is similar, the only difference being the equation

$$\frac{\partial^2}{\partial x^2}\int_{\gamma^+}\left[\frac{1}{\omega_1(\lambda)}-\frac{1}{\omega_2(\lambda)}\right]e^{i\lambda x}\hat{f}(-\lambda,t)\frac{d\lambda}{2i\rho(\lambda)}=\int_{-\infty}^{\infty}e^{i\lambda x}\hat{f}(-\lambda,t)d\lambda=0,$$

where $\gamma^+:=(-\infty,-1]-\kappa^++[1,\infty)$.

***Part 4*** We will prove the $C^\infty$ – extension of $u(x,t)|_{\{x<t\}}$ and $u(x,t)|_{\{x>t\}}$ to $Q\cap\{x\le t\}$ and $Q\cap\{x\ge t\}$, respectively. This follows essentially from the interpretations (3.19) – (3.31) of the various integrals. Let us give an example. In the case $x>t>0$, (3.19) gives:

$$(4.14)\quad \frac{\partial^{n+m}}{\partial x^n \partial t^m}\left(\int_{\mathbb{L}} \omega_2(\lambda) e^{i\lambda x-\omega_1(\lambda)t}\hat{u}_0(\lambda)\frac{d\lambda}{2i\rho(\lambda)}\right)$$

$$= \sum_{j=1}^{N} u_0^{(j-1)}(0) \int_{\kappa^+} (i\lambda)^n [-\omega_1(\lambda)]^m \omega_2(\lambda) e^{i\lambda x-\omega_1(\lambda)t}\frac{1}{(i\lambda)^j}\frac{d\lambda}{2i\rho(\lambda)}$$

$$+ \int_{\mathbb{L}} (i\lambda)^n [-\omega_1(\lambda)]^m \omega_2(\lambda) e^{i\lambda x-\omega_1(\lambda)t}\frac{1}{(i\lambda)^N}(u_0^{(N)})\hat{}(\lambda)\frac{d\lambda}{2i\rho(\lambda)},$$

provided that $N$ is chosen sufficiently large. Now it is clear that the function, defined (in $\{x>t\}$) by the RHS of (4.14), extends continuously to $\{x \geq t\}$.

Similar formulas can be written for all the integrals in (1.5), in both cases: $x>t>0$ and $0<x<t$, and this implies the conclusion of *Part 4*.

The proof is complete.

## 5. Proof of Theorem 2

***Step 1*** We claim that

$$(5.1)\qquad \lim_{t\to 0^+} \mathcal{U}_0^+(u_0)(x,t) = u_0(x).$$

Since, for $t<x$,

$$(5.2)\quad \mathcal{U}_0^+(u_0)(x,t) = \int_{-1}^{1} [\omega_1(\lambda)e^{i\lambda x-\omega_2(\lambda)t} - \omega_2(\lambda)e^{i\lambda x-\omega_1(\lambda)t}]\hat{u}_0(\lambda)\frac{d\lambda}{2i\rho(\lambda)}$$

$$+\int_{\mathbb{L}} [\omega_1(\lambda)e^{i\lambda x-\omega_2(\lambda)t} - \omega_2(\lambda)e^{i\lambda x-\omega_1(\lambda)t}]\left[\frac{u_0(0)}{i\lambda}+\frac{(u_0^{(1)})\hat{}(\lambda)}{i\lambda}\right]\frac{d\lambda}{2i\rho(\lambda)},$$

we see that, indeed,

$$\lim_{t\to 0^+} \mathcal{U}_0^+(u_0)(x,t) = \int_{-1}^{1} [\omega_1(\lambda)e^{i\lambda x} - \omega_2(\lambda)e^{i\lambda x}]\hat{u}_0(\lambda)\frac{d\lambda}{2i\rho(\lambda)}$$

$$+\lim_{t\to 0^+}\int_{\kappa^+} [\omega_1(\lambda)e^{i\lambda x-\omega_2(\lambda)t} - \omega_2(\lambda)e^{i\lambda x-\omega_1(\lambda)t}]\frac{u_0(0)}{i\lambda}\frac{d\lambda}{2i\rho(\lambda)}$$

$$+\lim_{t\to 0^+}\int_{\mathbb{L}} [\omega_1(\lambda)e^{i\lambda x-\omega_2(\lambda)t} - \omega_2(\lambda)e^{i\lambda x-\omega_1(\lambda)t}]\frac{(u_0^{(1)})\hat{}(\lambda)}{i\lambda}\frac{d\lambda}{2i\rho(\lambda)}$$

$$=\int_{-1}^{1} e^{i\lambda x}\hat{u}_0(\lambda)d\lambda + \int_{\kappa^+} [\omega_1(\lambda)e^{i\lambda x} - \omega_2(\lambda)e^{i\lambda x}]\frac{u_0(0)}{i\lambda}\frac{d\lambda}{2i\rho(\lambda)} + \int_{\mathbb{L}} [\omega_1(\lambda)e^{i\lambda x} - \omega_2(\lambda)e^{i\lambda x}]\frac{(u_0^{(1)})\hat{}(\lambda)}{i\lambda}\frac{d\lambda}{2i\rho(\lambda)}$$

$$=\int_{-1}^{1} e^{i\lambda x}\hat{u}_0(\lambda)d\lambda + \int_{\kappa^+} e^{i\lambda x}\frac{u_0(0)}{i\lambda}d\lambda + \int_{\mathbb{L}} e^{i\lambda x}\frac{(u_0^{(1)})\hat{}(\lambda)}{i\lambda}d\lambda$$

$$=\int_{-1}^{1} e^{i\lambda x}\hat{u}_0(\lambda)d\lambda + \int_{\mathbb{L}} e^{i\lambda x}\frac{u_0(0)}{i\lambda}d\lambda + \int_{\mathbb{L}} e^{i\lambda x}\frac{(u_0^{(1)})\hat{}(\lambda)}{i\lambda}d\lambda = \int_{-\infty}^{\infty} e^{i\lambda x}\hat{u}_0(\lambda)d\lambda = 2\pi u_0(x).$$

(In the above computation, it is crucial that the contours $[-1,1]$ and $\kappa^+$ are finite, and that the integrand of the integral taken on $\mathbb{L}$ and contains the term $(u_0^{(1)})\hat{}(\lambda)$ is $O(1/\lambda^2)$. It is because of these facts that the interchange of the integral with the limit $\lim_{t\to 0^+}$ is justified.)

***Step 2*** We claim that

$$(5.3)\qquad \lim_{t\to 0^+} \mathcal{U}_0^-(u_0)(x,t) = 0.$$

Indeed, working as in the previous step, we obtain that

$$\lim_{t\to 0^+} \mathcal{U}_0^-(u_0)(x,t) = \int_{-\infty}^{\infty} e^{i\lambda x}\hat{u}_0(-\lambda)d\lambda = 0,$$

where the last equation follows from the Cauchy-Jordan lemma.

***Step 3*** It is easy to see that

$$\lim_{t\to 0^+}\mathcal{U}_1^+(u_1)=0 \text{ and } \lim_{t\to 0^+}\mathcal{U}_1^-(u_1)=0 .$$

***Step 4*** We claim that

$$\lim_{t\to 0^+}\frac{\partial}{\partial t}[\mathcal{U}_0^+(u_0)(x,t)]=0 .$$

Indeed, writing

$$\mathcal{U}_0^+(u_0)(x,t)=\int_{-1}^{1}[\omega_1(\lambda)e^{i\lambda x-\omega_2(\lambda)t}-\omega_2(\lambda)e^{i\lambda x-\omega_1(\lambda)t}]\hat{u}_0(\lambda)\frac{d\lambda}{2i\rho(\lambda)}$$

$$+\int_{\mathbb{L}}[\omega_1(\lambda)e^{i\lambda x-\omega_2(\lambda)t}-\omega_2(\lambda)e^{i\lambda x-\omega_1(\lambda)t}]\left[\frac{u_0(0)}{i\lambda}+\frac{u_0^{(1)}(0)}{(i\lambda)^2}+\frac{(u_0^{(2)})\hat{}(\lambda)}{(i\lambda)^2}\right]\frac{d\lambda}{2i\rho(\lambda)}$$

$$=\int_{-1}^{1}[\omega_1(\lambda)e^{i\lambda x-\omega_2(\lambda)t}-\omega_2(\lambda)e^{i\lambda x-\omega_1(\lambda)t}]\hat{u}_0(\lambda)\frac{d\lambda}{2i\rho(\lambda)}$$

$$+\int_{\kappa^+}[\omega_1(\lambda)e^{i\lambda x-\omega_2(\lambda)t}-\omega_2(\lambda)e^{i\lambda x-\omega_1(\lambda)t}]\left[\frac{u_0(0)}{i\lambda}+\frac{u_0^{(1)}(0)}{(i\lambda)^2}\right]\frac{d\lambda}{2i\rho(\lambda)}$$

$$+\int_{\mathbb{L}}[\omega_1(\lambda)e^{i\lambda x-\omega_2(\lambda)t}-\omega_2(\lambda)e^{i\lambda x-\omega_1(\lambda)t}]\frac{(u_0^{(2)})\hat{}(\lambda)}{(i\lambda)^2}\frac{d\lambda}{2i\rho(\lambda)},$$

differentiating and letting $t\to 0^+$, we obtain

$$\lim_{t\to 0^+}\frac{\partial}{\partial t}[\mathcal{U}_0^+(u_0)(x,t)]=-\int_{-1}^{1}[\omega_1(\lambda)\omega_2(\lambda)e^{i\lambda x}-\omega_1(\lambda)\omega_2(\lambda)e^{i\lambda x}]\hat{u}_0(\lambda)\frac{d\lambda}{2i\rho(\lambda)}$$

$$-\int_{\kappa^+}[\omega_1(\lambda)\omega_2(\lambda)e^{i\lambda x}-\omega_1(\lambda)\omega_2(\lambda)e^{i\lambda x}]\left[\frac{u_0(0)}{i\lambda}+\frac{(u_0^{(1)})\hat{}(\lambda)}{i\lambda}\right]\frac{d\lambda}{2i\rho(\lambda)}$$

$$-\int_{\mathbb{L}}[\omega_1(\lambda)\omega_2(\lambda)e^{i\lambda x}-\omega_1(\lambda)\omega_2(\lambda)e^{i\lambda x}]\frac{(u_0^{(2)})\hat{}(\lambda)}{(i\lambda)^2}\frac{d\lambda}{2i\rho(\lambda)}=0 .$$

***Step 5*** Writing

$$\mathcal{U}_0^-(u_0)=\int_{\kappa^++[-1,1]}[\omega_2(\lambda)e^{i\lambda x-\omega_1(\lambda)t}-\omega_1(\lambda)e^{i\lambda x-\omega_2(\lambda)t}]\hat{u}_0(-\lambda)\frac{d\lambda}{2i\rho(\lambda)},$$

we easily obtain that $\lim_{t\to 0^+}\frac{\partial}{\partial t}[\mathcal{U}_0^-(u_0)(x,t)]=0$.

***Step 6*** We claim that

$$\lim_{t\to 0^+}\frac{\partial}{\partial t}[\mathcal{U}_1^+(u_1)(x,t)]=u_1(x) \text{ and } \lim_{t\to 0^+}\frac{\partial}{\partial t}[\mathcal{U}_1^-(u_1)(x,t)]=0 . \tag{5.4}$$

Indeed, the proof of (5.4) is similar to the proofs of (5.1) and (5.3).

*Proof of 1st assertion* In view of (3.32), it is trivial that

$$\lim_{t\to 0^+}\mathcal{G}(g_0)(x,t)=0 \text{ and } \lim_{t\to 0^+}\frac{\partial}{\partial t}[\mathcal{G}(g_0)(x,t)]=0 . \tag{5.5}$$

It is also clear that

$$\lim_{t\to 0^+}[\mathcal{F}^+(f)+\mathcal{F}^-(f)](x,t)=0 \text{ and } \lim_{t\to 0^+}\frac{\partial}{\partial t}[\mathcal{F}^+(f)+\mathcal{F}^-(f)](x,t)=0 , \tag{5.6}$$

since the integrands of both $\mathcal{F}^+(f)(x,t)$ and $\mathcal{F}^-(f)(x,t)$ are $O(1/\lambda^3)$, for $\lambda\to\pm\infty$ ($\lambda\in\mathbb{R}$). Thus, the 1st assertion follows from the results of *Steps 1- 6*, (5.5) and (5.6).

***Step 7*** We claim that

$$\lim_{x\to 0^+}\mathcal{G}(g_0)(x,t)=g_0(t) . \tag{5.7}$$

For $0<x<t$, we have

(5.8)
$$\mathcal{G}(g_0)(x,t) = \mathcal{G}_1(g_0)(x,t) + \mathcal{G}_2(g_0)(x,t),$$
where
$$\mathcal{G}_1(g_0)(x,t) = \int_{-1}^{1}[e^{i\lambda x-\omega_1(\lambda)t}\tilde{g}_0(\omega_1(\lambda),t) - e^{i\lambda x-\omega_2(\lambda)t}\tilde{g}_0(\omega_2(\lambda),t)]\frac{\lambda d\lambda}{\rho(\lambda)}$$
and
$$\mathcal{G}_2(g_0)(x,t) = \int_{\mathbb{L}}[e^{i\lambda x-\omega_1(\lambda)t}\tilde{g}_0(\omega_1(\lambda),t) - e^{i\lambda x-\omega_2(\lambda)t}\tilde{g}_0(\omega_2(\lambda),t)]\frac{\lambda d\lambda}{\rho(\lambda)}.$$
Now

(5.9)
$$\lim_{x\to 0^+}\mathcal{G}_1(g_0)(x,t) = \int_{-1}^{1}[e^{-\omega_1(\lambda)t}\tilde{g}_0(\omega_1(\lambda),t) - e^{-\omega_2(\lambda)t}\tilde{g}_0(\omega_2(\lambda),t)]\frac{\lambda d\lambda}{\rho(\lambda)} = 0,$$
where the last equation follows from the fact that the integrand in (5.9) is an odd function of $\lambda$.
On the other hand, we may write

(5.10)
$$\mathcal{G}_2(g_0)(x,t) = g_0(t)\int_{\kappa^+} e^{i\lambda x}\left[\frac{1}{\omega_1(\lambda)} - \frac{1}{\omega_2(\lambda)}\right]\frac{\lambda d\lambda}{\rho(\lambda)} - g_0(0)\left[\int_{-\kappa^-}\frac{e^{i\lambda x-\omega_1(\lambda)t}}{\omega_1(\lambda)}\frac{\lambda d\lambda}{\rho(\lambda)} - \int_{\kappa^+}\frac{e^{i\lambda x-\omega_2(\lambda)t}}{\omega_2(\lambda)}\frac{\lambda d\lambda}{\rho(\lambda)}\right]$$
$$-\int_{\mathbb{L}} e^{i\lambda x}\left[\frac{e^{-\omega_1(\lambda)t}}{\omega_1(\lambda)}(g_0^{(1)}\tilde{)}(\omega_1(\lambda),t) - \frac{e^{-\omega_2(\lambda)t}}{\omega_2(\lambda)}(g_0^{(1)}\tilde{)}(\omega_2(\lambda),t)\right]\frac{\lambda d\lambda}{\rho(\lambda)}.$$
But

(5.11)
$$\lim_{x\to 0^+}\int_{\kappa^+} e^{i\lambda x}\left[\frac{1}{\omega_1(\lambda)} - \frac{1}{\omega_2(\lambda)}\right]\frac{\lambda d\lambda}{\rho(\lambda)} = \int_{\kappa^+}\left[\frac{1}{\omega_1(\lambda)} - \frac{1}{\omega_2(\lambda)}\right]\frac{\lambda d\lambda}{\rho(\lambda)}$$
$$= \int_{\kappa^+}\frac{\omega_2(\lambda)-\omega_1(\lambda)}{\omega_1(\lambda)\omega_2(\lambda)}\frac{\lambda d\lambda}{\rho(\lambda)} = \int_{\kappa^+}\frac{-2i\rho(\lambda)}{\lambda^2}\frac{\lambda d\lambda}{\rho(\lambda)} = 2\pi.$$
Also,

(5.12)
$$\lim_{x\to 0^+}\left[\int_{-\kappa^-}\frac{e^{i\lambda x-\omega_1(\lambda)t}}{\omega_1(\lambda)}\frac{\lambda d\lambda}{\rho(\lambda)} - \int_{\kappa^+}\frac{e^{i\lambda x-\omega_2(\lambda)t}}{\omega_2(\lambda)}\frac{\lambda d\lambda}{\rho(\lambda)}\right] = \int_{-\kappa^-}\frac{e^{-\omega_1(\lambda)t}}{\omega_1(\lambda)}\frac{\lambda d\lambda}{\rho(\lambda)} - \int_{\kappa^+}\frac{e^{-\omega_2(\lambda)t}}{\omega_2(\lambda)}\frac{\lambda d\lambda}{\rho(\lambda)} = 0,$$
where the last equation follows by changing the variable in the integral over $\kappa^+$, setting $\mu = -\lambda$.
Similarly,

(5.13)
$$\lim_{x\to 0^+}\int_{\mathbb{L}} e^{i\lambda x}\left[\frac{e^{-\omega_1(\lambda)t}}{\omega_1(\lambda)}(g_0^{(1)}\tilde{)}(\omega_1(\lambda),t) - \frac{e^{-\omega_2(\lambda)t}}{\omega_2(\lambda)}(g_0^{(1)}\tilde{)}(\omega_2(\lambda),t)\right]\frac{\lambda d\lambda}{\rho(\lambda)} = 0,$$
since, the integrand in (5.13) becomes an odd function of $\lambda$, ***when*** $x = 0$.
Thus, (5.7) follows from (5.8) – (5.13).

***Step 8*** We claim that

(5.14)
$$\lim_{x\to 0^+}[\mathcal{U}_0^+(u_0) + \mathcal{U}_0^-(u_0)](x,t) = 0.$$

Indeed, since (5.2) holds also for $0 < x < t$, we obtain

(5.15)
$$\lim_{x\to 0^+}\mathcal{U}_0^+(u_0)(x,t) = \int_{-1}^{1}[\omega_1(\lambda)e^{-\omega_2(\lambda)t} - \omega_2(\lambda)e^{-\omega_1(\lambda)t}]\hat{u}_0(\lambda)\frac{d\lambda}{2i\rho(\lambda)}$$
$$+\int_{\kappa^+}\omega_1(\lambda)e^{-\omega_2(\lambda)t}\frac{u_0(0)}{i\lambda}\frac{d\lambda}{2i\rho(\lambda)} + \int_{\mathbb{L}}\omega_1(\lambda)e^{-\omega_2(\lambda)t}\frac{(u_0^{(1)}\hat{)}(\lambda)}{i\lambda}\frac{d\lambda}{2i\rho(\lambda)} + \int_{-\kappa^-}\omega_2(\lambda)e^{-\omega_1(\lambda)t}\hat{u}_0(\lambda)\frac{d\lambda}{2i\rho(\lambda)}$$
$$= \int_{-1}^{1}[\omega_1(\lambda)e^{-\omega_2(\lambda)t} - \omega_2(\lambda)e^{-\omega_1(\lambda)t}]\hat{u}_0(\lambda)\frac{d\lambda}{2i\rho(\lambda)} + \int_{\mathbb{L}}\omega_1(\lambda)e^{-\omega_2(\lambda)t}\frac{u_0(0)}{i\lambda}\frac{d\lambda}{2i\rho(\lambda)}$$
$$+\int_{\mathbb{L}}\omega_1(\lambda)e^{-\omega_2(\lambda)t}\frac{(u_0^{(1)}\hat{)}(\lambda)}{i\lambda}\frac{d\lambda}{2i\rho(\lambda)} + \int_{\mathbb{L}}\omega_2(\lambda)e^{-\omega_1(\lambda)t}\hat{u}_0(\lambda)\frac{d\lambda}{2i\rho(\lambda)}$$
$$= \int_{-\infty}^{\infty}[\omega_1(\lambda)e^{-\omega_2(\lambda)t} - \omega_2(\lambda)e^{-\omega_1(\lambda)t}]\hat{u}_0(\lambda)\frac{d\lambda}{2i\rho(\lambda)}.$$

(We had to use (5.2), since certain of the above integrals, in general, do not converge absolutely – for example the last one.)
Similarly, we show that

$$\lim_{x\to 0^+} \mathcal{U}_0^-(u_0)(x,t) = \int_{-\infty}^{\infty} [\omega_2(\lambda)e^{-\omega_1(\lambda)t} - \omega_1(\lambda)e^{-\omega_2(\lambda)t}]\hat{u}_0(-\lambda)\frac{d\lambda}{2i\rho(\lambda)}. \tag{5.16}$$

Now, changing the variable in the integral of (5.16), by setting $\mu=-\lambda$, and taking into consideration (5.15), we obtain (5.14).

*Proof of 2nd assertion* Since the integrals $\mathcal{U}_1^+(u_0)(x,t)$ and $\mathcal{U}_1^-(u_0)(x,t)$ converge absolutely, it is easy to see, by setting $\mu=-\lambda$ in the second one, that

$$\lim_{x\to 0^+} [\mathcal{U}_1^+(u_0)+\mathcal{U}_1^-(u_0)](x,t) = 0. \tag{5.17}$$

Similarly,

$$\lim_{x\to 0^+} [\mathcal{F}^+(f)+\mathcal{F}^-(f)](x,t)] = 0. \tag{5.18}$$

Thus, the 2nd conclusion follows from (5.7), (5.14), (5.17) and (5.18).

*Proof of 3rd assertion* Using the analysis of $\mathcal{U}_0^+(u_0)(x,t)$, given by (5.2), the analogous analysis of $\mathcal{U}_0^-(u_0)(x,t)$ and the analysis $\mathcal{G}_2(g_0)(x,t)$, given by (5.10), we see that the existence of the limit

$$\lim_{\substack{(x,t)\to(p,p)\\ x\neq t}} [\mathcal{U}_0^+(u_0)+\mathcal{U}_0^-(u_0)+\mathcal{G}(g_0)](x,t), \text{ for } p>0,$$

is equivalent to the existence of the limit

$$\lim_{\substack{(x,t)\to(p,p)\\ x\neq t}} \left[ u_0(0)\int_{\mathbb{L}} \frac{\omega_1(\lambda)e^{i\lambda x-\omega_2(\lambda)t}}{i\lambda}\frac{d\lambda}{i\rho(\lambda)} + g_0(0)\int_{\mathbb{L}} \frac{e^{i\lambda x-\omega_2(\lambda)t}}{\omega_2(\lambda)}\frac{\lambda d\lambda}{\rho(\lambda)} \right.$$
$$\left. - u_0(0)\int_{\mathbb{L}} \frac{\omega_2(\lambda)e^{i\lambda x-\omega_1(\lambda)t}}{i\lambda}\frac{d\lambda}{i\rho(\lambda)} - g_0(0)\int_{\mathbb{L}} \frac{e^{i\lambda x-\omega_1(\lambda)t}}{\omega_1(\lambda)}\frac{\lambda d\lambda}{\rho(\lambda)} \right]. \tag{5.19}$$

Now it is easy to see that the assumption "$u_0(0)=g_0(0)$" and the equality $\omega_1(\lambda)\omega_2(\lambda)=\lambda^2$ imply that the quantity inside the square brackets (of the quantity (5.19)), is zero. Thus, considering the analysis of $[\mathcal{U}_0^+(u_0)+\mathcal{U}_0^-(u_0)+\mathcal{G}(g_0)](x,t)$, which we used above, we easily see that (1.5) is reduced to formula (1.7) and this completes the proof of 3rd assertion.

*Proof of 4th assertion* To begin with, for every $p>0$, the limits

$$\lim_{\substack{(x,t)\to(p,p)\\ x\neq t}} \frac{\partial^{n+m}}{\partial x^n \partial x^m}[\mathcal{F}^{\pm}(f)(x,t)] \text{ exist, for } (n,m)\in\{(0,0),(1,0),(0,1)\}, \tag{5.20}$$

since the integrands of both $\mathcal{F}^+(f)(x,t)$ and $\mathcal{F}^-(f)(x,t)$ are $\mathrm{O}(1/\lambda^3)$, as $\lambda\to\pm\infty$ ($\lambda\in\mathbb{R}$).

Using (3.7) and (3.11), we obtain that

(5.21) *the integrand of* $[\mathcal{U}_0^+(u_0)+\mathcal{U}_0^-(u_0)](x,t) = -u_0(0)\frac{\omega_2(\lambda)}{i\lambda}\frac{1}{i\rho(\lambda)}e^{i\lambda x-\omega_1(\lambda)t}$

$$+u_0(0)\frac{\omega_1(\lambda)}{i\lambda}\frac{1}{i\rho(\lambda)}e^{i\lambda x-\omega_2(\lambda)t}+\mathrm{O}(1/\lambda^3),$$

(5.22) *the integrand of* $\mathcal{G}(g_0)(x,t) = \left\{-\frac{g_0(0)}{\omega_1(\lambda)}\frac{\lambda}{\rho(\lambda)}+\frac{g_0'(0)}{[\omega_1(\lambda)]^2}\frac{\lambda}{\rho(\lambda)}\right\}e^{i\lambda x-\omega_1(\lambda)t}$

$$-\left\{-\frac{g_0(0)}{\omega_2(\lambda)}\frac{\lambda}{\rho(\lambda)}+\frac{g_0'(0)}{[\omega_2(\lambda)]^2}\frac{\lambda}{\rho(\lambda)}\right\}e^{i\lambda x-\omega_2(\lambda)t}$$

$$+\left\{\frac{g_0(t)}{\omega_1(\lambda)}\frac{\lambda}{\rho(\lambda)}-\frac{g_0'(t)}{[\omega_1(\lambda)]^2}\frac{\lambda}{\rho(\lambda)}\right\}e^{i\lambda x}-\left\{\frac{g_0(t)}{\omega_2(\lambda)}\frac{\lambda}{\rho(\lambda)}-\frac{g_0'(t)}{[\omega_2(\lambda)]^2}\frac{\lambda}{\rho(\lambda)}\right\}e^{i\lambda x}+\mathrm{O}(1/\lambda^3)$$

and

(5.23) *the integrand of* $[\mathcal{U}_1^+(u_1)+\mathcal{U}_1^-(u_1)](x,t) = -\frac{u_1(0)}{i\lambda}\frac{1}{i\rho(\lambda)}e^{i\lambda x-\omega_1(\lambda)t}+\frac{u_1(0)}{i\lambda}\frac{1}{i\rho(\lambda)}e^{i\lambda x-\omega_2(\lambda)t}+\mathrm{O}(1/\lambda^3)$.

It follows from (5.21) – (5.23), in view of the assumption "$u_0(0)=g_0(0)$ & $u_1(0)=g_0'(0)$", that the

(5.24) *the integrand of* $[\mathcal{U}_0^+(u_0)+\mathcal{U}_0^-(u_0)+\mathcal{U}_1^+(u_1)+\mathcal{U}_1^-(u_1)+\mathcal{G}(g_0)](x,t)$

$$=\frac{u_1(0)e^{i\lambda x-\omega_1(\lambda)t}}{\lambda\omega_1(\lambda)\rho(\lambda)}-\frac{u_1(0)e^{i\lambda x-\omega_2(\lambda)t}}{\lambda\omega_2(\lambda)\rho(\lambda)}$$

$$+\,terms\ of\ the\ form\ (function\ of\ \lambda\ \&\ t)e^{i\lambda x}+\mathrm{O}(1/\lambda^3)=\mathrm{O}(1/\lambda^3),$$

where we used also the fact that

$$\frac{1}{(i\lambda)[i\rho(\lambda)]}-\frac{\lambda}{[\omega_j(\lambda)]^2\rho(\lambda)}=-\frac{\omega_j(\lambda)}{\lambda[\omega_j(\lambda)]^2\rho(\lambda)}=-\frac{1}{\lambda\omega_j(\lambda)\rho(\lambda)}.$$

(The two terms that we kept in the beginning of the RHS of (5.24) and which are $\mathrm{O}(1/\lambda^3)$, are needed in the proof of 5th assertion.)

It follows from (5.24) and the equation

$$\int_{\mathbb{L}}(function\ of\ \lambda\ \&\ t)e^{i\lambda x}d\lambda=\int_{\kappa^+}(function\ of\ \lambda\ \&\ t)e^{i\lambda x}d\lambda$$

(which holds for the *particular* functions of $\lambda,t$, appearing in (5.24)) that the limits

(5.25) $\lim\limits_{\substack{(x,t)\to(p,p)\\ x\neq t}}\frac{\partial^{n+m}}{\partial x^n\partial x^m}\{[\mathcal{U}_0^+(u_0)+\mathcal{U}_0^-(u_0)+\mathcal{U}_1^+(u_1)+\mathcal{U}_1^-(u_1)+\mathcal{G}(g_0)](x,t)\}$ exist, for $(n,m)\in\{(0,0),(1,0),(0,1)\}$.

Now assertion 4th follows from (5.20) and (5.25).

*Proof of 5th assertion* Using (3.7), (3.11) and (3.12), and working as previously, we find, in view of the assumption "$u_0(0)=g_0(0)$ & $u_1(0)=g_0'(0)$", that

*the integrand of* $[\mathcal{U}_0^+(u_0)+\mathcal{U}_0^-(u_0)+\mathcal{U}_1^+(u_1)+\mathcal{U}_1^-(u_1)+\mathcal{G}(g_0)+\mathcal{F}^+(f)+\mathcal{F}^-(f)](x,t)$

$$=\frac{1}{\rho(\lambda)}\left[-u_0''(0)\frac{\omega_2(\lambda)}{i(i\lambda)^3}+g_0'(0)\frac{1}{\lambda\omega_1(\lambda)}-g_0''(0)\frac{\lambda}{[\omega_1(\lambda)]^3}+f(0,0)\frac{1}{i(i\lambda)\omega_1(\lambda)}\right]e^{i\lambda x-\omega_1(\lambda)t}$$

$$-\frac{1}{\rho(\lambda)}\left[-u_0''(0)\frac{\omega_1(\lambda)}{i(i\lambda)^3}+g_0'(0)\frac{1}{\lambda\omega_2(\lambda)}-g_0''(0)\frac{\lambda}{[\omega_2(\lambda)]^3}-f(0,0)\frac{1}{i(-i\lambda)\omega_2(\lambda)}\right]e^{i\lambda x-\omega_2(\lambda)t}$$

$$+\,terms\ of\ the\ form\ (function\ of\ \lambda\ \&\ t)e^{i\lambda x}+\mathrm{O}(1/\lambda^4).$$

But each of the following quantities

$$-\frac{1}{\rho(\lambda)}\frac{\omega_j(\lambda)}{i(i\lambda)^3},\ -\frac{1}{\rho(\lambda)\lambda\omega_j(\lambda)},\ \frac{\lambda}{\rho(\lambda)[\omega_j(\lambda)]^3},\ \frac{1}{\rho(\lambda)i(i\lambda)\omega_j(\lambda)}\ (j=1,2)\text{ can be written as }(-1)^j\frac{i}{\lambda^3}+\mathrm{O}(1/\lambda^4)$$

and, therefore, the assumption "$g_0''(0)+g_0'(0)-u_0''(0)=f(0,0)$" implies the 5th conclusion.

This completes the proof of Theorem 2.

***Remarks*** **(1)** It is easy to write formulas analogous to (1.7) also in the cases of 4th and 5th parts of Theorem 2.

**(2)** Throughtout this paper, we are careful in every case where we deal with interchanging limits, in particular when we have to interchange differentiation with integration or a limit with integration. Indeed, the justification of such an interchange, whenever this is possible, is absolutely necessary. In our computations, this justification is clear or easy to give, although, sometimes, it may be implicit and not clearly stated.

To see what can happen if we do not pay the appropriate attention to such switchings, let us give the following example. The integral

$$\mathcal{G}(g_0)(x,t)=\int_{-\infty}^{\infty}[e^{i\lambda x-\omega_1(\lambda)t}\tilde{g}_0(\omega_1(\lambda),t)-e^{i\lambda x-\omega_2(\lambda)t}\tilde{g}_0(\omega_2(\lambda),t)]\frac{\lambda d\lambda}{\rho(\lambda)},\tag{5.26}$$

which exists as a generalized integral and appears in the solution (1.5), is the part of the solution which gives the boundary condition:

$$\lim_{x\to 0^+} u(x,t) = \lim_{x\to 0^+} \mathcal{G}(g_0)(x,t) = 2\pi g_0(t).$$

However, if we simply set "$x=0$" in the integral (5.26), we would get

$$\mathcal{G}(g_0)(x,t)\big|_{x=0} = \int_{-\infty}^{\infty} [e^{-\omega_1(\lambda)t}\tilde{g}_0(\omega_1(\lambda),t) - e^{-\omega_2(\lambda)t}\tilde{g}_0(\omega_2(\lambda),t)]\frac{\lambda d\lambda}{\rho(\lambda)} = 0, \tag{5.27}$$

since the integrant in (5.27) is an odd function of $\lambda$. In other words, for $t>0$,

$$\lim_{x\to 0^+}\left\{\lim_{\mathrm{A}\to\infty}\int_{-\mathrm{A}}^{\mathrm{A}} [e^{i\lambda x-\omega_1(\lambda)t}\tilde{g}_0(\omega_1(\lambda),t) - e^{i\lambda x-\omega_2(\lambda)t}\tilde{g}_0(\omega_2(\lambda),t)]\frac{\lambda d\lambda}{\rho(\lambda)}\right\} = 2\pi g_0(t)$$

while

$$\lim_{\mathrm{A}\to\infty}\left\{\int_{-\mathrm{A}}^{\mathrm{A}} \lim_{x\to 0^+}[e^{i\lambda x-\omega_1(\lambda)t}\tilde{g}_0(\omega_1(\lambda),t) - e^{i\lambda x-\omega_2(\lambda)t}\tilde{g}_0(\omega_2(\lambda),t)]\frac{\lambda d\lambda}{\rho(\lambda)}\right\} = 0.$$

## 6. Proof of Theorem 3

*Proof of 1st and 2nd parts* The existence of the limits of the derivatives, when $x\to 0^+$ or $t\to 0^+$, claimed in the 1st and 2nd parts, follows from the formulas (3.19) – (3.31), which give the interpretation of the various integrals. For example, differentiating (3.22), we obtain a formula for the derivative

$$\frac{\partial^{n+m}}{\partial x^n \partial t^m}\left[\int_{\mathbb{L}} e^{i\lambda x-\omega_1(\lambda)t}\tilde{\tilde{f}}(\lambda,\omega_1(\lambda),t)\frac{d\lambda}{2i\rho(\lambda)}\right], \text{ in the case for } x>t, \tag{6.1}$$

provided that $N$ and $M$ are chosen sufficiently large, and, from this formula, it follows that the limit of (6.1), as $t\to 0^+$, exists and defines the function $u_{n,m}$. Furthermore, it is easy to see that this limit is uniform for $x$ in compact of $(0,+\infty)$.

*Note* It is clear that integral representation formulas can be written for the derivatives of any part of the solution (1.5). Equation (4.14), used in the proof of part 4 of Theorem 1, is another formula expressing derivatives like (6.1). See also equation (7.2), below.

*Proof of 3rd part* The 2nd conclusion of the theorem guarantees that $u(x,t)$ extends to a $C^\infty$ function in $Q\cup\{(x,0):x>0\}$. Clearly, by the 1st conclusion of Theorem 2, for this extension of $u(x,t)$ we have

$$u(x,0) = u_0(x) \text{ and } \frac{\partial u}{\partial t}(x,0) = u_1(x),\ x>0. \tag{6.2}$$

Similarly, by the 1st conclusion, $u(x,t)$ extends to a $C^\infty$ function in $Q\cup\{(0,t):t>0\}$ and, by the 2nd conclusion of Theorem 2, we obtain that, for this extension of $u(x,t)$, we have $u(0,t)=g_0(t)$ for $t>0$.

*Proof of 4th and 5th parts* These follow from 3rd conclusion and the second equation in (6.2).

Finally, 6th assertion follows from the proof of 5th part of Theorem 2.

## 7. Proof of Theorem 4

Since, in this proof, $x\to+\infty$, we will use integral representation formulas which are valid for $x>t>0$.

***Step 1*** In view of (3.19) and (3.20),

$$\begin{aligned}\mathcal{U}_0^+(u_0)(x,t) &= \int_{-1}^{1} e^{i\lambda x}[\omega_1(\lambda)e^{-\omega_2(\lambda)t} - \omega_2(\lambda)e^{-\omega_1(\lambda)t}]\hat{u}_0(\lambda)\frac{d\lambda}{2i\rho(\lambda)} \\ &+ \sum_{j=1}^{N} u_0^{(j-1)}(0)\int_{\kappa^+} e^{i\lambda x}[\omega_1(\lambda)e^{-\omega_2(\lambda)t} - \omega_2(\lambda)e^{-\omega_1(\lambda)t}]\frac{1}{(i\lambda)^j}\frac{d\lambda}{2i\rho(\lambda)} \\ &+ \int_{\mathbb{L}} e^{i\lambda x}[\omega_1(\lambda)e^{-\omega_2(\lambda)t} - \omega_2(\lambda)e^{-\omega_1(\lambda)t}]\frac{1}{(i\lambda)^N}(u_0^{(N)}\hat{)}(\lambda)\frac{d\lambda}{2i\rho(\lambda)}.\end{aligned} \tag{7.1}$$

Differentiating (7.1) with sufficiently large $N$, we see that

$$(7.2)\quad \frac{\partial^n \mathcal{U}_0^+(u_0)(x,t)}{\partial x^n} = \int_{-1}^{1} (i\lambda)^n e^{i\lambda x}[\omega_1(\lambda)e^{-\omega_2(\lambda)t} - \omega_2(\lambda)e^{-\omega_1(\lambda)t}]\hat{u}_0(\lambda)\frac{d\lambda}{2i\rho(\lambda)}$$
$$+\sum_{j=1}^{N} u_0^{(j-1)}(0)\int_{\kappa^+}(i\lambda)^n e^{i\lambda x}[\omega_1(\lambda)e^{-\omega_2(\lambda)t} - \omega_2(\lambda)e^{-\omega_1(\lambda)t}]\frac{1}{(i\lambda)^j}\frac{d\lambda}{2i\rho(\lambda)}$$
$$+\int_{\mathbb{L}}(i\lambda)^n e^{i\lambda x}[\omega_1(\lambda)e^{-\omega_2(\lambda)t} - \omega_2(\lambda)e^{-\omega_1(\lambda)t}]\frac{1}{(i\lambda)^N}(u_0^{(N)})\hat{}(\lambda)\frac{d\lambda}{2i\rho(\lambda)}.$$

Using the equation $e^{i\lambda x} = \frac{1}{ix}\frac{\partial}{\partial\lambda}[e^{i\lambda x}]$ and integrating by parts we obtain

$$(7.3)\quad ix\frac{\partial^n \mathcal{U}_0^+(u_0)(x,t)}{\partial x^n} = -\int_{-1}^{1} e^{i\lambda x}\frac{\partial}{\partial\lambda}\left\{(i\lambda)^n[\omega_1(\lambda)e^{-\omega_2(\lambda)t} - \omega_2(\lambda)e^{-\omega_1(\lambda)t}]\hat{u}_0(\lambda)\frac{1}{2i\rho(\lambda)}\right\}d\lambda$$
$$-\sum_{j=1}^{N} u_0^{(j-1)}(0)\int_{\kappa^+} e^{i\lambda x}\frac{\partial}{\partial\lambda}\left\{(i\lambda)^n[\omega_1(\lambda)e^{-\omega_2(\lambda)t} - \omega_2(\lambda)e^{-\omega_1(\lambda)t}]\frac{1}{(i\lambda)^j}\frac{1}{2i\rho(\lambda)}\right\}d\lambda$$
$$-\int_{\mathbb{L}} e^{i\lambda x}\frac{\partial}{\partial\lambda}\left\{(i\lambda)^n[\omega_1(\lambda)e^{-\omega_2(\lambda)t} - \omega_2(\lambda)e^{-\omega_1(\lambda)t}]\frac{1}{(i\lambda)^N}(u_0^{(N)})\hat{}(\lambda)\frac{1}{2i\rho(\lambda)}\right\}d\lambda,$$

since the "intermediate boundary terms" (i.e., the evaluations at the points $\lambda = -1$ and $\lambda = 1$) cancel each other while the evaluations at $\lambda = -\infty$ and $\lambda = +\infty$ vanish.
It follows from (7.3) that

$$(7.4)\quad \sup\left\{\left|x\frac{\partial^n \mathcal{U}_0^+(u_0)(x,t)}{\partial x^n}\right| : x \geq 1,\, 0 < t \leq T\right\} < \infty.$$

Let us point out that, in order to conclude (7.4), it is crucial that $\left|e^{i\lambda x}\right| \leq 1$ when $\lambda$ lies on any of the contours of the integrals in (7.3).
Further integration by parts show that

$$(7.5)\quad \sup\left\{\left|x^\ell\frac{\partial^n \mathcal{U}_0^+(u_0)(x,t)}{\partial x^n}\right| : x \geq 1,\, 0 < t \leq T\right\} < \infty\ \ (\forall \ell).$$

In order to carry out these integration by parts processes in the integrals $\int_{-1}^{1}\cdots d\lambda$ which appear in the RHS of (7.3) and of its analogues, obtained for the representation of the quantities $x^\ell[\partial^n \mathcal{U}_0^+(u_0)(x,t)]/\partial x^n$, it is crucial that $\mathcal{K}_0(\lambda,t)$ is analytic in $\lambda$, in view of Lemma 1.

***Step 2*** Writing

$$\mathcal{U}_0^-(u_0)(x,t) = \left(\int_{-1}^{1} + \int_{\kappa^+}\right)\left\{[\omega_2(\lambda)e^{i\lambda x-\omega_1(\lambda)t} - \omega_1(\lambda)e^{i\lambda x-\omega_2(\lambda)t}]\hat{u}_0(-\lambda)\frac{d\lambda}{2i\rho(\lambda)}\right\}$$

differentiating and integrating by parts, we obtain

$$(ix)^\ell\frac{\partial^n \mathcal{U}_0^-(u_0)(x,t)}{\partial x^n} = (-1)^\ell\left(\int_{-1}^{1} + \int_{\kappa^+}\right)\left[e^{i\lambda x}\frac{\partial^\ell}{\partial\lambda^\ell}\left\{(i\lambda)^\ell[\omega_2(\lambda)e^{-\omega_1(\lambda)t} - \omega_1(\lambda)e^{-\omega_2(\lambda)t}]\hat{u}_0(-\lambda)\frac{d\lambda}{2i\rho(\lambda)}\right\}\right].$$

It follows that (7.5) holds with $\mathcal{U}_0^+$ replaced by $\mathcal{U}_0^-$.

***Step 3*** Writing

$$\mathcal{G}(g_0)(x,t) = \int_{-1}^{1} e^{i\lambda x}[e^{-\omega_1(\lambda)t}\tilde{g}_0(\omega_1(\lambda),t) - e^{-\omega_2(\lambda)t}\tilde{g}_0(\omega_2(\lambda),t)]\frac{\lambda d\lambda}{\rho(\lambda)}$$
$$+\sum_{j=1}^{N}(-1)^{j-1}\left\{g_0^{(j-1)}(t)\int_{\kappa^+} e^{i\lambda x}\left(\frac{1}{[\omega_1(\lambda)]^j} - \frac{1}{[\omega_2(\lambda)]^j}\right)\frac{\lambda d\lambda}{\rho(\lambda)} - g_0^{(j-1)}(0)\int_{\kappa^+} e^{i\lambda x}\left(\frac{e^{-\omega_1(\lambda)t}}{[\omega_1(\lambda)]^j} - \frac{e^{-\omega_2(\lambda)t}}{[\omega_2(\lambda)]^j}\right)\frac{\lambda d\lambda}{\rho(\lambda)}\right\}$$

$$+(-1)^N \int_{\mathbb{L}} e^{i\lambda x}\left\{\frac{1}{[\omega_1(\lambda)]^N} e^{-\omega_1(\lambda)t}(g_0^{(N)})^{\sim}(\omega_1(\lambda),t)-\frac{1}{[\omega_2(\lambda)]^N} e^{-\omega_2(\lambda)t}(g_0^{(N)})^{\sim}(\omega_2(\lambda),t)\right\}\frac{\lambda d\lambda}{\rho(\lambda)}$$

and working as previously, we see that (7.5) holds with $\mathcal{U}_0^+$ replaced by $\mathcal{G}(g_0)$. Again, the analyticity of $\mathcal{L}_{g_0}(\lambda,t)$, by Lemma 1, is crucial in handling the integrals $\int_{-1}^{1}\cdots d\lambda$.

***Step 4*** Working as previously, we easily see that (7.5) holds with $\mathcal{U}_0^+$ replaced by any of the functions $\mathcal{U}_1^+$, $\mathcal{U}_1^-$, $\mathcal{F}^+(f)$ or $\mathcal{F}^-(f)$.

*Completion of the proof* The conclusion of the theorem follows from the results of steps 1-4.

***Remark*** Extending the above proof, it is easy to see that, more generally,

$$\lim_{x\to+\infty}\left[x^{\ell}\frac{\partial^{n+m}u(x,t)}{\partial x^n\partial t^m}\right]=0,$$

for nonnegative integers $n$, $m$ and $\ell$, uniformly for $t$ in compact subsets of $[0,\infty)$.

## 8. A uniqueness theorem

**Theorem 5** *If we assume, in addition, that the data satisfy (1.4), then the solution $u(x,t)$, defined by (1.2), is the unique solution of (1.1) in the following sense: If $v(x,t)$ is $C^2$ for $(x,t)\in Q-\{(0,0)\}$, satisfies the equation $v_{tt}+v_t-v_{xx}=f$ in $Q$, the conditions*

$$v(x,0)=u_0(x) \text{ and } v_t(x,0)=u_1(x), \text{ for } x>0, \text{ and } v(0,t)=g_0(t), \text{ for } t>0,$$

(8.1) $$\lim_{x\to\infty} v_x(x,t)=0,\ \sup_{x\geq 1}|v_t(x,t)|<\infty, \text{ for } t>0,$$

*and, for every $T>0$, the functions $|v_t(x,t)|^2$, $|v_{tt}(x,t)|^2$, $|v_x(x,t)|^2$, $|v_{xx}(x,t)|^2$, $|v_{tx}(x,t)|^2$ are, uniformly for $0<t\leq T$, integrable with respect to $x\in[0,\infty)$, i.e., there exists a positive function $\mathrm{B}_T(x)$ such that $\int_0^\infty \mathrm{B}_T(x)dx<+\infty$ and, for $0<t\leq T$,*

(8.2) $$|v_t(x,t)|^2\leq \mathrm{B}_T(x),\ |v_{tt}(x,t)|^2\leq \mathrm{B}_T(x),\ |v_x(x,t)|^2\leq \mathrm{B}_T(x),\ |v_{xx}(x,t)|^2\leq \mathrm{B}_T(x),\ |v_{tx}(x,t)|^2\leq \mathrm{B}_T(x)\ (x>0),$$

*then $v\equiv u$.*

**Proof** Assuming that the functions $u$ and $v$ are real-valued, we set $w(x,t):=u(x,t)-v(x,t)$. Firstly, the solution $u(x,t)$ satisfies the conditions (8.1) and (8.2), by Theorems 3 and 4. (In particular, by 6th part of Theorem 3, it follows that the derivatives of $u(x,t)$, of order $\leq 2$, are bounded for $(x,t)$ close to the point $(0,0)$.) Therefore, so does the function $w(x,t)$. (We may have to choose $\mathrm{B}_T(x)$ larger.)
Also

(8.3) $$w(x,0)=0 \text{ and } w_t(x,0)=0, \text{ for } x>0,\ w(0,t)=0, \text{ for } t>0,$$

and

(8.4) $$w_t(x,t)w_{tt}(x,t)+[w_t(x,t)]^2-w_t(x,t)w_{xx}(x,t)=0, \text{ for } x\geq 0,\ t>0.$$

Fixing $t>0$ and integrating (8.4), we obtain

(8.5) $$\frac{1}{2}\int_{x=0}^{\infty}\left(\frac{\partial}{\partial t}[w_t(x,t)]^2\right)dx+\int_{x=0}^{\infty}[w_t(x,t)]^2dx-\int_{x=0}^{\infty}w_t(x,t)w_{xx}(x,t)dx=0.$$

Note that the integrals in (8.5) are absolutely convergent, since $w$ satisfies (8.2).
Integrating by parts, we see that

(8.6) $$\int_{x=0}^{\infty}w_t(x,t)w_{xx}(x,t)dx=w_t(x,t)w_x(x,t)\Big|_{x=0}^{x=\infty}-\int_{x=0}^{\infty}w_{tx}(x,t)w_x(x,t)dx=-\int_{x=0}^{\infty}w_{tx}(x,t)w_x(x,t)dx,$$

since $w_t(x,t)w_x(x,t)\Big|_{x=0}^{x=\infty}=0$. This last equation follows from the last equation in (8.3) and the fact that (8.1) holds with with $w$ in place of $v$.
Substituting (8.6) in (8.5), we obtain

(8.7) $$\frac{1}{2}\int_{x=0}^{\infty}\left(\frac{\partial}{\partial t}[w_t(x,t)]^2\right)dx+\int_{x=0}^{\infty}[w_t(x,t)]^2\,dx+\int_{x=0}^{\infty}w_{tx}(x,t)w_x(x,t)dx=0.$$

By Lebesgue's dominated convergence theorem and the fact that (8.2) holds with $w$ in place of $v$, (8.7) gives

$$\frac{1}{2}\frac{\partial}{\partial t}\int_{x=0}^{\infty}\left([w_t(x,t)]^2+[w_x(x,t)]^2\right)dx=-\int_{x=0}^{\infty}[w_t(x,t)]^2\,dx \;\Rightarrow\; \frac{\partial}{\partial t}\int_{x=0}^{\infty}\left([w_t(x,t)]^2+[w_x(x,t)]^2\right)dx\le 0 \text{ for } t>0.$$

Thus the function

$$\int_{x=0}^{\infty}\left([w_t(x,t)]^2+[w_x(x,t)]^2\right)dx \text{ is non-increasing for } t>0\text{, i.e.,}$$

(8.8) $$\int_{x=0}^{\infty}\left([w_t(x,t)]^2+[w_x(x,t)]^2\right)dx\le\int_{x=0}^{\infty}\left([w_t(x,\tau)]^2+[w_x(x,\tau)]^2\right)dx \text{ for } 0<\tau<t.$$

Letting $\tau\to 0^+$, again by Lebesgue's dominated convergence theorem, we obtain

$$\int_{x=0}^{\infty}\left([w_t(x,t)]^2+[w_x(x,t)]^2\right)dx\le 0 \text{ for } t>0,$$

since, by the first two equations in (8.3), $w_t(x,0)=w_x(x,0)=0$.
Therefore,

$$\int_{x=0}^{\infty}\left([w_t(x,t)]^2+[w_x(x,t)]^2\right)dx=0 \;\Rightarrow\; w_t(x,t)=0 \text{ and } w_x(x,t)=0.$$

But

$$w_t(x,t)=0 \;\Rightarrow\; w(x,t)=w(x,0)=0\text{, i.e., } v\equiv u.$$

(Analogous results on well-posedness have also been obtained for other linear PDEs in [14,6,7,10].

## 9. A non-controlability theorem

**Theorem 6** *With the notation of Lemma 1, suppose that the data* $u_0$ *and* $u_1$ *are such that the function* $\mathcal{K}_0(\lambda,t)\hat{u}_0(\lambda)+\mathcal{K}_1(\lambda,t)\hat{u}_1(\lambda)$ *does not extend to an analytic function of* $\lambda\in\mathbb{C}$*, while* $f$ *is such that, for some* $T>0$*, the function* $\hat{f}(\lambda,T)$ *extends to an entire function, i.e. analytic for* $\lambda\in\mathbb{C}$*. Then, for every* $g_0$ *with*

(9.1) $$g_0(0)=u_0(0),\; g_0'(0)=u_1(0) \text{ and } g_0''(0)=u_0''(0)-u_1(0)+f(0,0),$$

*we have that the solution of (1.1), given by (1.5), satisfies* $u(x,T)\not\equiv 0$.
*In particular, the homogenous version of problem (1.5) (i.e., with* $f\equiv 0$*) is not null controllable.*

**Proof** Let $T$, $g_0$ and $f$ be as in the statement of the theorem. Suppose – to reach a contradiction – that $u(x,T)\equiv 0$. Then $\hat{u}(\lambda,T)=0$, and (2.6) gives the following equation

(9.2) $$\mathcal{K}_0(\lambda,T)\hat{u}_0(\lambda)+\mathcal{K}_1(\lambda,T)\hat{u}_1(\lambda)+\mathcal{L}_{g_0}(\lambda,T)+\mathcal{L}_{g_1}(\lambda,T)+\mathcal{M}_f^{+}(\lambda,T)=0, \text{ for } \lambda\in\mathbb{C} \text{ with } \operatorname{Im}\lambda\le 0.$$

We emphasize that it is because of the assumptions on the data, namely (9.1), that the derivation of (9.2) is ***rigorous***. Indeed, this follows from Theorems 1, 2, 3, and 4, which guarantee the behavior of the solution $u(x,t)$ on the boundary of $Q$, i.e., when $x=0$ or $t=0$, including the point $(0,0)$, and when $x\to\infty$. In particular, the existence of the limits

(9.3) $$\lim_{\substack{(x,t)\to(0,0)\\(x,t)\in\overline{Q}-\{(0,0)\}}}\frac{\partial^{n+m}u(x,t)}{\partial x^n\partial t^m}\text{, for } n+m\le 2,$$

proved in part 6, of Theorem 3 (using (9.1)), imply that the function $\mathcal{L}_{g_0}(\lambda,T)+\mathcal{L}_{g_1}(\lambda,T)$ is entire.

Also, the use of Green's theorem, which was used in the derivation of (2.6), is legitimate, since $u(x,t)$ is $C^2$ in an open neighbourhood of $\overline{Q}-\{(0,0)\}$ and the existence of the limits (6.3) imply that $\partial^{n+m}u(x,t)/\partial x^n \partial t^m$, $n+m\le 2$, are bounded for $(x,t)$ close to $(0,0)$.
Therefore, (9.2) implies that the function $\mathcal{K}_0(\lambda,t)\hat{u}_0(\lambda)+\mathcal{K}_1(\lambda,t)\hat{u}_1(\lambda)$ extends to an analytic function of $\lambda\in\mathbb{C}$. This contradicts our assumption and proves the conclusion of the theorem.

(See [9] for a *generic* property – in the sense of Baire's category theorem – of the Fourier-Laplace transform for functions in the half-line Schwartz space.)

## 10. The equation $\alpha u_{tt}+\beta u_t-u_{xx}=f$ with $\alpha>0$ and $\beta\in\mathbb{R}-\{0\}$

We can easily extend the results of the previous sections to the more general problem

$$(10.1)\qquad \begin{cases} \alpha\dfrac{\partial^2 u}{\partial t^2}+\beta\dfrac{\partial u}{\partial t}-\dfrac{\partial^2 u}{\partial x^2}=f,\ (x,t)\in Q:=\mathbb{R}^+\times\mathbb{R}^+ \\ \lim\limits_{t\to 0^+} u(x,t)=u_0(x),\ x\in\mathbb{R}^+ \\ \lim\limits_{t\to 0^+}\dfrac{\partial u(x,t)}{\partial t}=u_1(x),\ x\in\mathbb{R}^+ \\ \lim\limits_{x\to 0^+} u(x,t)=g_0(t),\ t\in\mathbb{R}^+. \end{cases}$$

Indeed, the function $u(x,t)=e^{i\lambda x-\omega(\lambda)t}$ satisfies the equation in (10.1) if $\alpha[\omega(\lambda)]^2-\beta\omega(\lambda)+\lambda^2=0$. The roots of this quadratic equation are

$$\omega=\frac{\beta}{2\alpha}\pm\frac{i}{\sqrt{\alpha}}\sqrt{\lambda^2-\left(\frac{\beta}{2\sqrt{\alpha}}\right)^2}=\frac{\beta}{2\alpha}\pm\frac{i}{\sqrt{\alpha}}\rho(\lambda).$$

(The notation in this section is appropriately adjusted to the introduction of the parameters $\alpha$ and $\beta$.) We choose the following branch of the square root:

$$\rho(\lambda)=\rho_{\alpha,\beta}(\lambda):=\sqrt{\left|\lambda^2-\left(\frac{\beta}{2\sqrt{\alpha}}\right)^2\right|}\exp\left[i\left(\frac{\theta_1(\lambda)+\theta_2(\lambda)}{2}\right)\right]\ \text{for}\ \lambda\in\mathbb{C}-[-\tfrac{\beta}{2\sqrt{\alpha}},\tfrac{\beta}{2\sqrt{\alpha}}].$$

(The angles $\theta_1(\lambda)=\theta_{1,\alpha,\beta}(\lambda)$ and $\theta_2(\lambda)=\theta_{2,\alpha,\beta}(\lambda)$ are depicted in fig. 7.)

We extend $\rho(\lambda)=\rho_{\alpha,\beta}(\lambda)$ also for $\lambda\in[-\tfrac{\beta}{2\sqrt{\alpha}},\tfrac{\beta}{2\sqrt{\alpha}}]$ as follows:

$$\rho(\lambda):=\lim_{\substack{\mu\to\lambda\\ \operatorname{Im}\mu<0}}\rho(\mu)=-i\sqrt{\left|\lambda^2-\left(\frac{\beta}{2\sqrt{\alpha}}\right)^2\right|},\ \lambda\in[-\tfrac{\beta}{2\sqrt{\alpha}},\tfrac{\beta}{2\sqrt{\alpha}}].$$

Then

$$\omega_1(\lambda)=\omega_{1,\alpha,\beta}(\lambda):=\frac{\beta}{2\alpha}+\frac{i}{\sqrt{\alpha}}\rho_{\alpha,\beta}(\lambda)\ \text{and}\ \omega_2(\lambda)=\omega_{2,\alpha,\beta}(\lambda):=\frac{\beta}{2\alpha}-\frac{i}{\sqrt{\alpha}}\rho(\lambda)$$

are analytic for $\lambda\in\mathbb{C}-[-\tfrac{\beta}{2\sqrt{\alpha}},\tfrac{\beta}{2\sqrt{\alpha}}]$, their restriction to the (closed) lower half-plane $\{\lambda\in\mathbb{C}:\operatorname{Im}\lambda\le 0\}$ is continuous and $C^\infty$ in $\{\lambda\in\mathbb{C}:\operatorname{Im}\lambda\le 0\}-\{-\tfrac{\beta}{2\sqrt{\alpha}},\tfrac{\beta}{2\sqrt{\alpha}}\}$.

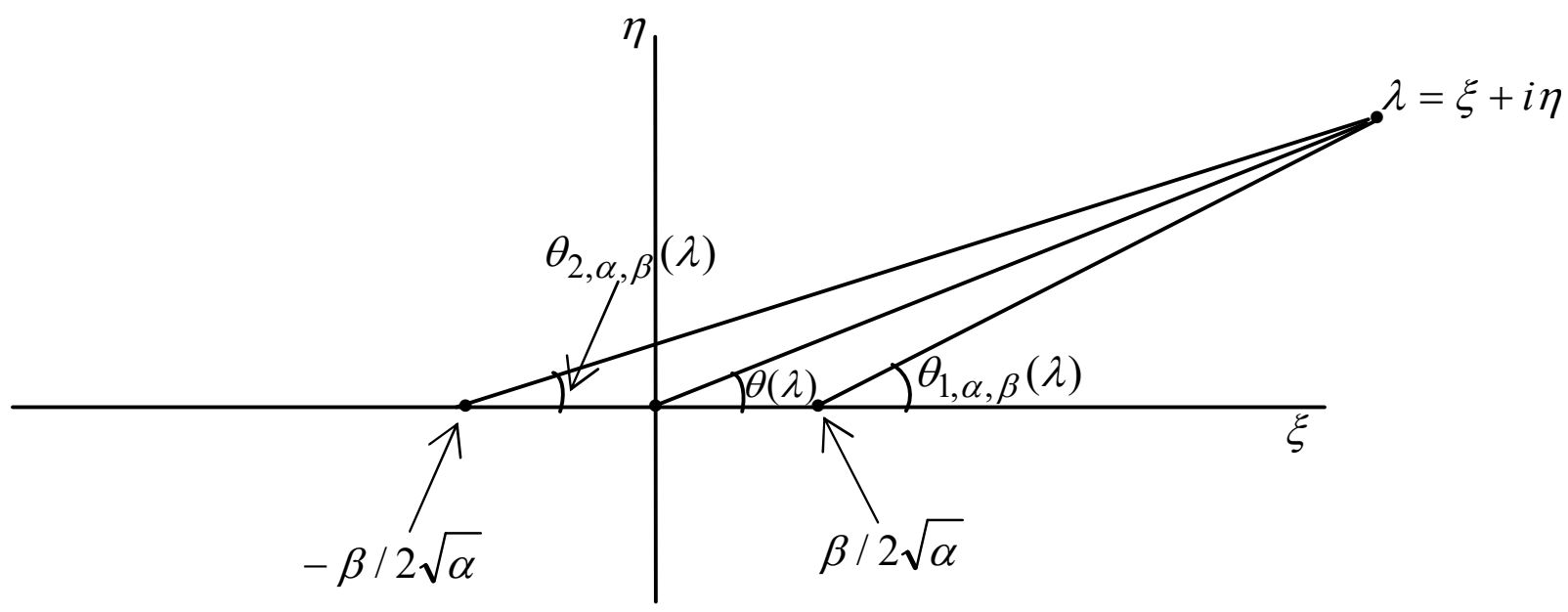

**Fig. 7** $\lambda-(1/2\sqrt{\alpha})=\left|\lambda-(1/2\sqrt{\alpha})\right|e^{i\theta_{1,\alpha,\beta}(\lambda)}$, $\lambda-(-1/2\sqrt{\alpha})=\left|\lambda-(-1/2\sqrt{\alpha})\right|e^{i\theta_{2,\alpha,\beta}(\lambda)}$, $-\pi<\theta_{1,\alpha,\beta}(\lambda),\ \theta_{2,\alpha,\beta}(\lambda)\le\pi$ .

With this notation, the Fokas solution $u(x,t)=u_{\alpha,\beta}(x,t)$ of Problem (10.1) is given by the following formula: For $x>0$ and $t>0$,

$$
\begin{aligned}
(10.2)\quad 2\pi u_{\alpha,\beta}(x,t) &= \sqrt{\alpha}\int_{-\infty}^{\infty}[\omega_1(\lambda)e^{i\lambda x-\omega_2(\lambda)t}-\omega_2(\lambda)e^{i\lambda x-\omega_1(\lambda)t}]\hat{u}_0(\lambda)\frac{d\lambda}{2i\rho(\lambda)} \\
&+\sqrt{\alpha}\int_{-\infty}^{\infty}[\omega_2(\lambda)e^{i\lambda x-\omega_1(\lambda)t}-\omega_1(\lambda)e^{i\lambda x-\omega_2(\lambda)t}]\hat{u}_0(-\lambda)\frac{d\lambda}{2i\rho(\lambda)} \\
&+\sqrt{\alpha}\int_{-\infty}^{\infty}[e^{i\lambda x-\omega_2(\lambda)t}-e^{i\lambda x-\omega_1(\lambda)t}]\hat{u}_1(\lambda)\frac{d\lambda}{2i\rho(\lambda)}+\sqrt{\alpha}\int_{-\infty}^{\infty}[e^{i\lambda x-\omega_1(\lambda)t}-e^{i\lambda x-\omega_2(\lambda)t}]\hat{u}_1(-\lambda)\frac{d\lambda}{2i\rho(\lambda)} \\
&+\frac{1}{\sqrt{\alpha}}\int_{-\infty}^{\infty}[e^{i\lambda x-\omega_1(\lambda)t}\tilde{g}_0(\omega_1(\lambda),t)-e^{i\lambda x-\omega_2(\lambda)t}\tilde{g}_0(\omega_2(\lambda),t)]\frac{\lambda d\lambda}{\rho(\lambda)} \\
&+\frac{1}{\sqrt{\alpha}}\int_{-\infty}^{\infty}\left[e^{i\lambda x-\omega_2(\lambda)t}\tilde{\hat{f}}(\lambda,\omega_2(\lambda),t)-e^{i\lambda x-\omega_1(\lambda)t}\tilde{\hat{f}}(\lambda,\omega_1(\lambda),t)\right]\frac{d\lambda}{2i\rho(\lambda)} \\
&+\frac{1}{\sqrt{\alpha}}\int_{-\infty}^{\infty}\left[e^{i\lambda x-\omega_1(\lambda)t}\tilde{\hat{f}}(-\lambda,\omega_1(\lambda),t)-e^{i\lambda x-\omega_2(\lambda)t}\tilde{\hat{f}}(-\lambda,\omega_2(\lambda),t)\right]\frac{d\lambda}{2i\rho(\lambda)}.
\end{aligned}
$$

***Formal derivation of the solution*** The differential equation in (10.1) can be written in divergence form as follows:

$$
(10.3)\qquad \frac{\partial}{\partial t}\left\{\left([\beta-\alpha\omega(\lambda)]u+\alpha u_t\right)e^{-i\lambda x+\omega(\lambda)t}\right\}-\frac{\partial}{\partial x}\left[(i\lambda u+u_x)e^{-i\lambda x+\omega(\lambda)t}\right]=f(x,t)e^{-i\lambda x+\omega(\lambda)t}.
$$

Using Green's theorem, (10.3) gives

$$
\begin{aligned}
(10.4)\quad &\int_{x=0}^{\infty}\left\{\left([\beta-\alpha\omega(\lambda)]u+\alpha u_t\right)e^{-i\lambda x+\omega(\lambda)t}\right\}\Big|_{t=0}dx-\int_{x=0}^{\infty}\left\{\left([\beta-\alpha\omega(\lambda)]u+\alpha u_t\right)e^{-i\lambda x+\omega(\lambda)\tau}\right\}\Big|_{\tau=t}dx \\
&-\int_{\tau=0}^{t}\left[(i\lambda u+u_x)e^{-i\lambda x+\omega(\lambda)\tau}\right]\Big|_{x=0}d\tau=-\int_{\tau=0}^{t}e^{\omega(\lambda)\tau}\left[\int_{x=0}^{\infty}f(x,\tau)e^{-i\lambda x}dx\right]d\tau,\quad \operatorname{Im}\lambda\le 0.
\end{aligned}
$$

Then (10.4) leads to the equation

$$
\begin{aligned}
&\{[\beta-\alpha\omega(\lambda)]\hat{u}_0(\lambda)+\alpha\hat{u}_1(\lambda)\}-\{[\beta-\alpha\omega(\lambda)]\hat{u}(\lambda,t)e^{\omega(\lambda)t}+\alpha\hat{(u_t)}(\lambda,t)e^{\omega(\lambda)t}\} \\
&\qquad=(i\lambda)\tilde{g}_0(\omega(\lambda),t)-\tilde{g}_1(\omega(\lambda),t)-\tilde{\hat{f}}(\lambda,\omega(\lambda),t).
\end{aligned}
$$

Thus,

$$
\begin{aligned}
(10.5)\quad &\{[\beta-\alpha\omega(\lambda)]e^{-\omega(\lambda)t}\hat{u}_0(\lambda)+\alpha e^{-\omega(\lambda)t}\hat{u}_1(\lambda)\}-\{[\beta-\alpha\omega(\lambda)]\hat{u}(\lambda,t)+\alpha\hat{(u_t)}(\lambda,t)\} \\
&\qquad=(i\lambda)e^{-\omega(\lambda)t}\tilde{g}_0(\omega(\lambda),t)-e^{-\omega(\lambda)t}\tilde{g}_1(\omega(\lambda),t)-e^{-\omega(\lambda)t}\tilde{\hat{f}}(\lambda,\omega(\lambda),t).
\end{aligned}
$$

Writing (10.5) for each $\omega$ with $\omega\in\{\omega_1,\omega_2\}$ and taking into consideration that $\omega_1(\lambda)+\omega_2(\lambda)=\beta/\alpha$, we obtain

$$
\begin{aligned}
(10.6)\quad &[\alpha\omega_2 e^{-\omega_1 t}\hat{u}_0(\lambda)+\alpha e^{-\omega_1 t}\hat{u}_1(\lambda)]-[\alpha\omega_2\hat{u}(\lambda,t)+\alpha\hat{(u_t)}(\lambda,t)] \\
&\qquad=(i\lambda)e^{-\omega_1 t}\tilde{g}_0(\omega_1,t)-e^{-\omega_1 t}\tilde{g}_1(\omega_1,t)-e^{-\omega_1 t}\tilde{\hat{f}}(\lambda,\omega_1,t).
\end{aligned}
$$

$$(10.7)\quad [\alpha\omega_1 e^{-\omega_2 t}\hat{u}_0(\lambda)+\alpha e^{-\omega_2 t}\hat{u}_1(\lambda)]-[\alpha\omega_1\hat{u}(\lambda,t)+\alpha\hat{(u_t)}(\lambda,t)]$$
$$=(i\lambda)e^{-\omega_2 t}\tilde{g}_0(\omega_2,t)-e^{-\omega_2 t}\tilde{g}_1(\omega_2,t)-e^{-\omega_2 t}\tilde{\hat{f}}(\lambda,\omega_2,t).$$

Subtracting (10.7) from (10.6), we have that, for $\operatorname{Im}\lambda\le 0$,

$$\alpha(\omega_2 e^{-\omega_1 t}-\omega_1 e^{-\omega_2 t})\hat{u}_0(\lambda)+\alpha(e^{-\omega_1 t}-e^{-\omega_2 t})\hat{u}_1(\lambda)+\alpha(\omega_1-\omega_2)\hat{u}(\lambda,t)$$
$$=(i\lambda)e^{-\omega_1 t}\tilde{g}_0(\omega_1,t)-(i\lambda)e^{-\omega_2 t}\tilde{g}_0(\omega_2,t)-e^{-\omega_1 t}\tilde{g}_1(\omega_1,t)+e^{-\omega_2 t}\tilde{g}_1(\omega_2,t)-e^{-\omega_1 t}\tilde{\hat{f}}(\lambda,\omega_1,t)+e^{-\omega_2 t}\tilde{\hat{f}}(\lambda,\omega_2,t),$$

and, therefore,

$$(10.8)\quad \frac{1}{\omega_1-\omega_2}(\omega_2 e^{-\omega_1 t}-\omega_1 e^{-\omega_2 t})\hat{u}_0(\lambda)+\frac{1}{\omega_1-\omega_2}(e^{-\omega_1 t}-e^{-\omega_2 t})\hat{u}_1(\lambda)+\hat{u}(\lambda,t)$$
$$=\frac{i\lambda}{\alpha(\omega_1-\omega_2)}e^{-\omega_1 t}\tilde{g}_0(\omega_1,t)-\frac{i\lambda}{\alpha(\omega_1-\omega_2)}e^{-\omega_2 t}\tilde{g}_0(\omega_2,t)-\frac{1}{\alpha(\omega_1-\omega_2)}e^{-\omega_1 t}\tilde{g}_1(\omega_1,t)+\frac{1}{\alpha(\omega_1-\omega_2)}e^{-\omega_2 t}\tilde{g}_1(\omega_2,t)$$
$$-\frac{1}{\alpha(\omega_1-\omega_2)}e^{-\omega_1 t}\tilde{\hat{f}}(\lambda,\omega_1,t)+\frac{1}{\alpha(\omega_1-\omega_2)}e^{-\omega_2 t}\tilde{\hat{f}}(\lambda,\omega_2,t),\ \text{for } \operatorname{Im}\lambda\le 0,\ \lambda\ne\frac{\beta}{2\sqrt{\alpha}}.$$

Multiplying (10.8) by $e^{i\lambda x}$ and integrating we obtain

$$(10.9)\quad \int_{-\infty}^{\infty}\frac{[\omega_2(\lambda)e^{i\lambda x-\omega_1(\lambda)t}-\omega_1(\lambda)e^{i\lambda x-\omega_2(\lambda)t}]}{\omega_1(\lambda)-\omega_2(\lambda)}\hat{u}_0(\lambda)d\lambda+\int_{-\infty}^{\infty}\frac{e^{i\lambda x-\omega_1(\lambda)t}-e^{i\lambda x-\omega_2(\lambda)t}}{\omega_1(\lambda)-\omega_2(\lambda)}\hat{u}_1(\lambda)d\lambda+\int_{-\infty}^{\infty}e^{i\lambda x}\hat{u}(\lambda,t)d\lambda$$
$$=\int_{-\infty}^{\infty}\left\{\frac{i\lambda e^{i\lambda x-\omega_1(\lambda)t}}{\alpha[\omega_1(\lambda)-\omega_2(\lambda)]}\tilde{g}_0(\omega_1(\lambda),t)-\frac{i\lambda e^{i\lambda x-\omega_2(\lambda)t}}{\alpha[\omega_1(\lambda)-\omega_2(\lambda)]}\tilde{g}_0(\omega_2(\lambda),t)\right\}d\lambda$$
$$-\int_{-\infty}^{\infty}\left\{\frac{e^{i\lambda x-\omega_1(\lambda)t}}{\alpha[\omega_1(\lambda)-\omega_2(\lambda)]}\tilde{g}_1(\omega_1(\lambda),t)-\frac{e^{i\lambda x-\omega_2(\lambda)t}}{\alpha[\omega_1(\lambda)-\omega_2(\lambda)]}\tilde{g}_1(\omega_2(\lambda),t)\right\}d\lambda$$
$$-\int_{-\infty}^{\infty}\left\{\frac{e^{i\lambda x-\omega_1(\lambda)t}}{\alpha[\omega_1(\lambda)-\omega_2(\lambda)]}\tilde{\hat{f}}(\lambda,\omega_1(\lambda),t)-\frac{e^{i\lambda x-\omega_2(\lambda)t}}{\alpha[\omega_1(\lambda)-\omega_2(\lambda)]}\tilde{\hat{f}}(\lambda,\omega_2(\lambda),t)\right\}d\lambda,\ \text{for } \operatorname{Im}\lambda\le 0.$$

Now for $\lambda\in\mathbb{R}$, we have $\operatorname{Im}(-\lambda)\le 0$ and, therefore, (10.8) gives

$$(10.10)\quad \frac{1}{\omega_1(\lambda)-\omega_2(\lambda)}\{\omega_2(\lambda)e^{-\omega_1(\lambda)t}-\omega_1(\lambda)e^{-\omega_2(\lambda)t}\}\hat{u}_0(-\lambda)+\frac{1}{\omega_1(\lambda)-\omega_2(\lambda)}(e^{-\omega_1(\lambda)t}-e^{-\omega_2(\lambda)t})\hat{u}_1(-\lambda)+\hat{u}(-\lambda,t)$$
$$=\frac{-i\lambda}{\alpha[\omega_1(\lambda)-\omega_2(\lambda)]}e^{-\omega_1(\lambda)t}\tilde{g}_0(\omega_1(\lambda),t)-\frac{-i\lambda}{\alpha[\omega_1(\lambda)-\omega_2(\lambda)]}e^{-\omega_2(\lambda)t}\tilde{g}_0(\omega_2(\lambda),t)$$
$$-\frac{1}{\alpha[\omega_1(\lambda)-\omega_2(\lambda)]}e^{-\omega_1(\lambda)t}\tilde{g}_1(\omega_1(\lambda),t)+\frac{1}{\alpha[\omega_1(\lambda)-\omega_2(\lambda)]}e^{-\omega_2(\lambda)t}\tilde{g}_1(\omega_2(\lambda),t)$$
$$-\frac{1}{\alpha[\omega_1(\lambda)-\omega_2(\lambda)]}e^{-\omega_1(\lambda)t}\tilde{\hat{f}}(-\lambda,\omega_1(\lambda),t)+\frac{1}{\alpha[\omega_1(\lambda)-\omega_2(\lambda)]}e^{-\omega_2(\lambda)t}\tilde{\hat{f}}(-\lambda,\omega_2(\lambda),t),\ \text{for } \lambda\in\mathbb{R}.$$

Multiplying (10.10) by $e^{i\lambda x}$, integrating it from $\lambda=-\infty$ to $\lambda=\infty$, and subtracting the resulting equation from (10.9), we obtain:

$$\int_{-\infty}^{\infty}\frac{\omega_2(\lambda)e^{i\lambda x-\omega_1(\lambda)t}-\omega_1(\lambda)e^{i\lambda x-\omega_2(\lambda)t}}{\omega_1(\lambda)-\omega_2(\lambda)}\hat{u}_0(\lambda)d\lambda-\int_{-\infty}^{\infty}\frac{\omega_2(\lambda)e^{i\lambda x-\omega_1(\lambda)t}-\omega_1(\lambda)e^{i\lambda x-\omega_2(\lambda)t}}{\omega_1(\lambda)-\omega_2(\lambda)}\hat{u}_0(-\lambda)d\lambda$$
$$+\int_{-\infty}^{\infty}\frac{e^{i\lambda x-\omega_1(\lambda)t}-e^{i\lambda x-\omega_2(\lambda)t}}{\omega_1(\lambda)-\omega_2(\lambda)}\hat{u}_1(\lambda)d\lambda-\int_{-\infty}^{\infty}\frac{e^{i\lambda x-\omega_1(\lambda)t}-e^{i\lambda x-\omega_2(\lambda)t}}{\omega_1(\lambda)-\omega_2(\lambda)}\hat{u}_1(-\lambda)d\lambda+2\pi u(x,t)$$
$$=2i\int_{-\infty}^{\infty}\left\{\frac{\lambda e^{i\lambda x-\omega_1(\lambda)t}}{\alpha[\omega_1(\lambda)-\omega_2(\lambda)]}\tilde{g}_0(\omega_1(\lambda),t)-\frac{\lambda e^{i\lambda x-\omega_2(\lambda)t}}{\alpha[\omega_1(\lambda)-\omega_2(\lambda)]}\tilde{g}_0(\omega_2(\lambda),t)\right\}d\lambda$$

$$-\int_{-\infty}^{\infty}\left\{\frac{e^{i\lambda x-\omega_1(\lambda)t}}{\alpha[\omega_1(\lambda)-\omega_2(\lambda)]}\tilde{\tilde{f}}(\lambda,\omega_1(\lambda),t)-\frac{e^{i\lambda x-\omega_2(\lambda)t}}{\alpha[\omega_1(\lambda)-\omega_2(\lambda)]}\tilde{\tilde{f}}(\lambda,\omega_2(\lambda),t)\right\}d\lambda$$
$$+\int_{-\infty}^{\infty}\left\{\frac{e^{i\lambda x-\omega_1(\lambda)t}}{\alpha[\omega_1(\lambda)-\omega_2(\lambda)]}\tilde{\tilde{f}}(-\lambda,\omega_1(\lambda),t)-\frac{e^{i\lambda x-\omega_2(\lambda)t}}{\alpha[\omega_1(\lambda)-\omega_2(\lambda)]}\tilde{\tilde{f}}(-\lambda,\omega_2(\lambda),t)\right\}d\lambda.$$

This gives (10.2) and completes the formal derivation of the solution.

## 11. The limit of the solution of (10.1) as $\beta\to 0^+$

In this section, we apply formula (10.2) with $\alpha=1$ and $\beta>0$, and we show that the limit, as $\beta\to 0^+$, of the solution $u_\beta(x,t):=u_{1,\beta}(x,t)=u_{\alpha,\beta}(x,t)\big|_{\alpha=1}$, is the solution of the following problem for the wave equation:

(11.1) $$u_{tt}=u_{xx}+f,\ u(x,0)=u_0(x)\ (x>0),\ u_t(x,0)=u_1(x)\ (x>0),\ u(0,t)=g_0(t)\ (t>0).$$

Thus,

(11.2) $$2\pi u_\beta(x,t)=[\mathcal{U}^+_{0,\beta}(u_0)+\mathcal{U}^-_{0,\beta}(u_0)+\mathcal{U}^+_{1,\beta}(u_1)+\mathcal{U}^-_{1,\beta}(u_1)+\mathcal{G}_\beta(g_0)+\mathcal{F}^+_\beta(f)+\mathcal{F}^-_\beta(f)](x,t),$$

is the analogue of (1.6), i.e.,

(11.3) $$\mathcal{U}^+_{0,\beta}(u_0),\ \mathcal{U}^-_{0,\beta}(u_0),\ \mathcal{U}^+_{1,\beta}(u_1),\ \mathcal{U}^-_{1,\beta}(u_1),\ \mathcal{G}_\beta(g_0),\ \mathcal{F}^+_\beta(f),\ \mathcal{F}^-_\beta(f)$$

are the functions of $(x,t)$, defined by the integrals in (10.2) with $\alpha=1$, which involve the terms

$$\hat{u}_0(\lambda),\ \hat{u}_0(-\lambda),\ \hat{u}_1(\lambda),\ \hat{u}_1(-\lambda),\ \tilde{g}_0(\lambda,\cdot),\ \tilde{\tilde{f}}(\lambda,\cdot,t),\ \tilde{\tilde{f}}(-\lambda,\cdot,t),$$

respectively.
Since the integrals $\mathcal{U}^+_{1,\beta}(u_1),\ \mathcal{U}^-_{1,\beta}(u_1),\mathcal{F}^+_\beta(f),\ \mathcal{F}^-_\beta(f)$ converge absolutely, we immediately obtain that

$$\lim_{\beta\to0^+}\mathcal{U}^+_{1,\beta}(u_1)=\mathcal{U}^+_{1,\beta}(u_1)\Big|_{\beta=0},\ \lim_{\beta\to0^+}\mathcal{U}^-_{1,\beta}(u_1)=\mathcal{U}^-_{1,\beta}(u_1)\Big|_{\beta=0},$$
$$\lim_{\beta\to0^+}\mathcal{F}^+_\beta(f)=\mathcal{F}^+_\beta(f)\Big|_{\beta=0},\ \lim_{\beta\to0^+}\mathcal{F}^-_\beta(f)=\mathcal{F}^-_\beta(f)\Big|_{\beta=0}.$$

The remaining integrals in (11.3) may not converge absolutely but they exist in the generalized sense. With appropriate deformation of the contours (using the analogues of the formulas of section 3), we prove that, in the case of these integrals, we also have

$$\lim_{\beta\to0^+}\mathcal{U}^+_{0,\beta}(u_0)=\mathcal{U}^+_{0,\beta}(u_0)\Big|_{\beta=0},\ \lim_{\beta\to0^+}\mathcal{U}^-_{0,\beta}(u_0)=\mathcal{U}^-_{0,\beta}(u_0)\Big|_{\beta=0},\ \lim_{\beta\to0^+}\mathcal{G}_\beta(g_0)=\mathcal{G}_\beta(g_0)\Big|_{\beta=0}.$$

Thus, letting $\beta\to 0^+$ in (11.2), we obtain

(11.4) $$2\pi\lim_{\beta\to0^+}u_\beta(x,t)=\frac{1}{2}\int_{-\infty}^{\infty}[e^{i\lambda(x+t)}+e^{i\lambda(x-t)}]\hat{u}_0(\lambda)d\lambda+\frac{1}{2}\int_{-\infty}^{\infty}[-e^{i\lambda(x+t)}-e^{i\lambda(x-t)}]\hat{u}_0(-\lambda)d\lambda$$
$$+\frac{1}{2}\int_{-\infty}^{\infty}[e^{i\lambda(x+t)}-e^{i\lambda(x-t)}]\hat{u}_1(\lambda)\frac{d\lambda}{i\lambda}+\frac{1}{2}\int_{-\infty}^{\infty}[e^{i\lambda(x-t)}-e^{i\lambda(x+t)}]\hat{u}_1(-\lambda)\frac{d\lambda}{i\lambda}$$
$$+\int_{-\infty}^{\infty}[e^{i\lambda(x-t)}\int_{\tau=0}^{t}e^{i\lambda\tau}g_0(\tau)d\tau-e^{i\lambda(x+t)}\int_{\tau=0}^{t}e^{-i\lambda\tau}g_0(\tau)d\tau]d\lambda$$
$$+\frac{1}{2}\int_{-\infty}^{\infty}\left[e^{i\lambda(x+t)}\int_{\tau=0}^{t}e^{-i\lambda\tau}\hat{f}(\lambda,\tau)d\tau-e^{i\lambda(x-t)}\int_{\tau=0}^{t}e^{i\lambda\tau}\hat{f}(\lambda,\tau)d\tau\right]\frac{d\lambda}{i\lambda}$$
$$+\frac{1}{2}\int_{-\infty}^{\infty}\left[e^{i\lambda(x-t)}\int_{\tau=0}^{t}e^{i\lambda\tau}\hat{f}(-\lambda,\tau)d\tau-e^{i\lambda(x+t)}\int_{\tau=0}^{t}e^{-i\lambda\tau}\hat{f}(-\lambda,\tau)d\tau\right]\frac{d\lambda}{i\lambda}.$$

We will show that the limit $\lim_{\beta\to0^+}u_\beta(x,t)$ is the solution of (11.1).

*Note* We have to be cautious about the fact that, not only does the function

$$\rho_{1,\beta}(\lambda)=\sqrt{\left|\lambda^2-\left(\frac{\beta}{2}\right)^2\right|}\exp\left[i\left(\frac{\theta_{1,1,\beta}(\lambda)+\theta_{2,1,\beta}(\lambda)}{2}\right)\right]$$

depend on $\beta$, but also the *branch cut*, namely the interval $[-\frac{\beta}{2},\frac{\beta}{2}]$, depends on $\beta$. However,

$$\lim_{\beta\to 0^+}\rho_{1,\beta}(\lambda)=\lambda\text{, for }\lambda\in\mathbb{C}-[-1,1]. \tag{11.5}$$

Thus, in order to obtain the limits of the integrals (11.3), as $\beta\to 0^+$, we may split these integrals according to the formula $\int_{-\infty}^{\infty}=\int_{-\infty}^{-1}+\int_{-1}^{1}+\int_{1}^{\infty}$, and work separately with the integrals $\int_{\mathbb{L}}=\int_{-\infty}^{-1}+\int_{1}^{\infty}$ and $\int_{-1}^{1}$. Then, for the case of the integrals over $\mathbb{L}$, we may use (11.5), while for the integrals over the interval $[-1,1]$ we may use the analogues of the formulas of Lemma 1. We point out that the interchange of the limit as $\beta\to 0^+$ with the integrals over the interval $[-1,1]$ is immediate and that the *square root* function does not appear in the analogues of the kernels (2.7) and (2.8).
Finally, we observe that, in view of (11.5),

$$\lim_{\beta\to 0^+}\omega_{1,1,\beta}(\lambda)=i\lambda \text{ and } \lim_{\beta\to 0^+}\omega_{2,1,\beta}(\lambda)=-i\lambda\text{, for }\lambda\in\mathbb{C}-[-1,1].$$

***Writing the function*** $\lim_{\beta\to 0^+}u_\beta(x,t)$ ***in d'Alembert form*** Computing the integrals in the RHS of (11.4), we will show that

$$\lim_{\beta\to 0^+}u_\beta(x,t)=\frac{1}{2}u_0(x+t)+\frac{1}{2}u_0(x-t)+\frac{1}{2}\int_{y=x-t}^{x+t}u_1(y)dy+\frac{1}{2}\int_{\tau=0}^{t}\left[\int_{y=x-(t-\tau)}^{x+(t-\tau)}f(y,\tau)dy\right]d\tau\text{, for }x>t, \tag{11.6}$$

and

$$\begin{aligned}\lim_{\beta\to 0^+}u_\beta(x,t)&=\frac{1}{2}u_0(t+x)-\frac{1}{2}u_0(t-x)+g_0(t-x)+\frac{1}{2}\int_{y=t-x}^{t+x}u_1(y)dy\\&+\frac{1}{2}\int_{\tau=0}^{t-x}\left[\int_{y=(t-x)-\tau}^{(t+x)-\tau}f(y,\tau)dy\right]d\tau+\frac{1}{2}\int_{\tau=t-x}^{t}\left[\int_{y=(x-t)+\tau}^{(x+t)-\tau}f(y,\tau)dy\right]d\tau\text{, for }t>x.\end{aligned} \tag{11.7}$$

Firstly, extending the notation (11.3) for $\beta=0$, let

$$\mathcal{U}_{0,0}^{+}(u_0),\ \mathcal{U}_{0,0}^{-}(u_0),\ \mathcal{U}_{1,0}^{+}(u_1),\ \mathcal{U}_{1,0}^{-}(u_1),\ \mathcal{G}_0(g_0),\ \mathcal{F}_0^{+}(f),\ \mathcal{F}_0^{-}(f)$$

be the functions of $(x,t)$, defined by the integrals in (11.4), which involve the terms

$$\hat{u}_0(\lambda),\ \hat{u}_0(-\lambda),\ \hat{u}_1(\lambda),\ \hat{u}_1(-\lambda),\ \tilde{g}_0(\lambda,\cdot),\ \tilde{\hat{f}}(\lambda,\cdot,t),\ \tilde{\hat{f}}(-\lambda,\cdot,t),$$

respectively.

***In the case*** $x>t$***,*** we have:

$$\mathcal{U}_{0,0}^{+}(u_0)+\mathcal{U}_{0,0}^{-}(u_0)=\frac{1}{2}u_0(x+t)+\frac{1}{2}u_0(x-t), \tag{11.8}$$

$$\mathcal{U}_{1,0}^{+}(u_0)+\mathcal{U}_{1,0}^{-}(u_0)=\frac{1}{2}\int_{y=x-t}^{x+t}u_1(y)dy, \tag{11.9}$$

$$\mathcal{G}_0(g_0)=0, \tag{11.10}$$

$$\mathcal{F}_0^{+}(f)+\mathcal{F}_0^{-}(f)=\frac{1}{2}\int_{\tau=0}^{t}\left[\int_{y=x-(t-\tau)}^{x+(t-\tau)}f(y,\tau)dy\right]d\tau. \tag{11.11}$$

Indeed, (11.8) follows by Fourier's inversion formula.

*Proof of (11.9)* Working with $x>t$ and differentiating with respect to $x$, we obtain, by Fourier's inversion formula,

$$\frac{\partial[\mathcal{U}_{1,0}^{+}(u_0)+\mathcal{U}_{1,0}^{-}(u_0)]}{\partial x}=\frac{1}{2}u_1(x+t)-\frac{1}{2}u_1(x-t)=\frac{1}{2}\frac{\partial}{\partial x}\left[\int_{y=x-t}^{x+t}u_1(y)dy\right].$$

Therefore,

(11.12) $$[\mathcal{U}_{1,0}^{+}(u_0)+\mathcal{U}_{1,0}^{-}(u_0)]-\frac{1}{2}\left[\int_{y=x-t}^{x+t}u_1(y)dy\right]=\varphi(t)\text{, for some function }\varphi(t).$$

Differentiating (11.12), with respect to $t$, we see that $\frac{d\varphi(t)}{dt}=0$, i.e., $\varphi(t)=\varphi(0)=0$. Thus, (11.12) gives (11.9).

By Fourier's inversion formula,

$$\int_{-\infty}^{\infty}[e^{i\lambda(x+t)}\int_{\tau=0}^{t}e^{-i\lambda\tau}g_0(\tau)d\tau]d\lambda=0\text{ and }\int_{-\infty}^{\infty}[e^{i\lambda(x-t)}\int_{\tau=0}^{t}e^{i\lambda\tau}g_0(\tau)d\tau]d\lambda=0,$$

and (11.10) follows.

*Proof of (11.11)* By Fubini's theorem,

$$\mathcal{F}_0^{+}(f)=\frac{1}{2}\int_{\tau=0}^{t}\left[\int_{\lambda=-\infty}^{\infty}\{e^{i\lambda[x+(t-\tau)]}-e^{i\lambda[x-(t-\tau)]}\}\frac{\hat{f}(\lambda,\tau)}{i\lambda}d\lambda\right]d\tau=\frac{1}{2}\int_{\tau=0}^{t}\left[\int_{y=x-(t-\tau)}^{x+(t-\tau)}f(y,\tau)dy\right]d\tau,$$

where for the last equation, we used (11.9). (*Note*. The inner integral, i.e., the $d\lambda-integral$, is similar to $\mathcal{U}_{1,0}^{+}(u_0)$.) Working similarly, we see that $\mathcal{F}_0^{-}(f)=0$, and (11.11) follows.

***In the case*** $x<t$***,*** we have:

(11.13) $\mathcal{U}_{0,0}^{+}(u_0)+\mathcal{U}_{0,0}^{-}(u_0)=\frac{1}{2}u_0(x+t)-\frac{1}{2}u_0(t-x)$,

(11.14) $\mathcal{U}_{1,0}^{+}(u_0)+\mathcal{U}_{1,0}^{-}(u_0)=\frac{1}{2}\int_{y=t-x}^{x+t}u_1(y)dy$,

(11.15) $\mathcal{G}_0(g_0)=g_0(t-x)$,

(11.16) $\mathcal{F}_0^{+}(f)+\mathcal{F}_0^{-}(f)=\frac{1}{2}\int_{\tau=0}^{t-x}\left[\int_{y=(t-x)-\tau}^{(t+x)-\tau}f(y,\tau)dy\right]d\tau+\frac{1}{2}\int_{\tau=t-x}^{t}\left[\int_{y=(x-t)+\tau}^{(x+t)-\tau}f(y,\tau)dy\right]d\tau$.

Indeed, the term $-\frac{1}{2}u_0(t-x)$ in (11.13), this time, comes from the integral $\mathcal{U}_{0,0}^{-}(u_0)$. The proof of (11.14) is similar to the proof of (11.9), while (11.15) follows from the formula

$$\int_{-\infty}^{\infty}[e^{i\lambda(x-t)}\int_{\tau=0}^{t}e^{i\lambda\tau}g_0(\tau)d\tau]d\lambda=\int_{-\infty}^{\infty}[e^{i\lambda(x-t)}\int_{\tau=0}^{t}e^{-i\lambda(-\tau)}g_0(\tau)d\tau]d\lambda=\int_{-\infty}^{\infty}[e^{i\lambda(x-t)}\int_{\tau=-t}^{0}e^{-i\lambda\tau}g_0(-\tau)d\tau]d\lambda=g_0(t-x).$$

*Proof of (11.16)* Let us define $f(x,t)=0$, for $x<0$. As in the proof of (11.11), with $x<t$, we compute:

$$\begin{aligned}\mathcal{F}_0^{+}(f)=&\frac{1}{2}\int_{\tau=0}^{t-x}\left[\int_{\lambda=-\infty}^{\infty}\{e^{i\lambda[x+(t-\tau)]}-e^{i\lambda[x-(t-\tau)]}\}\frac{\hat{f}(\lambda,\tau)}{i\lambda}d\lambda\right]d\tau\\&+\frac{1}{2}\int_{\tau=t-x}^{t}\left[\int_{\lambda=-\infty}^{\infty}\{e^{i\lambda[x+(t-\tau)]}-e^{i\lambda[x-(t-\tau)]}\}\frac{\hat{f}(\lambda,\tau)}{i\lambda}d\lambda\right]d\tau\\=&\frac{1}{2}\int_{\tau=0}^{t-x}\left[\int_{y=x-(t-\tau)}^{x+(t-\tau)}f(y,\tau)dy\right]d\tau+\frac{1}{2}\int_{\tau=t-x}^{t}\left[\int_{y=x-(t-\tau)}^{x+(t-\tau)}f(y,\tau)dy\right]d\tau\\=&\frac{1}{2}\int_{\tau=0}^{t-x}\left[\int_{y=0}^{x+(t-\tau)}f(y,\tau)dy\right]d\tau+\frac{1}{2}\int_{\tau=t-x}^{t}\left[\int_{y=x-(t-\tau)}^{x+(t-\tau)}f(y,\tau)dy\right]d\tau\end{aligned}$$

and

$$\mathcal{F}_0^{-}(f)=\frac{1}{2}\int_{\tau=0}^{t-x}\left[\int_{\lambda=-\infty}^{\infty}\{e^{i\lambda[x-(t-\tau)]}-e^{i\lambda[x+(t-\tau)]}\}\frac{\hat{f}(-\lambda,\tau)}{i\lambda}d\lambda\right]d\tau$$

$$+\frac{1}{2}\int_{\tau=t-x}^{t}\left[\int_{\lambda=-\infty}^{\infty}\{e^{i\lambda[x-(t-\tau)]}-e^{i\lambda[x+(t-\tau)]}\}\frac{\hat{f}(-\lambda,\tau)}{i\lambda}d\lambda\right]d\tau=-\frac{1}{2}\int_{\tau=0}^{t-x}\left[\int_{y=0}^{-x+t-\tau}f(y,\tau)dy\right]d\tau,$$

and (11.16) follows.

Thus, equations (11.8) – (11.11) give (11.6), while equations (11.13) – (11.16) give (11.7). On the other hand, independently of the above computations and formula (11.4), it is straightforward to verify that, indeed, (11.6) and (11.7) solve problem (11.1).

## 12. The behaviour of the solution, as $t\to\infty$, with periodic data

In this section, we consider the solution (1.5) of problem (1.1), in the case the data $g_0(t)$ and $f(x,t)$ are periodic in $t$. More precisely, assuming, in addition to (1.2), that, for a fixed $T>0$,

(12.1) $$g_0(t+T)=g_0(t) \text{ and } f(x,t)=f(x,t+T) \textit{ for every } t\ge 0 \textit{ and } x\ge 0,$$

we will prove the following theorem.

**Theorem 7** *For the solution $u(x,t)$, of problem (1.1) given by (1.5), we have*

(12.2) $$u(x,t+T)-u(x,t)=\mathrm{O}(1/t), \textit{ as } t\to\infty,$$

*uniformly for $x$ in compact subsets of $[0,\infty)$, provided that $g_0$ and $f$ satisfy (12.1).*

**Proof** We will use the notation of (1.6). Also, throughout this proof, $t>>x$.

***Step 1*** We claim that

(12.3) $$\mathcal{G}(g_0)(x,t+T)-\mathcal{G}(g_0)(x,t)=\mathrm{O}(1/t), \textit{ as } t\to\infty, \textit{ uniformly for } x \textit{ in compact subsets of } [0,\infty).$$

First, let us note that, in view of (12.1),

(12.4) $$e^{-\omega(\lambda)(t+T)}\tilde{g}_0(\omega(\lambda),t+T)-e^{-\omega(\lambda)t}\tilde{g}_0(\omega(\lambda),t)$$
$$=e^{-\omega(\lambda)(t+T)}\int_{\tau=0}^{t+T}e^{\omega(\lambda)\tau}g_0(\tau)d\tau-e^{-\omega(\lambda)t}\int_{\tau=0}^{t}e^{\omega(\lambda)\tau}g_0(\tau)d\tau=e^{-\omega(\lambda)t}\int_{\tau=-T}^{0}e^{\omega(\lambda)\tau}g_0(\tau)d\tau.$$

It follows that

(12.5) $$\mathcal{G}(g_0)(x,t+T)-\mathcal{G}(g_0)(x,t)=\int_{-1/2}^{1/2}\left[e^{i\lambda x-\omega_1(\lambda)t}\int_{\tau=-T}^{0}e^{\omega_1(\lambda)\tau}g_0(\tau)d\tau-e^{i\lambda x-\omega_2(\lambda)t}\int_{\tau=-T}^{0}e^{\omega_2(\lambda)\tau}g_0(\tau)d\tau\right]\frac{\lambda d\lambda}{\rho(\lambda)}$$
$$+\left(\int_{-1}^{-1/2}+\int_{1/2}^{1}\right)\left[e^{i\lambda x-\omega_1(\lambda)t}\int_{\tau=-T}^{0}e^{\omega_1(\lambda)\tau}g_0(\tau)d\tau-e^{i\lambda x-\omega_2(\lambda)t}\int_{\tau=-T}^{0}e^{\omega_2(\lambda)\tau}g_0(\tau)d\tau\right]\frac{\lambda d\lambda}{\rho(\lambda)}$$
$$+\int_{\mathbb{L}}\left[e^{i\lambda x-\omega_1(\lambda)t}\int_{\tau=-T}^{0}e^{\omega_1(\lambda)\tau}g_0(\tau)d\tau\right]\frac{\lambda d\lambda}{\rho(\lambda)}-\int_{\mathbb{L}}\left[-e^{i\lambda x-\omega_2(\lambda)t}\int_{\tau=-T}^{0}e^{\omega_2(\lambda)\tau}g_0(\tau)d\tau\right]\frac{\lambda d\lambda}{\rho(\lambda)}:=J_1+J_2+J_3+J_4.$$

($J_1, J_2, J_3, J_4$ are the 1st, 2nd, 3rd and 4th integrals in (12.5), respectively.)

In order to deal with the 2nd integral in the RHS of (12.5), namely $J_2$, let us note that, as in the proof of Lemma 1,

(12.6) $$\left[\int_{\tau=-T}^{0}e^{-\omega_1(\lambda)(t-\tau)}g_0(\tau)d\tau-\int_{\tau=-T}^{0}e^{-\omega_2(\lambda)(t-\tau)}g_0(\tau)d\tau\right]\frac{1}{\rho(\lambda)}=e^{-t/2}\int_{\tau=-T}^{0}e^{\tau/2}\frac{e^{-i\rho(\lambda)(t-\tau)}-e^{i\rho(\lambda)(t-\tau)}}{\rho(\lambda)}g_0(\tau)d\tau$$
$$=-2ie^{-t/2}\int_{\tau=-T}^{0}e^{\tau/2}\frac{\sin[\rho(\lambda)(t-\tau)]}{\rho(\lambda)}g_0(\tau)d\tau,$$

which implies that

(12.7) $$|J_2(x,t)|\le 2(t+T)e^{-t/2}\left(\int_{-1}^{-1/2}+\int_{1/2}^{1}\right)\left[\left\{\int_{\tau=-T}^{0}e^{\tau/2}\left|\frac{\sin[\rho(\lambda)(t-\tau)]}{\rho(\lambda)(t-\tau)}\right||g_0(\tau)|d\tau\right\}\lambda d\lambda\right]=\mathrm{O}(te^{-t/2}),$$

since for $\lambda \in [-1,-\frac{1}{2}] \cup [\frac{1}{2},1]$, we have $\rho(\lambda) \in \mathbb{R}$ and, therefore,

(12.8) $$\sup\left\{ \frac{|\sin[\rho(\lambda)(t-\tau)]|}{|\rho(\lambda)(t-\tau)|} : \lambda \in [-1,-\tfrac{1}{2}] \cup [\tfrac{1}{2},1],\ \tau \in [-T,0] \ and\ t \geq 0 \right\} < \infty .$$

The 1st integral in the RHS of (12.5) is equal to

$$J_1(x,t) = -2ie^{-t/2} \int_{-1/2}^{1/2} \left\{ e^{i\lambda x} \int_{\tau=-T}^{0} \left[ e^{\tau/2} \sum_{n=1}^{\infty} \frac{(-1)^{n-1}}{(2n-1)!} (t-\tau)^{2n-1} \left( \lambda^2 - \frac{1}{4} \right)^{n-1} \right] g_0(\tau) d\tau \right\} \lambda d\lambda .$$

Indeed, the proof of this formula is similar to the proof of (2.11).

Therefore,

(12.9) $$|J_1(x,t)| \leq 2e^{-t/2} \int_{0}^{1/2} \left\{ \int_{\tau=-T}^{0} \left[ e^{\tau/2} \sum_{n=1}^{\infty} \frac{1}{(2n-1)!} (t-\tau)^{2n-1} \left( \frac{1}{4} - \lambda^2 \right)^{n-1} \right] g_0(\tau) d\tau \right\} d(\lambda^2)$$

$$= 4e^{-t/2} \int_{\tau=-T}^{0} \left[ \frac{1}{t-\tau} e^{\tau/2} \sum_{n=1}^{\infty} \frac{1}{(2n)!} \left( \frac{t-\tau}{2} \right)^{2n} \right] g_0(\tau) d\tau = \mathrm{O}(1/t) ,\ as\ t \to \infty ,$$

where we used also the identity

(12.10) $$\sum_{n=1}^{\infty} \frac{1}{(2n)!} \left( \frac{t-\tau}{2} \right)^{2n} = \frac{1}{2} [e^{(t-\tau)/2} + e^{-(t-\tau)/2}] - 1 .$$

Since

(12.11) $$\int_{\tau=-T}^{0} e^{\omega_1(\lambda)\tau} g_0(\tau) d\tau = \frac{g_0(0)}{\omega_1(\lambda)} - \frac{e^{-\omega_1(\lambda)T} g_0(-T)}{\omega_1(\lambda)}$$

$$- \frac{g_0^{(1)}(0)}{[\omega_1(\lambda)]^2} + \frac{e^{-\omega_1(\lambda)T} g_0^{(1)}(-T)}{[\omega_1(\lambda)]^2} + \frac{1}{[\omega_1(\lambda)]^2} \int_{\tau=-T}^{0} e^{\omega_1(\lambda)\tau} g_0^{(2)}(\tau) d\tau ,$$

the 3rd integral in the RHS of (12.5) becomes

(12.12) $$J_3(x,t) = \int_{\mathbb{L}} \left[ e^{i\lambda x - \omega_1(\lambda)t} \int_{\tau=-T}^{0} e^{\omega_1(\lambda)\tau} g_0(\tau) d\tau \right] \frac{\lambda d\lambda}{\rho(\lambda)}$$

$$= \int_{-\kappa^-} e^{i\lambda x - \omega_1(\lambda)t} \frac{g_0(0)}{\omega_1(\lambda)} \frac{\lambda d\lambda}{\rho(\lambda)} - \int_{-\kappa^-} e^{i\lambda x - \omega_1(\lambda)t} \frac{e^{-\omega_1(\lambda)T} g_0(-T)}{\omega_1(\lambda)} \frac{\lambda d\lambda}{\rho(\lambda)}$$

$$- \int_{-\kappa^-} e^{i\lambda x - \omega_1(\lambda)t} \frac{g_0^{(1)}(0)}{[\omega_1(\lambda)]^2} \frac{\lambda d\lambda}{\rho(\lambda)} + \int_{-\kappa^-} e^{i\lambda x - \omega_1(\lambda)t} \frac{e^{-\omega_1(\lambda)T} g_0^{(1)}(-T)}{[\omega_1(\lambda)]^2} \frac{\lambda d\lambda}{\rho(\lambda)}$$

$$+ \int_{\mathbb{L}} \left\{ e^{i\lambda x - \omega_1(\lambda)t} \frac{1}{[\omega_1(\lambda)]^2} \int_{\tau=-T}^{0} e^{\omega_1(\lambda)\tau} g_0^{(2)}(\tau) d\tau \right\} \frac{\lambda d\lambda}{\rho(\lambda)} .$$

Indeed, this follows from from (3.5) and Jordan's lemma.

Now, using (1) of (12.13) (below),

(12.13) $$(1)\ \sin\left( \frac{\theta_1(\lambda) + \theta_2(\lambda)}{2} \right) \leq 0\ for\ \lambda \in \kappa^- ,\ and\ (2)\ \sin\left( \frac{\theta_1(\lambda) + \theta_2(\lambda)}{2} \right) \geq 0\ for\ \lambda \in \kappa^+ ,$$

we see, also in view of (3.5), that

(12.14) *the integrals in the RHS of (12.12), taken over the semicircle* $-\kappa^-$, *are* $\mathrm{O}(e^{-t/2})$, *as* $t \to \infty$.

On the other hand, since $\left| e^{i\lambda x - \omega_1(\lambda)t} \right| = e^{-t/2}$ for $\lambda \in \mathbb{L}$, we have

(12.15) $$\left| \int_{\mathbb{L}} \left\{ e^{i\lambda x - \omega_1(\lambda)t} \frac{1}{[\omega_1(\lambda)]^2} \int_{\tau=-T}^{0} e^{\omega_1(\lambda)\tau} g_0^{(2)}(\tau) d\tau \right\} \frac{\lambda d\lambda}{\rho(\lambda)} \right|$$

$$\le e^{-t/2}\left[\int_{\tau=-T}^{0} e^{\tau/2}\left|g_0^{(2)}(\tau)\right|d\tau\right]\left(\int_{\mathbb{L}}\frac{1}{|\omega_1(\lambda)|^2}\frac{|\lambda|d\lambda}{|\rho(\lambda)|}\right)=\mathrm{O}(e^{-t/2}).$$

It follows from (12.14) and (12.15) that

(12.16) $J_3(x,t)=\mathrm{O}(e^{-t/2})$, *as* $t\to\infty$, *uniformly for* $x$ *in compact subsets of* $[0,\infty)$.

Finally, writing

$$(12.17)\qquad J_4=\int_{\mathbb{L}}\left[-e^{i\lambda x-\omega_2(\lambda)t}\int_{\tau=-T}^{0}e^{\omega_2(\lambda)\tau}g_0(\tau)d\tau\right]\frac{\lambda d\lambda}{\rho(\lambda)}$$

$$=\int_{\kappa^+}e^{i\lambda x-\omega_2(\lambda)t}\frac{g_0(0)}{\omega_2(\lambda)}\frac{\lambda d\lambda}{\rho(\lambda)}-\int_{\kappa^+}e^{i\lambda x-\omega_2(\lambda)t}\frac{e^{-\omega_2(\lambda)T}g_0(-T)}{\omega_2(\lambda)}\frac{\lambda d\lambda}{\rho(\lambda)}$$

$$-\int_{\kappa^+}e^{i\lambda x-\omega_2(\lambda)t}\frac{g_0^{(1)}(0)}{[\omega_2(\lambda)]^2}\frac{\lambda d\lambda}{\rho(\lambda)}+\int_{\kappa^+}e^{i\lambda x-\omega_2(\lambda)t}\frac{e^{-\omega_2(\lambda)T}g_0^{(1)}(-T)}{[\omega_2(\lambda)]^2}\frac{\lambda d\lambda}{\rho(\lambda)}$$

$$+\int_{\mathbb{L}}\left\{e^{i\lambda x-\omega_2(\lambda)t}\frac{1}{[\omega_2(\lambda)]^2}\int_{\tau=-T}^{0}e^{\omega_2(\lambda)\tau}g_0^{(2)}(\tau)d\tau\right\}\frac{\lambda d\lambda}{\rho(\lambda)}$$

and working as previously, we see, also in view of (3.6) and (2) of (12.13), that

(12.18) $J_4(x,t)=\mathrm{O}(e^{-t/2})$, *as* $t\to\infty$, *uniformly for* $x$ *in compact subsets of* $[0,\infty)$.

Thus, (12.3) follows from (12.5), (12.7), (12.9), (12.16) and (12.18).

***Step 2*** We claim that

(12.19) $\mathcal{U}_1^+(u_1)(x,t)+\mathcal{U}_1^-(u_1)(x,t)=\mathrm{O}(1/t)$, *as* $t\to\infty$, *uniformly for* $x$ *in compact subsets of* $[0,\infty)$.

Changing the variable in the part $\int_{-1/2}^{1/2}\cdots d\lambda$ of the integral $\mathcal{U}_1^-(u_1)(x,t)$, by setting $\mu=-\lambda$, we obtain that

$$(12.20)\quad \mathcal{U}_1^+(u_1)(x,t)+\mathcal{U}_1^-(u_1)(x,t)=\int_{-1/2}^{1/2}(e^{i\lambda x}-e^{-i\lambda x})[e^{-\omega_2(\lambda)t}-e^{-\omega_1(\lambda)t}]\hat{u}_1(\lambda)\frac{d\lambda}{2i\rho(\lambda)}$$

$$+\left(\int_{-1}^{-1/2}+\int_{1/2}^{1}\right)[e^{i\lambda x-\omega_2(\lambda)t}-e^{i\lambda x-\omega_1(\lambda)t}][\hat{u}_1(\lambda)-\hat{u}_1(-\lambda)]\frac{d\lambda}{2i\rho(\lambda)}$$

$$+\int_{\mathbb{L}}[e^{i\lambda x-\omega_2(\lambda)t}-e^{i\lambda x-\omega_1(\lambda)t}][\hat{u}_1(\lambda)-\hat{u}_1(-\lambda)]\frac{d\lambda}{2i\rho(\lambda)}:=I_{[-\frac{1}{2},\frac{1}{2}]}+I_{(-1,-\frac{1}{2}]\cup[\frac{1}{2},1)}+I_{\mathbb{L}}.$$

Now, in view of (2.9),

$$(12.21)\qquad \left|I_{[-\frac{1}{2},\frac{1}{2}]}(x,t)\right|=e^{-t/2}\left|\int_{-1/2}^{1/2}\sin(\lambda x)[e^{i\rho(\lambda)t}-e^{-i\rho(\lambda)t}]\hat{u}_1(\lambda)\frac{d\lambda}{\rho(\lambda)}\right|$$

$$=e^{-t/2}\left|\int_{-1/2}^{1/2}\sin(\lambda x)\left[\sum_{n=1}^{\infty}\frac{(-1)^{n-1}}{(2n-1)!}t^{2n-1}\left(\lambda^2-\frac{1}{4}\right)^{n-1}\right]\hat{u}_1(\lambda)d\lambda\right|$$

$$\le e^{-t/2}\int_{-1/2}^{1/2}\left|\frac{\sin(\lambda x)}{\lambda}\right|\left[\sum_{n=1}^{\infty}\frac{t^{2n-1}}{(2n-1)!}\left(\frac{1}{4}-\lambda^2\right)^{n-1}\right]|\hat{u}_1(\lambda)||\lambda|d\lambda$$

$$\le 2e^{-t/2}\sup_{-\frac{1}{2}\le\lambda\le\frac{1}{2}}\left[\left|\frac{\sin(\lambda x)}{\lambda}\right||\hat{u}_1(\lambda)|\right]\int_0^{1/2}\left[\sum_{n=1}^{\infty}\frac{t^{2n-1}}{(2n-1)!}\left(\frac{1}{4}-\lambda^2\right)^{n-1}\right]\lambda d\lambda=\mathrm{O}(1/t),$$

where we used also formula (12.10) with $\tau=0$.

On the other hand, again in view of (2.9),

$$(12.22)\qquad \left|I_{(-1,-\frac{1}{2}]\cup[\frac{1}{2},1)}(x,t)\right|\le e^{-t/2}\left(\int_{-1}^{-1/2}+\int_{1/2}^{1}\right)\left|\frac{e^{i\rho(\lambda)t}-e^{-i\rho(\lambda)t}}{2\rho(\lambda)}\right||\hat{u}_1(\lambda)-\hat{u}_1(-\lambda)|d\lambda$$

$$\leq te^{-t/2}\left(\int_{-1}^{-1/2}+\int_{1/2}^{1}\right)\left|\frac{\sin[\rho(\lambda)t]}{t\rho(\lambda)}\right|\left|\hat{u}_1(\lambda)-\hat{u}_1(-\lambda)\right|d\lambda=\mathrm{O}(te^{-t/2}),$$

where we used also the fact that, since $\rho(\lambda)\in\mathbb{R}$ for $\lambda\in\mathbb{R}-[-\frac{1}{2},\frac{1}{2}]$,

$$\sup\left\{\frac{|\sin[\rho(\lambda)t]|}{|\rho(\lambda)t|}:\ \lambda\in[-1,-\tfrac{1}{2}]\cup[\tfrac{1}{2},1]\ \textit{and}\ t\geq 0\right\}<\infty .$$

Also, keeping in mind that the integral $I_{\mathbb{L}}$ converges absolutely, we have

$$\left|I_{\mathbb{L}}(x,t)\right|\leq e^{-t/2}\left(\int_{-\infty}^{-1}+\int_{1}^{\infty}\right)\left|\frac{\sin[\rho(\lambda)t]}{\rho(\lambda)}\right|\left|\hat{u}_1(\lambda)-\hat{u}_1(-\lambda)\right|d\lambda=\mathrm{O}(e^{-t/2}). \tag{12.23}$$

Thus, (12.19) follows from (12.20) – (12.23).

***Step 3*** We claim that

$$\mathcal{U}_0^{+}(u_0)(x,t)+\mathcal{U}_0^{-}(u_0)(x,t)=\mathrm{O}(1/t),\ \textit{as}\ t\to\infty,\ \textit{uniformly for}\ x\ \textit{in compact subsets of}\ [0,\infty). \tag{12.24}$$

We have

$$\begin{aligned}\mathcal{U}_0^{+}(u_0)(x,t)+\mathcal{U}_0^{-}(u_0)(x,t)=&\int_{-1/2}^{1/2}(e^{i\lambda x}-e^{-i\lambda x})[\omega_1(\lambda)e^{-\omega_2(\lambda)t}-\omega_2(\lambda)e^{-\omega_1(\lambda)t}]\hat{u}_0(\lambda)\frac{d\lambda}{2i\rho(\lambda)}\\&+\int_{-1}^{-1/2}(e^{i\lambda x}+e^{-i\lambda x})[\omega_1(\lambda)e^{-\omega_2(\lambda)t}-\omega_2(\lambda)e^{-\omega_1(\lambda)t}]\hat{u}_0(\lambda)\frac{d\lambda}{2i\rho(\lambda)}\\&+\int_{\mathbb{L}}[\omega_1(\lambda)e^{i\lambda x-\omega_2(\lambda)t}-\omega_2(\lambda)e^{i\lambda x-\omega_1(\lambda)t}]\hat{u}_0(\lambda)\frac{d\lambda}{2i\rho(\lambda)}\\&+\int_{\mathbb{L}}[\omega_2(\lambda)e^{i\lambda x-\omega_1(\lambda)t}-\omega_1(\lambda)e^{i\lambda x-\omega_2(\lambda)t}]\hat{u}_0(-\lambda)\frac{d\lambda}{2i\rho(\lambda)}.\end{aligned} \tag{12.25}$$

Since this time, the integrals over $\mathbb{L}$ do not converge absolutely (in general), we have to deform their contours as follows:

$$\int_{\mathbb{L}}\omega_1(\lambda)e^{i\lambda x-\omega_2(\lambda)t}\hat{u}_0(\lambda)\frac{d\lambda}{2i\rho(\lambda)}=u_0(0)\int_{\kappa^+}\omega_1(\lambda)e^{i\lambda x-\omega_2(\lambda)t}\frac{1}{i\lambda}\frac{d\lambda}{2i\rho(\lambda)}+\int_{\mathbb{L}}\omega_1(\lambda)e^{i\lambda x-\omega_2(\lambda)t}\frac{1}{i\lambda}(u_0^{(1)})\hat{}(\lambda)\frac{d\lambda}{2i\rho(\lambda)},$$

$$\int_{\mathbb{L}}\omega_2(\lambda)e^{i\lambda x-\omega_1(\lambda)t}\hat{u}_0(\lambda)\frac{d\lambda}{2i\rho(\lambda)}=\int_{-\kappa^-}\omega_2(\lambda)e^{i\lambda x-\omega_1(\lambda)t}\hat{u}_0(\lambda)\frac{d\lambda}{2i\rho(\lambda)},$$

$$\begin{aligned}\int_{\mathbb{L}}\omega_2(\lambda)e^{i\lambda x-\omega_1(\lambda)t}\hat{u}_0(-\lambda)\frac{d\lambda}{2i\rho(\lambda)}=&\,u_0(0)\int_{-\kappa^-}\omega_2(\lambda)e^{i\lambda x-\omega_1(\lambda)t}\frac{1}{(-i\lambda)}\frac{d\lambda}{2i\rho(\lambda)}\\&+\int_{\mathbb{L}}\omega_2(\lambda)e^{i\lambda x-\omega_1(\lambda)t}\frac{1}{(-i\lambda)}(u_0^{(1)})\hat{}(-\lambda)\frac{d\lambda}{2i\rho(\lambda)},\end{aligned}$$

$$\int_{\mathbb{L}}\omega_1(\lambda)e^{i\lambda x-\omega_2(\lambda)t}\hat{u}_0(-\lambda)\frac{d\lambda}{2i\rho(\lambda)}=\int_{\kappa^+}\omega_1(\lambda)e^{i\lambda x-\omega_2(\lambda)t}\hat{u}_0(-\lambda)\frac{d\lambda}{2i\rho(\lambda)}.$$

Now, using the above formulas, it is easy to see that the integrals over $\mathbb{L}$, in the RHS of (12.25), are $\mathrm{O}(e^{-t/2})$.
On the other hand, working as in the proof of (12.21) and using (2.10), we see that the integral in (12.25), which is taken over the interval $[-\frac{1}{2},\frac{1}{2}]$, is $\mathrm{O}(1/t)$, while, as in the proof of (12.22) and in view of (2.10) again, we see that the integral in (12.25), taken over the interval $[-1,-\frac{1}{2}]$, is $\mathrm{O}(te^{-t/2})$.
Combining the above conclusions, we see that (12.24) follows from (12.25).

***Step 4*** Setting

$$\mathcal{H}(f)(x,t):=\mathcal{F}^{+}(f)(x,t)+\mathcal{F}^{-}(f)(x,t),$$

we claim that

(12.26) $\mathcal{H}(f)(x,t+T)-\mathcal{H}(f)(x,t)=\mathrm{O}(1/t)$, *as* $t\to\infty$, *uniformly for* $x$ *in compact subsets of* $[0,\infty)$.

By the periodicity assumptions on $f(x,\cdot)$, we have $\hat{f}(\lambda,t+T)=\hat{f}(\lambda,t)$ and

$$e^{-\omega(t+T)}\tilde{\hat{f}}(\lambda,\omega(\lambda),t+T)-e^{-\omega t}\tilde{\hat{f}}(\lambda,\omega(\lambda),t)=e^{-\omega t}\int_{\tau=-T}^{0}e^{\omega(\lambda)\tau}\hat{f}(\lambda,\tau)d\tau . \quad (12.27)$$

Working as in *step 2* and in view of (12.27), we have

$$\begin{aligned}(12.28)\quad &\mathcal{H}(f)(x,t+T)-\mathcal{H}(f)(x,t)\\ &=\int_{-1/2}^{1/2}\left\{\left(e^{i\lambda x}-e^{-i\lambda x}\right)\left[e^{-\omega_2(\lambda)t}\int_{\tau=-T}^{0}e^{\omega_2(\lambda)\tau}\hat{f}(\lambda,\tau)d\tau-e^{-\omega_1(\lambda)t}\int_{\tau=-T}^{0}e^{\omega_1(\lambda)\tau}\hat{f}(\lambda,\tau)d\tau\right]\right\}\frac{d\lambda}{2i\rho(\lambda)}\\ &+\left(\int_{-\infty}^{-1/2}+\int_{1/2}^{\infty}\right)\left[e^{i\lambda x-\omega_2(\lambda)t}\int_{\tau=-T}^{0}e^{\omega_2(\lambda)\tau}\hat{f}(\lambda,\tau)d\tau-e^{i\lambda x-\omega_1(\lambda)t}\int_{\tau=-T}^{0}e^{\omega_1(\lambda)\tau}\hat{f}(\lambda,\tau)d\tau\right]\frac{d\lambda}{2i\rho(\lambda)}\\ &+\left(\int_{-\infty}^{-1/2}+\int_{1/2}^{\infty}\right)\left[e^{i\lambda x-\omega_1(\lambda)t}\int_{\tau=-T}^{0}e^{\omega_1(\lambda)\tau}\hat{f}(-\lambda,\tau)d\tau-e^{i\lambda x-\omega_2(\lambda)t}\int_{\tau=-T}^{0}e^{\omega_2(\lambda)\tau}\hat{f}(-\lambda,\tau)d\tau\right]\frac{d\lambda}{2i\rho(\lambda)}.\end{aligned}$$

Using formulas similar to (2.12) and (2.13), we see that, from the above integrals, the one taken over the interval $[-\frac{1}{2},\frac{1}{2}]$ is equal to

$$e^{-t/2}\int_{\tau=-T}^{0}\left\{\int_{-1/2}^{1/2}\sin(\lambda x)\left[e^{\tau/2}\sum_{n=1}^{\infty}\frac{1}{(2n-1)!}(t-\tau)^{2n-1}\left(\frac{1}{4}-\lambda^2\right)^{n-1}\right]\hat{f}(\lambda,\tau)d\lambda\right\}d\tau ,$$

which is $\mathrm{O}(1/t)$. Indeed, the proof of the last assertion is similar to the proof of (12.21), using again (12.8).
On the other hand, the integrals in (12.27), which are taken over $(-\infty,-\frac{1}{2}]\cup[\frac{1}{2},\infty)$, converge absolutely and are easily seen to be $\mathrm{O}(te^{-t/2})$, as in (12.22).
Therefore, (12.28) implies (12.26).

*Completion of the proof of the theorem* Combining the results of *steps 1 – 4*, we obtain the estimate (12.2).

## 13. The behaviour of the solution, as $t\to\infty$, with weakly periodic data

In this section we will prove the following extension of *Theorem 7*.

**Theorem 8** *Let* $T>0$ *be fixed. Suppose that, in addition to (1.2),* $g_0(t)$ *and* $f(x,t)$ *are weakly* $T-$*periodic in the following sense:*

(13.1) (1) $\int_0^\infty\left|g_0^{(j)}(t+T)-g_0^{(j)}(t)\right|dt<+\infty$, for $j=0,1,2$,

(2) $\lim_{t\to\infty}\left|g_0(t+T)-g_0(t)\right|=0$, (3) $\sup_{t\geq 0}\left|g_0^{(1)}(t+T)-g_0^{(1)}(t)\right|<+\infty$,

(13.2) (1) $\int_0^\infty\int_0^\infty\left|f(x,t+T)-f(x,t)\right|dtdx<+\infty$ , (2) $\int_0^\infty\int_0^\infty\left|f_x(x,t+T)-f_x(x,t)\right|dtdx<+\infty$,

(3) $\int_0^\infty\left|f(0,t+T)-f(0,t)\right|dt<+\infty$.

*Then the solution* (1.5) *of problem* (1.1) *satisfies the following*:

$$\lim_{t\to\infty}[u(x,t+T)-u(x,t)]=0, \quad (13.3)$$

*uniformly for* $x$ *in compact subsets of* $[0,\infty)$.

**Proof** Firstly, in view of the proof of *Theorem 7*, it suffices to deal only with the parts of the solution which involve the data $g_0(t)$ and $f(x,t)$.

Secondly, in this case the LHSs of (12.4) and (12.27) become

$$(13.4)\qquad e^{-\omega(t+T)}\tilde{g}_0(\omega(\lambda),t+T)-e^{-\omega t}\tilde{g}_0(\omega(\lambda),t)=e^{-\omega(t+T)}\int_0^{t+T}e^{\omega(\lambda)\tau}g_0(\tau)d\tau-e^{-\omega t}\int_0^t e^{\omega(\lambda)\tau}g_0(\tau)d\tau$$

$$=e^{-\omega t}\int_0^t e^{\omega(\lambda)\tau}[g_0(\tau+T)-g_0(\tau)]d\tau+e^{-\omega t}\int_{-T}^0 e^{\omega(\lambda)\tau}g_0(\tau+T)d\tau$$

and

$$(13.5)\qquad e^{-\omega(t+T)}\tilde{\hat{f}}(\lambda,\omega(\lambda),t+T)-e^{-\omega t}\tilde{\hat{f}}(\lambda,\omega(\lambda),t)$$

$$=e^{-\omega t}\int_0^t e^{\omega(\lambda)\tau}[\hat{f}(\lambda,\tau+T)-\hat{f}(\lambda,\tau)]d\tau+e^{-\omega t}\int_{-T}^0 e^{\omega(\lambda)\tau}\hat{f}(\lambda,\tau+T)d\tau\,.$$

Also, again in view of the proof of *Theorem 7*, we see that the parts of $u(x,t+T)-u(x,t)$, which involve the last terms in (13.4) and (13.5), namely the integrals of the form

$$e^{-\omega t}\int_{-T}^0 e^{\omega(\lambda)\tau}g_0(\tau+T)d\tau\,,\ e^{-\omega t}\int_{-T}^0 e^{\omega(\lambda)\tau}\hat{f}(\pm\lambda,\tau+T)d\tau\,,$$

tend to zero, as $t\to\infty$.

***Step 1*** We claim that

$$(13.6)\qquad \lim_{t\to\infty}\int_{-\infty}^{\infty}\left\{e^{i\lambda x-\omega_1(\lambda)t}\int_0^t e^{\omega_1(\lambda)\tau}[g_0(\tau+T)-g_0(\tau)]d\tau-e^{i\lambda x-\omega_2(\lambda)t}\int_0^t e^{\omega_2(\lambda)\tau}[g_0(\tau+T)-g_0(\tau)]d\tau\right\}\frac{\lambda d\lambda}{\rho(\lambda)}=0\,,$$

uniformly for $x$ in compact subsets of $[0,\infty)$.

To prove this, we write the integral in (13.6) as follows:

$$(13.7)\qquad \int_{-1/2}^{1/2}\left\{e^{i\lambda x}\int_0^{\infty}\left(\chi_{[0,t]}(\tau)\left[\frac{e^{-\omega_1(\lambda)(t-\tau)}-e^{-\omega_2(\lambda)(t-\tau)}}{\rho(\lambda)}\right]\varpi(\tau)\right)d\tau\right\}\lambda d\lambda$$

$$+\left(\int_{-1}^{-1/2}+\int_{1/2}^{1}\right)\left\{e^{i\lambda x}\int_0^{\infty}\left(\chi_{[0,t]}(\tau)\left[\frac{e^{-\omega_1(\lambda)(t-\tau)}-e^{-\omega_2(\lambda)(t-\tau)}}{\rho(\lambda)}\right]\varpi(\tau)\right)d\tau\right\}\lambda d\lambda$$

$$+\int_{\mathbb{L}}\left\{e^{i\lambda x-\omega_1(\lambda)t}\int_0^{\infty}\chi_{[0,t]}(\tau)e^{\omega_1(\lambda)\tau}\varpi(\tau)d\tau\right\}\frac{\lambda d\lambda}{\rho(\lambda)}$$

$$-\int_{\mathbb{L}}\left\{e^{i\lambda x-\omega_2(\lambda)t}\int_0^{\infty}\chi_{[0,t]}(\tau)e^{\omega_2(\lambda)\tau}\varpi(\tau)d\tau\right\}\frac{\lambda d\lambda}{\rho(\lambda)}:=\mathrm{I}_1+\mathrm{I}_2+\mathrm{I}_3+\mathrm{I}_4\,,$$

where we have set $\varpi(t):=g_0(t+T)-g_0(t)$. ($\mathrm{I}_1,\mathrm{I}_2,\mathrm{I}_3,\mathrm{I}_4$ are the 1st, 2nd, 3rd and 4th integrals in (13.7), respectively.)

Now, we recall the formula

$$\frac{e^{-\omega_1(\lambda)(t-\tau)}-e^{-\omega_2(\lambda)(t-\tau)}}{\rho(\lambda)}=-2ie^{-(t-\tau)/2}\frac{1}{t-\tau}\sum_{n=1}^{\infty}\frac{(-1)^{n-1}}{(2n-1)!}\left(\lambda^2-\frac{1}{4}\right)^{n-1}(t-\tau)^{2n}\,,$$

hence (using also (12.10))

$$\left|\int_{-1/2}^{1/2}e^{i\lambda x}\left[\frac{e^{-\omega_1(\lambda)(t-\tau)}-e^{-\omega_2(\lambda)(t-\tau)}}{\rho(\lambda)}\right]\lambda d\lambda\right|\le\left|e^{-(t-\tau)/2}\frac{1}{t-\tau}[e^{(t-\tau)/2}+e^{-(t-\tau)/2}-2]\right|$$

and, therefore,

$$\left|\mathrm{I}_1(x,t)\right|\le\int_{\tau=0}^{\infty}\left|e^{-(t-\tau)/2}\frac{1}{t-\tau}[e^{(t-\tau)/2}+e^{-(t-\tau)/2}-2]\right|\chi_{[0,t]}(\tau)\left|\varpi(\tau)\right|d\tau\,.$$

Since

$$\lim_{t\to\infty}\left|e^{-(t-\tau)/2}\frac{1}{t-\tau}[e^{(t-\tau)/2}+e^{-(t-\tau)/2}-2]\chi_{[0,t]}(\tau)\varpi(\tau)\right|=0\,,\textit{for every fixed }\tau\,,$$

$$\sup\left\{\left|e^{-(t-\tau)/2}\frac{1}{t-\tau}[e^{(t-\tau)/2}+e^{-(t-\tau)/2}-2]\chi_{[0,t]}(\tau)\right| : t\in[0,\infty),\ \tau\in[0,\infty)\right\}<\infty$$

and, in view of the assumption (13.1)(1)( $j=0$ ), Lebesgue's dominated convergence theorem implies that

(13.8) $$\lim_{t\to\infty}\mathrm{I}_1(x,t)=0 .$$

Working as in the proof of (12.7) and using again the assumption (13.1)(1) ( $j=0$ ), we obtain

(13.9) $$\lim_{t\to\infty}\mathrm{I}_2(x,t)=0 .$$

Next, using the formula

$$e^{i\lambda x-\omega_1(\lambda)t}\int_0^\infty\chi_{[0,t]}(\tau)e^{\omega_1(\lambda)\tau}\varpi(\tau)d\tau=\frac{e^{i\lambda x}}{\omega_1(\lambda)}\varpi(t)-\frac{e^{i\lambda x-\omega_1(\lambda)t}}{\omega_1(\lambda)}\varpi(0)-\frac{e^{i\lambda x-\omega_1(\lambda)t}}{\omega_1(\lambda)}\int_0^\infty\chi_{[0,t]}(\tau)e^{\omega_1(\lambda)\tau}\varpi'(\tau)d\tau ,$$

we write

(13.10) $$\mathrm{I}_3=\int_{\mathbb{L}}\left\{e^{i\lambda x-\omega_1(\lambda)t}\int_0^\infty\chi_{[0,t]}(\tau)e^{\omega_1(\lambda)\tau}\varpi(\tau)d\tau\right\}\frac{\lambda d\lambda}{\rho(\lambda)}$$

$$=\left[\int_{\kappa^+}\frac{e^{i\lambda x}}{\omega_1(\lambda)}\frac{\lambda d\lambda}{\rho(\lambda)}\right]\varpi(t)-\left[\int_{-\kappa^-}\frac{e^{i\lambda x-\omega_1(\lambda)t}}{\omega_1(\lambda)}\frac{\lambda d\lambda}{\rho(\lambda)}\right]\varpi(0)-\int_{\mathbb{L}}\left\{\frac{e^{i\lambda x-\omega_1(\lambda)t}}{\omega_1(\lambda)}\int_0^\infty\chi_{[0,t]}(\tau)e^{\omega_1(\lambda)\tau}\varpi'(\tau)d\tau\right\}\frac{\lambda d\lambda}{\rho(\lambda)} .$$

Now we show that

(13.11) $$\lim_{t\to\infty}\left\{\left(\frac{e^{i\lambda x-\omega_1(\lambda)t}}{\omega_1(\lambda)}\int_0^\infty\chi_{[0,t]}(\tau)e^{\omega_1(\lambda)\tau}\varpi'(\tau)d\tau\right)\frac{\lambda}{\rho(\lambda)}\right\}=0 \text{, for every fixed } \lambda\in\mathbb{L} .$$

Indeed, this follows from Lebesgue's dominated convergence theorem and (13.1)(1)( $j=1$ ), since

$$\lim_{t\to\infty}\left(\frac{e^{i\lambda x-\omega_1(\lambda)t}}{\omega_1(\lambda)}\chi_{[0,t]}(\tau)e^{\omega_1(\lambda)\tau}\varpi'(\tau)\frac{\lambda}{\rho(\lambda)}\right)=0\ \ (\forall\tau\geq 0)$$

and

$$\left|\frac{e^{i\lambda x-\omega_1(\lambda)t}}{\omega_1(\lambda)}\chi_{[0,t]}(\tau)e^{\omega_1(\lambda)\tau}\varpi'(\tau)\frac{\lambda}{\rho(\lambda)}\right|\leq\left|\varpi'(\tau)\frac{\lambda}{\omega_1(\lambda)\rho(\lambda)}\right|\ \ (\forall t\geq 0,\forall\tau\geq 0).$$

Also,

$$\frac{e^{i\lambda x-\omega_1(\lambda)t}}{\omega_1(\lambda)}\int_0^\infty\chi_{[0,t]}(\tau)e^{\omega_1(\lambda)\tau}\varpi'(\tau)d\tau=\frac{e^{i\lambda x}}{[\omega_1(\lambda)]^2}\varpi'(t)-\frac{e^{i\lambda x-\omega_1(\lambda)t}}{[\omega_1(\lambda)]^2}\varpi'(0)-\frac{e^{i\lambda x-\omega_1(\lambda)t}}{[\omega_1(\lambda)]^2}\int_0^\infty\chi_{[0,t]}(\tau)e^{\omega_1(\lambda)\tau}\varpi''(\tau)d\tau$$

and, in view of (13.1)(3) and (13.1)(1)( $j=2$ ), we obtain that

(13.12) $$\left\{\frac{e^{i\lambda x-\omega_1(\lambda)t}}{\omega_1(\lambda)}\int_0^\infty\chi_{[0,t]}(\tau)e^{\omega_1(\lambda)\tau}\varpi'(\tau)d\tau\right\}\frac{\lambda}{\rho(\lambda)}=\mathrm{O}(1/\lambda^2)\text{, as } \lambda\to\pm\infty\text{, uniformly for } t\geq 0 .$$

It follows from (13.11) and (13.12) that the last integral in the RHS of (13.10) tends to zero, as $t\to\infty$. Since the integrals over $\kappa^+$ (in view of (13.1)(2)) and $\kappa^-$, in the RHS of (13.10), tend also to zero, we see that

(13.13) $$\lim_{t\to\infty}\mathrm{I}_3(x,t)=0 .$$

Finally,

(13.14) $$\mathrm{I}_4=\int_{\kappa^+}\left\{e^{i\lambda x-\omega_2(\lambda)t}\int_0^\infty\chi_{[0,t]}(\tau)e^{\omega_2(\lambda)\tau}\varpi(\tau)d\tau\right\}\frac{\lambda d\lambda}{\rho(\lambda)}=\iint_{(\lambda,\tau)\in\kappa^+\times[0,\infty)}\left\{e^{i\lambda x-\omega_2(\lambda)(t-\tau)}\chi_{[0,t]}(\tau)\varpi(\tau)\frac{\lambda}{\rho(\lambda)}\right\}d\lambda d\tau .$$

Since the integrand of the integral in RHS of (13.14), tends to zero, as $t\to\infty$, $\forall(\lambda,\tau)\in\kappa^+\times[0,\infty)$, and $\int_{\kappa^+\times[0,\infty)}|\varpi(\tau)|d|\lambda|d\tau<+\infty$, Lebesgue's dominated convergence theorem implies that

(13.15) $$\lim_{t\to\infty}\mathrm{I}_4(x,t)=0 .$$

Now, (13.6) follows from (13.7), (13.8), (13.9), (13.13) and (13.15).

***Step 2*** Setting

$$\mathcal{D}_0(f)(x,t) := \int_{-\infty}^{\infty}\left[e^{i\lambda x-\omega_2(\lambda)t}\int_0^t e^{\omega_2(\lambda)\tau}[\hat f(\lambda,\tau+T)-\hat f(\lambda,\tau)]d\tau - e^{i\lambda x-\omega_1(\lambda)t}\int_0^t e^{\omega_1(\lambda)\tau}[\hat f(\lambda,\tau+T)-\hat f(\lambda,\tau)]d\tau\right]\frac{d\lambda}{2i\rho(\lambda)}$$
$$+\int_{-\infty}^{\infty}\left[e^{i\lambda x-\omega_1(\lambda)t}\int_0^t e^{\omega_1(\lambda)\tau}[\hat f(-\lambda,\tau+T)-\hat f(-\lambda,\tau)]d\tau - e^{i\lambda x-\omega_2(\lambda)t}\int_0^t e^{\omega_2(\lambda)\tau}[\hat f(-\lambda,\tau+T)-\hat f(-\lambda,\tau)]d\tau\right]\frac{d\lambda}{2i\rho(\lambda)},$$

we claim that

(13.16) $\lim_{t\to\infty}\mathcal{D}_0(f)(x,t)=0$, uniformly for $x$ in compact subsets of $[0,\infty)$.

For its proof, we split the integrals as follows:

$$(13.17)\quad \mathcal{D}_0(f)(x,t)=\int_{\lambda=-1/2}^{1/2}\left\{(e^{i\lambda x}-e^{-i\lambda x})\left[e^{-\omega_2(\lambda)t}\int_0^t e^{\omega_2(\lambda)\tau}[\hat f(\lambda,\tau+T)-\hat f(\lambda,\tau)]d\tau\right.\right.$$
$$\left.\left.-e^{-\omega_1(\lambda)t}\int_0^t e^{\omega_1(\lambda)\tau}[\hat f(\lambda,\tau+T)-\hat f(\lambda,\tau)]d\tau\right]\frac{d\lambda}{2i\rho(\lambda)}\right\}$$
$$+\left(\int_{-\infty}^{-1/2}+\int_{1/2}^{\infty}\right)\left[e^{i\lambda x-\omega_2(\lambda)t}\int_0^t e^{\omega_2(\lambda)\tau}[\hat f(\lambda,\tau+T)-\hat f(\lambda,\tau)]d\tau - e^{i\lambda x-\omega_1(\lambda)t}\int_0^t e^{\omega_1(\lambda)\tau}[\hat f(\lambda,\tau+T)-\hat f(\lambda,\tau)]d\tau\right]\frac{d\lambda}{2i\rho(\lambda)}$$
$$+\left(\int_{-\infty}^{-1/2}+\int_{1/2}^{\infty}\right)\left\{\left[e^{i\lambda x-\omega_1(\lambda)t}\int_0^t e^{\omega_1(\lambda)\tau}[\hat f(-\lambda,\tau+T)-\hat f(-\lambda,\tau)]d\tau\right.\right.$$
$$\left.\left.-e^{i\lambda x-\omega_2(\lambda)t}\int_0^t e^{\omega_2(\lambda)\tau}[\hat f(-\lambda,\tau+T)-\hat f(-\lambda,\tau)]d\tau\right]\frac{d\lambda}{2i\rho(\lambda)}\right\}.$$

Now the integral in (13.17) taken over $[-\frac{1}{2},\frac{1}{2}]$ is equal to

$$(13.18)\quad 4\int_{\lambda=-1/2}^{1/2}\sin(\lambda x)\left[\int_{\tau=0}^{t}\left\{e^{-(t-\tau)/2}\frac{1}{t-\tau}\sum_{n=1}^{\infty}\frac{(-1)^{n-1}}{(2n-1)!}\left(\lambda^2-\frac{1}{4}\right)^{n-1}(t-\tau)^{2n}\int_0^{\infty}e^{-i\lambda y}[f(y,\tau+T)-f(y,\tau)]dy\right\}d\tau\right]d\lambda .$$

The absolute value of the above integral is

$$\le 4\int_{\lambda=-1/2}^{1/2}\left|\sin(\lambda x)\right|\left[\int_{\tau=0}^{\infty}\left\{\chi_{[0,t]}(\tau)e^{-(t-\tau)/2}\frac{1}{t-\tau}\sum_{n=1}^{\infty}\frac{1}{(2n-1)!}\left(\frac{1}{4}-\lambda^2\right)^{n-1}(t-\tau)^{2n}\int_0^{\infty}\left|f(y,\tau+T)-f(y,\tau)\right|dy\right\}d\tau\right]d\lambda .$$

Working as in previous cases (in particular, replacing $\left|\sin(\lambda x)\right|$ by $\left|\lambda\right|x$ and integrating, first, we respect to $\lambda$), we see, using the assumption (13.2)(1) and Lebesgue's dominated convergence theorem, that the limit, as $t\to\infty$, of the integral (13.18) is zero.

For the 2nd integral in RHS of (13.17), we write it as follows:

$$(13.19)\quad\left(\int_{-\infty}^{-1/2}+\int_{1/2}^{\infty}\right)\left[\left\{e^{i\lambda x-\omega_2(\lambda)t}\int_0^{\infty}\chi_{[0,t]}(\tau)e^{\omega_2(\lambda)\tau}\left[\frac{f(0,\tau+T)-f(0,\tau)}{i\lambda}\right]d\tau\right.\right.$$
$$\left.\left.-e^{i\lambda x-\omega_1(\lambda)t}\int_0^{\infty}\chi_{[0,t]}(\tau)e^{\omega_1(\lambda)\tau}\left[\frac{f(0,\tau+T)-f(0,\tau)}{i\lambda}\right]d\tau\right\}\frac{d\lambda}{2i\rho(\lambda)}\right]$$
$$+\left(\int_{-\infty}^{-1/2}+\int_{1/2}^{\infty}\right)\left[\left\{e^{i\lambda x-\omega_2(\lambda)t}\int_0^{t}e^{\omega_2(\lambda)\tau}\left[\frac{(\hat{f_x})(\lambda,\tau+T)-(\hat{f_x})(\lambda,\tau)}{i\lambda}\right]d\tau\right.\right.$$
$$\left.\left.-e^{i\lambda x-\omega_1(\lambda)t}\int_0^{t}e^{\omega_1(\lambda)\tau}\left[\frac{(\hat{f_x})(\lambda,\tau+T)-(\hat{f_x})(\lambda,\tau)}{i\lambda}\right]d\tau\right\}\frac{d\lambda}{2i\rho(\lambda)}\right].$$

The absolute value of the 1st of the above integrals is

$$\leq\left(\int_{-1}^{-1/2}+\int_{1/2}^{1}\right)\left[\left\{\int_{\tau=0}^{\infty}\chi_{[0,t]}(\tau)e^{-(t-\tau)/2}\left|\frac{\sin[\rho(\lambda)(t-\tau)]}{\rho(\lambda)}\right|\left|f(0,\tau+T)-f(0,\tau)\right|d\tau\right\}\frac{d\lambda}{|\lambda|}\right]$$

$$+2\left(\int_{-\infty}^{-1}+\int_{1}^{\infty}\right)\left\{\left[\int_{\tau=0}^{\infty}\chi_{[0,t]}(\tau)e^{-(t-\tau)/2}\left|f(0,\tau+T)-f(0,\tau)\right|d\tau\right]\frac{d\lambda}{|\lambda\rho(\lambda)|}\right\}$$

$$=\iint_{(\tau,\lambda)\in[0,\infty)\times\{[-1,-\frac{1}{2}]\cup[\frac{1}{2},1]\}}\chi_{[0,t]}(\tau)e^{-(t-\tau)/2}|t-\tau|\left|\frac{\sin[\rho(\lambda)(t-\tau)]}{\rho(\lambda)(t-\tau)}\right|\left|f(0,\tau+T)-f(0,\tau)\right|\frac{d\tau d\lambda}{|\lambda|}$$

$$+2\iint_{(\tau,\lambda)\in[0,\infty)\times\{(-\infty,-1]\cup[1,\infty)\}}\chi_{[0,t]}(\tau)e^{-(t-\tau)/2}\left|f(0,\tau+T)-f(0,\tau)\right|\frac{d\tau d\lambda}{|\lambda\rho(\lambda)|}$$

and, therefore, this integral tends to zero, as $t\to\infty$, by (13.2)(3).

The absolute value of the 2nd integral in (13.19) is

$$\leq\left(\int_{-1}^{-1/2}+\int_{1/2}^{1}\right)\left[\left\{\int_{\tau=0}^{\infty}\chi_{[0,t]}(\tau)e^{-(t-\tau)/2}\left|\frac{\sin[\rho(\lambda)(t-\tau)]}{\rho(\lambda)}\right|\left(\int_{y=0}^{\infty}\left|f_x(y,\tau+T)-f_x(y,\tau)\right|dy\right)d\tau\right\}\frac{d\lambda}{|\lambda|}\right]$$

$$+2\left(\int_{-\infty}^{-1}+\int_{1}^{\infty}\right)\left\{\int_{\tau=0}^{\infty}\chi_{[0,t]}(\tau)e^{-(t-\tau)/2}\left(\int_{y=0}^{\infty}\left|f_x(y,\tau+T)-f_x(y,\tau)\right|dy\right)d\tau\frac{d\lambda}{|\lambda\rho(\lambda)|}\right\}$$

$$\leq\iiint_{(\tau,y,\lambda)\in[0,\infty)\times[0,\infty)\times\{[-1,-\frac{1}{2}]\cup[\frac{1}{2},1]\}}\chi_{[0,t]}(\tau)e^{-(t-\tau)/2}|t-\tau|\left|\frac{\sin[\rho(\lambda)(t-\tau)]}{\rho(\lambda)(t-\tau)}\right|\left|f_x(y,\tau+T)-f_x(y,\tau)\right|\frac{d\tau dy d\lambda}{|\lambda|}$$

$$+2\iiint_{(\tau,y,\lambda)\in[0,\infty)\times[0,\infty)\times\{(-\infty,-1]\cup[1,\infty)\}}\chi_{[0,t]}(\tau)e^{-(t-\tau)/2}\left|f_x(y,\tau+T)-f_x(y,\tau)\right|\frac{d\tau dy d\lambda}{|\lambda\rho(\lambda)|}$$

and, therefore, this integral tends to zero, as $t\to\infty$, by (13.2)(2).

Thus, the sum of integrals (13.19) tends to zero, as $t\to\infty$. Similarly we show that the 3rd integral in (13.17) also tends to zero, as $t\to\infty$. Finally, these conclusions imply (13.16).

*Completion of the proof* Combining the remarks made at the begining of the proof with (13.7) and (13.16), we obtain (13.3).

**Declarations**

*Ethical Approval*: Not applicable.
*Competing Interests*: Not applicable.
*Availability of data and materials*: Not applicable.